\documentclass[12pt]{amsart}
\usepackage{amsmath}
\usepackage{amssymb}
\usepackage{upgreek}
\usepackage{url}
\usepackage{amsthm}
\usepackage{enumerate}
\usepackage[all]{xy}
\usepackage[margin=.5in]{geometry}
\usepackage{verbatim}
\usepackage{setspace}
\usepackage{etoolbox}
\usepackage{environ}
\usepackage[utf8]{inputenc}
\usepackage{bm}
\usepackage{thmtools}

\usepackage{dirtytalk}

\usepackage{nameref,hyperref,cleveref}

\hypersetup{
    pdftitle={Proj. Duality, Unexpected Hypersurfaces and Hyperplane Deriv.},
    pdfauthor={Bill Trok},
    pdfsubject={Algebraic Geometry},
    pdfkeywords={Hyperplanes, Points, Algebra, Geometry},
    bookmarksnumbered=true,     
    bookmarksopen=true,         
    bookmarksopenlevel=1,       
    colorlinks=true,            
    pdfstartview=Fit,           
    pdfpagemode=UseOutlines,    
    pdfpagelayout=TwoPageRight
}

\newcommand{\bb}[1]{\mathbb{#1}}
\newcommand{\cal}[1]{\mathcal{#1}}
\newcommand{\cl}[1]{\mathcal{#1}}
\newcommand{\by}{\times}
\newcommand{\tensor}{\otimes}
\newcommand{\chr}{\mathrm{char}}

\newcommand{\proj}[1]{\mathrm{Proj}\left(#1\right)}
\newcommand{\partialfrac}[2]{\frac{\partial #1}{\partial #2}}
\newcommand{\angvec}[1]{\left\langle #1 \right\rangle}

\newcommand{\iso}{\cong}
\newcommand{\cross}{\times}

\DeclareMathOperator{\Der}{Der}
\DeclareMathOperator{\Spec}{Spec}
\DeclareMathOperator{\Proj}{Proj}

\DeclareMathOperator{\Sym}{Sym}
\DeclareMathOperator{\Gr}{Gr}

\DeclareMathOperator{\res}{res}
\DeclareMathOperator{\codim}{codim}
\DeclareMathOperator{\reg}{reg}
\DeclareMathOperator{\Pl}{PL}
\DeclareMathOperator{\Spn}{Span}
\DeclareMathOperator{\Bloc}{B.loc}
\DeclareMathOperator{\Cl}{Cl}
\DeclareMathOperator{\rk}{rk}

\DeclareMathOperator{\rank}{rank}
\DeclareMathOperator{\pow}{Pow}
\DeclareMathOperator{\Exc}{Ex.C}
\DeclareMathOperator{\Exbl}{Ex.Bl}

\declaretheorem[name=Theorem,parent = section,
refname={theorem,theorems}]{theorem}
\declaretheorem[name=Lemma,sibling = theorem,
refname={lemma,lemmas}]{lemma}
\declaretheorem[name=Propostion,sibling = theorem,
refname={proposition,propositions}]{prop}
\declaretheorem[name=Remark,sibling = theorem,
refname={remark,remarks},style=remark]{remark}
\declaretheorem[name=Claim,sibling = theorem,
refname={claim,claims},style=remark]{claim}
\declaretheorem[name=Definition,sibling = theorem,
refname={definition,definitions},style=definition]{defn}
\declaretheorem[name=Corollary,sibling = theorem,
refname={corollary,corollaries}]{cor}

\declaretheorem[name=Example,sibling = theorem,
refname={example,examples}]{example}
\declaretheorem[name=Conjecture,sibling = theorem,
refname={conjecture,conjectures}]{conj}

\numberwithin{equation}{theorem}

\newtoggle{outline}
\togglefalse{outline}

\NewEnviron{otProof}{\iftoggle{outline}{}{\emph{Proof.} \BODY
    \hfill\(\square\)}}

\newtoggle{PrintForEdit}

\togglefalse{PrintForEdit}

\iftoggle{PrintForEdit}{ \singlespacing }

\newcommand{\PFEIndent}
{ \iftoggle{PrintForEdit}{ \pagebreak }
  {
  }
}
\begin{document}

\title[Proj. Duality, Unexpected Hypersurfaces and Hyperplane Deriv.]{Projective Duality, Unexpected Hypersurfaces and Logarithmic
  Derivations of Hyperplane Arrangements}

\date{\today}
\author{Bill Trok}
\address{Department of Mathematics, University of Kentucky, Lexington,
  Kentucky, 40508}
\curraddr{Department of Mathematics, University of Kentucky, Lexington,
  Kentucky, 40508}

\email{william.trok@uky.edu}

\begin{abstract} Several papers have been written studying
  unexpected hypersurfaces. We say a finite set of points \(Z\)
  admits unexpected hypersurfaces if a general union of
  fat linear subspaces imposes less that the expected number of conditions on the
  ideal of \(Z\). In this paper, we introduce the concept of
  a very unexpected hypersurface. This is a stronger condition
  which takes into account an explanation for some hypersurfaces
  previously considered unexpected.

  We then develop a duality theory to relate the study of very unexpected
  hypersurfaces to the derivations of
  dual hyperplane arrangements. This allows us to generalize results
  in the plane of Cook, Harbourne, Migliore, Nagel \cite{CHMN} Faenzi, and
  Vall\'{e}s \cite{FV} to higher dimensions. In
  particular, we give a criterion to determine if a set \(Z\) admits
  very unexpected hypersurfaces in the case \(Z\) is invariant
  under the action of an irreducible reflection group on the ambient
  projective space.

  Our approach has new applications even in
 \(\bb{P}^2\) where we are able to place strong conditions on sets of points \(Z\)
 which admit certain types of unexpected curves. We close relating
 Terao's Freeness Conjecture for line arrangements to a conjecture due to G. Dirac \cite{Dirac}
 on configurations of points in \(\bb{R}^2\). 
\end{abstract}

\maketitle
\section{Introduction}\label{sec:Introduction}

Determining the dimension of a linear system is a fundamental problem
in algebraic geometry, and many classic results
are devoted to determining the dimensions of certain linear
systems. Frequently, when working with families of linear systems,
there is a naive dimension count, or an \emph{expected dimension}, that
should hold in the general case. The difficulty then lies in
determining when this dimension count fails.

As an illustration, we consider the projective
coordinate ring \(R = \bb{C}[X_0,X_1,X_2,X_3]\) of
\(\bb{P}^3_{\bb{C}}\). For a proper linear subspace \(L \subset \bb{P}^3\)
vanishing on \(L\) imposes \(\binom{\dim L + d}{d}\) conditions on forms of degree
\(d\). This gives a dimension count for the ideal \(I(L)\) of \(\dim [I(L)]_d =
\dim[R_d]-\binom{\dim L +d}{d} = \binom{d+3}{d}-\binom{\dim L +
  d}{d}.\) 

By extension, if \(Z \subseteq \bb{P}^3\) is any proper subscheme we
might expect that \[\dim[I(Z) \cap I(L)]_d = \max\left\{0,\dim[I(Z)]_d -
    \binom{\dim L + d}{d}\right\}.\]
In fact the above equality holds for all
sufficiently large \(d\) as long as \(L\) and \(Z\) have empty
intersection. Let \(Z\) be a \(0\)-dimensional scheme and let \(L\) be a
general linear subspace, to ensure that \(L \cap
Z=\emptyset\). We then get the expected dimension:

\begin{equation}\label{eq:ExcDim}
  \tag{1}[I(Z) \cap I(L)^m]_d = \min\{0, \dim[I(Z)]_d - \dim [R/I(L)^m]_d\}.
\end{equation}

Many papers have been written exploring when the above
formula fails to give the actual dimension. If \(Z\) is a set of double points with general support and \(L\) is a general point with
\(m=2\), the celebrated theorem of Alexander and Hirschowitz
\cite{AH} gives a
complete characterization of when equation \ref{eq:ExcDim} fails to
give the correct count. Where \(m=1\) and \(L \subseteq
\bb{P}^3\) is a line, the paper \cite{Migl08} relates equation
\ref{eq:ExcDim} to Lefschetz properties of general Artinian reductions 
of \(R/I(Z)\).

More recently, a number of papers have been released (see for instance
\cite{HMT},\cite{DFHMST},\cite{BMSS} and \cite{HMUZ}), which studied failure
of expected dimension under the name unexpected hypersurfaces
or unexpected curves. These papers all
study the above linear systems, where \(Z\) is a reduced sets of
points in \(\bb{P}^n\), and \(L\) is (possibly multiple) general
linear subspace with an associated multiplicity. Many of these papers
took inspiration from the paper \cite{CHMN}.
The authors of \cite{CHMN} built off of earlier work in \cite{DIV} and introduced the concept of unexpected curves in \(\bb{P}^2\).
Namely, they said that \(Z\) admits
unexpected curves in degree \(d\) if for a general point \(X\),
\[ \dim [I(Z) \cap I(X)^{d-1}]_d > \max\left\{0,\dim
  [I(Z)]_d - \dim [R/I(X)^{d-1}]_d\right\}.\]

The authors of \cite{CHMN} were able to give a full characterization of the degrees in which
a set of points admits unexpected curves. Surprisingly, this characterization does
not directly depend on the dimensions of either \([I(Z)]_d\) or \([I(Z) \cap I(X)^{d-1}]_d\). Namely, this
information can be replaced with combinatorial information about
\(Z\), and data coming from the (reduced) module of derivations,
\(D_0(\cl{A}_Z)\), of the line arrangement, \(\cl{A}_Z\),
dual to \(Z\). See \cref{defn:SplittingType} for the definition of
splitting type.

\begin{restatable*}[\cite{CHMN}]{theorem}{CHMNthm}\label{thm:UnexpectedCurveChar}
  For a finite set of points  \(Z \subseteq
  \bb{P}^2\), let \(\cl{A}_Z\) denote the dual
  line arrangement, and let \((a_1,a_2)\) denote the splitting type
  of the bundle defined by \(D_0(\cl{A}_Z)\). Then exactly one of the following
  statements holds:
  \begin{enumerate}[(i)]
  \item There is some line \(L \subseteq
    \bb{P}^2\) with \(|L \cap Z| > a_1+1\), in which case \(|L \cap
    Z| = a_2+1\) and \(Z\) never admits unexpected curves.
  \item \(Z\) admits unexpected curves in degree \(d\) for precisely
  those \(d\) with
  \(a_1 < d < a_2\).
\end{enumerate}
\end{restatable*} 

This result allowed researchers to discover many new examples of
unexpected curves by taking advantage of decades of prior research on
line arrangements.

Given the observed connection between certain line arrangements and unexpected curves,
it is natural to wonder if a similar connection exists in higher
dimensions. In this paper we show that this is true at least to a
certain extent. More
specifically, if \(Z \subseteq \bb{P}^n\) is a finite set of points
and \(L\) is a general codimension 
\(2\) linear subspace, we establish a general duality connecting the module of
derivations \(D_0(\cl{A}_Z)\) of the dual hyperplane 
arrangement to the intersection of ideals \([I(Z) \cap
I(L)^{d}]_{d+1}\). In particular this allows us to recover \(\dim
[I(Z) \cap I(L)^{d}]_{d+1}\) from knowledge of the splitting type of \(D_0(\cl{A}_Z)\).

In order to generalize \ref{thm:UnexpectedCurveChar}, we introduce a
modified definition of unexpected hypersurface which we call
very unexpected hypersurfaces. Given a generic linear subspace \(L\), we say a finite set of points \(Z
\subseteq \bb{P}^n\) admits very unexpected \(L\)-hypersurfaces
if the intersection \([I(Z) \cap I(L)^{d-1}]_d\)  is larger than
expected, as long as this failure is not ``easily
explained'' (see \cref{defn:VeryUnexpectedHypersurface}). Our definition
of very unexpected hypersurfaces is more technical than
that of unexpected hypersurfaces. However
the two definitions agree in \(\bb{P}^2\), in that a set of points
 \(Z \subseteq \bb{P}^2\) admits very unexpected curves if and only if
 it admits unexpected curves.

This new definition has a few
advantages compared with the definition for unexpected
  hypersurfaces. The first is that with the standard definition of
unexpected hypersurfaces a
generalization of \cref{thm:UnexpectedCurveChar} to higher dimensions
is impossible. The second is that, as we mentioned, in certain cases
the ``unexpectedness'' can be relatively easily explained. 
For instance, if all of the points in \(Z\) lie on a proper
subspace \(H\), ``unexpectedness'' may simply be a consequence of the fact
that \([I(H) \cap
I(L)^{d-1}]_d \subseteq [I(Z) \cap I(L)^{d-1}]_d\)  (for further
discussion, see \cref{example:TotallyExpected}).
It is then somewhat surprising that by merely accounting for cases where
``unexpectedness'' is well explained, we are able to recover a
generalization of \cref{thm:UnexpectedCurveChar}. More specifically, if
\(L\) is a generic codimension \(2\) subspace, the
degrees in which \(Z\) admits very unexpected \(L\)-hypersurfaces can again be
characterized by the combinatorial data of \(Z\) in conjunction with the
splitting type of the Derivation Bundle of \(\cl{A}_Z\). We define, for
every integer \(d \geq 0\), a number \(\Exc(Z,d)\) 
via a combinatorial optimization problem on the matroid of
\(Z\) (see \cref{defn:BestDefUHyp}). Using this number, \(\Exc(Z,d)\), we obtain the result below.

\begin{restatable*}{theorem}{BestDefUHyp}\label{cor:BestDefnUHyp} Let \(Z \subseteq \bb{P}^n\) be
  a finite set of points, and suppose that \(D_0(\cl{A}_Z)\) has splitting
  type \((a_1,..,a_n)\). Then for a fixed integer \(d\), 
  \[\begin{aligned}
      &\sum_{i=1}^n\max\{0,d-a_i\} \leq& nd+1 - \Exc(Z,d) \\
    \end{aligned}  \]
  and the inequality is strict if and only if \(Z\) admits very unexpected
  hypersurfaces in degree \(d\).
\end{restatable*}

In the case the points of \(Z\) are not too concentrated on some proper
subspace we obtain the following result which mimics \cref{thm:UnexpectedCurveChar}. 

\begin{restatable*}{theorem}{CHMNGeneralization}\label{thm:CHMNGeneralization}
  Let \(Z \subseteq \bb{P}^n\) and let \((a_1,..,a_n)\) be the splitting type of
  \(D_0(\cl{A}_Z)\), where \(a_i \leq a_{i+1}\).
  Suppose for all positive dimensional linear subspaces \(H \subseteq
  \bb{P}^n\), we have that
  \[ \frac{|Z \cap H| -1}{\dim H} \leq \frac{|Z| - 1}{n}.\]
  Then for an integer \(d\) the following are equivalent:
  \begin{enumerate}[(a)]
  \item \(Z\) admits very unexpected hypersurfaces in degree \(d\).
  \item \(Z\) admits unexpected hypersurfaces in degree \(d\).
  \item \(a_1 < d < a_n\).
\end{enumerate}
\end{restatable*}

Moreover we show in \cref{prop:ReflectIsGood} that this condition holds if an irreducible
reflection group \(G \subseteq \bb{PGL}(\bb{K},2)\) which acts on
\(Z\).

After discussing the theory of Unexpected
Hypersurfaces in general.
We apply this duality between \(I(Z)\) and \(D_0(\cl{A}_Z)\) to establish some structural results about both very unexpected \(Q\)-hypersurfaces
and the module of derivations \(D_0(\cl{A}_Z)\). Unlike the first part
of the paper where there are few dimension and field constraints,
these results focus on unexpected curves
in \(\bb{P}^2_{\bb{C}}\). In particular, we establish the following
bound sharp for all \(d \geq 1\).

\begin{restatable*}{theorem}{PointCardBound}\label{thm:PointCardBound}
  Let \(Z \subseteq \bb{P}^2_{\bb{C}}\) and suppose that \(|Z|\)
  admits an unexpected curve in degree \(d \geq 1\), then \(|Z| \leq 3d-3\).
\end{restatable*}

We note that if \(d\) is the smallest such degree in which \(Z\)
admits unexpected curves then it follows from
\cref{thm:UnexpectedCurveChar} that \(2d+1 \leq |Z|\). Consequently, no \(Z \subseteq \bb{P}^2_{\bb{C}}\) admits unexpected
curves in degree \(3\) or lower. 

Additionally, we show the splitting type of
\(D_0(\cl{A}_Z)\) can be easily determined using only the initial
degree of \(D_0(\cl{A}_Z)\). We note that the splitting type is
determined by the initial degree of the restriction of
\(D_0(\cl{A}_z)\) to a general line.

\begin{restatable*}{theorem}{SplittingTypeChar}\label{thm:SplittingTypeChar}
  Let \(Z\) be a finite set of points in \(\bb{P}^2_{\bb{C}}\) and
  let \(\alpha(D_0(\cl{A}_Z))\) denote the initial degree of
  \(D_0(\cl{A}_Z)\). Define \(a =
  \min\left\{\alpha(D_0(\cl{A}_Z)),
    \left\lfloor\frac{|Z|-1}{2}\right\rfloor\right\}\) 
  then \(D_0(\cl{A}_Z)\) has splitting type \((a,|Z|-a-1)\).
\end{restatable*}

The paper proceeds as follows. After defining some notation in
\cref{sec:Notation}. We  discuss some needed background
on the module of logarithmic derivations of a hyperplane
arrangement in \cref{sec:Derivations}. The reader familiar with Hyperplane Arrangements can
likely skip this section with perhaps the exception of
some non-standard notation found in \cref{defn:StanleysRho} and
\cref{defn:StanleysRestrictedRho}.

We proceed in \cref{sec:DerivationIdealCorrespondence},
expanding on the Faenzi-Vall\'{e}s duality between the module of derivations
\(D_0(\cl{A}_Z)\) of a hyperplane arrangement and certain elements of
the ideal \(I(Z)\) of points dual to \(\cl{A}_Z\). The results of this
section are
not wholly original as much of this is implicit in the first section of
\cite{FV}. Our approach however, is much more explicit and amenable to computation it also has
the advantage of working in arbitrary characteristic.  We state two
versions of this correspondence, the first (\cref{thm:GlobalIso})
applies to the module, \(D_0(\cl{A}_Z)\) itself, and we do not
believe it has been stated before in this form.
The second correspondence (\cref{thm:restrictedModIso}) applies to the
restriction of \(D_0(\cl{A}_Z)\) to a general line, generalizes the
duality found in \cite{FV}. We note that despite the
similarities in results, our method of proof and presentation is quite
different from the one given in \cite{FV}. Additionally, the results
here are not dependent on the characteristic of the ground field
\(\bb{K}\).

Section \ref{sec:UnexpectedHypersurfaces} introduces our definition of
an very unexpected \(Q\)-hypersurface (see
\cref{defn:VeryUnexpectedHypersurface}) for \(Q\)
a generic subspace. We then look at the case when \(Q\) has
codimension \(2\), establishing in
\cref{thm:BestDefUHyp} that the degrees in which \(Z\) admits
very unexpected \(Q\)-hypersurfaces
depends only on the splitting type of \(D_0(\cl{A}_Z)\) and
a combinatorial optimization problem involving \(Z\).

Section \ref{sec:LiftingAndCurves} starts by establishing a lifting
criterion for the restriction of
\(D_0(\cl{A}_Z)\) to a general line (see \cref{prop:DivisibilityProp}). We then recall some results on vector
bundles on \(\bb{P}^2_{\bb{C}}\) and apply these to show that
\cref{prop:DivisibilityProp} has
especially strong consequences in \(\bb{P}^2_{\bb{C}}\) (see
\cref{thm:UnexpectednessIsGlobal} and
\cref{cor:PolynomialStructure}). 

In Section \ref{sec:CombinatoricsInP2} we
give strong combinatorial constraints on the sets of points \(Z
\subseteq \bb{P}^2_{\bb{C}}\) which 
can admit unexpected curves. In particular, if \(Z \subseteq \bb{P}^2_{\bb{C}}\) admits an
unexpected curve in degree \(d\), then
\cref{thm:LinearlySpecial} shows that no more than \(d+1\) points
of \(Z\) are in linearly general position and
\cref{thm:PointCardBound} establishes a sharp bound on the number of
points in \(Z\) showing that \(|Z| \leq 3d-3\).

In section \ref{sec:Regularity} we show that if \(Z \subseteq
\bb{P}^2_{\bb{C}}\) admits unexpected curves in degree \(d\), then
\(Z\) imposes independent condition on
\((d-1)\) forms. We then briefly discuss generalizations to higher
dimensions and some consequences.

We close with section \ref{sec:ApplicationsToTeraos},
which discusses a few applications of these results to the field of
Hyperplane arrangements. We focus on Terao's Freeness Conjecture
mainly in \(\bb{P}^2_{\bb{C}}\). In particular, we look at the
conjecture for real line arrangements and connect it to the Weak Dirac
Conjecture on real point configurations.

We have attempted to keep this paper as self contained and elementary
as possible. This is largely true for the first 5 sections. However,
in later sections we do apply some results from the theory of Vector 
Bundles and from the combinatorics of line arrangements in
\(\bb{P}^n_{\bb{C}}\).

\begin{center}
  {\large Acknowledgements}
\end{center}

I would like to thank Brian Harbourne, Juan Migliore, Samuel
Herschenfeld, Uwe Nagel and Alexandra Seceleanu for the helpful
discussions. I'd also like to thank  Uwe Nagel and Alexandra Seceleanu
for the encouragement and advice that helped with the exposition.

\section{Notation and Conventions}\label{sec:Notation}
Throughout this paper \(\bb{K}\) will denote an algebraically closed of arbitrary
characteristic, unless specified otherwise. However, most of these results hold
as long as \(\bb{K}\) is infinite. \(V\) and \(W\) will be
dual \(\bb{K}\)-vector spaces. 
That is we suppose that there is a non-degenerate bilinear
pairing \(B(\; , \,): V \times W \to \bb{K}\), inducing isomorphisms
\(V \iso W^{*}\) and \(W \iso V^{*}\) 

If \(V\)  is a \(\bb{K}\) vector space, then \(V^{*}\)
will denote the dual vector space of linear maps \(V \to \bb{K}\). Our
pairing gives isomorphisms \(V \iso W^{*}\) and \(W \iso V^{*}\), we
denote these isomorphisms \(v \mapsto \ell_v\) and \(w \mapsto
\ell_w\), respectively. Here \(\ell_v(w) = \ell_w(v) = B(v,w)\). If
\(H \subseteq V\) is a linear subspace, then \(H^{\perp} = \{ w \in W
 \mid \ell_w(H) = \{0\}\}\). We similarly define  \(L^{\perp}
 \subseteq V\), for \(L \subseteq W\).

\(\Sym(V^{*})\) will denote the graded \(\bb{K}\)-algebra of symmetric
tensors. Given a choice of basis \(\{Y_0,Y_1,..,Y_n\}\) of \(V^{*}\), \(\Sym(V^{*})\),
is naturally isomorphic to the polynomial algebra
\(\bb{K}[Y_0,..,Y_n]\).

Moreover, the graded ring \(R = \Sym(V^*)\) is naturally identifiable with the projective
coordinate ring of \(\bb{P}(V)\). Dually, \(S = \Sym(W^{*})\) is the
projective coordinate ring of \(\bb{P}(W)\). 

The goal of this paper, is to relate properties of a finite set of
points in \(\bb{P}(V)\) to their dual hyperplanes in
\(\bb{P}(W)\).

\PFEIndent

\section{Derivations of Hyperplane Arrangements}\label{sec:Derivations}

In this section we recall some facts about the module of logarithmic
derivations \(D(\cl{A})\) of a hyperplane arrangement \(\cl{A}\). In
particular we state a few different known criteria for a general \(S\)
derivation to lie in \(D(\cl{A})\). We also give the definition (\cref{defn:SplittingType})
of the splitting type of \(D(\cl{A})\) which is used heavily in the sequel.

\begin{defn}
  A (central) \emph{Subspace Arrangement}, \(\cl{A}\), is a
  finite collection of linear subspaces \(\{H_0,..,H_s\}\) of a vector
  space \(W\).

  If each \(H_i\) is a hyperplane, we say that
  \(\cl{A}\) is a \emph{Hyperplane Arrangement}. We say \(\cl{A}\) is
  \emph{essential} if the only subspace contained in all the
  hyperplanes in \(\cl{A}\) is the \(0\)-subspace.
\end{defn}

\begin{remark} All subspace arrangements in this paper will be
  central. We make this restriction in order to identify a subspace arrangement
  \(\cl{A}\) in \(W\) with it's image in \(\bb{P}(W)\),
  something we will do freely and often without comment.

  A hyperplane arrangement is often defined in terms of a defining
  polynomial \(Q_\cl{A} = \prod_{H \in \cl{A}} \ell_H\). This is the
  product of linear forms each one defining a unique hyperplane in
  \(\cl{A}\).
\end{remark}

\begin{defn} If \(S\) is our graded polynomial ring and \(M\) is a graded \(S\)-module, then 
  a \(\bb{K}\)-derivation of \(S\) into \(M\) is a graded
  \(\bb{K}\)-linear map \(\theta:S \to M\) which satisfies the Leibniz
  product rule. Namely for \(f,g \in S\)
  \[ \theta(f \cdot g) = \theta(f) \cdot g + f \cdot \theta(g)\].

  These form a graded \(S\)-module, denoted \(\Der(S,M)\), obtained by setting \((f \cdot \theta) (g) =
  f (\theta(g))\).

  We grade \(\Der(S,M)\) by the \emph{polynomial degree}, namely, we set
  \(\deg \theta = \deg(\theta(\ell))\), where \(\ell \in [S]_1\).
\end{defn}

  In the case that \(M = S\), we set \(\Der(S) := \Der(S,S)\). In this
  paper our module \(M\) will either be \(S\) or a quotient ring of
  \(S\).

\begin{defn} If \(\cl{A} \subseteq \bb{P}(W)\) is a Hyperplane
Arrangement, we define the \emph{module of \(\cl{A}\)-derivations}, denoted
\(D(\cl{A})\), as submodule of \(\Der(S)\) via
\[ D(\cl{A}) := \{ \theta \in \Der(S) \mid \theta(I(H)) \subseteq I(H)\text{ for all } H
  \in \cl{A}\}\] 
\end{defn}

\begin{remark}Each element \(\alpha \in W\) defines a
  \(\bb{K}\)-derivation, \(\theta_\alpha\), of \(S = \Sym(W^{*})\). Namely, for \(\ell \in
  [S]_1\) we set \(\theta_\alpha(\ell) = \ell(\alpha)\) and extended
  to all of \(S\) via the Leibniz product rule.
\end{remark}

\begin{prop} Let \(S = \Sym(W^{*})\) and let \(M\) be a graded
  \(S\)-module, then there's an isomorphism of graded \(S\)-modules \(\Der(S,M) \iso M
  \tensor_{\bb{K}} W\). Here the grading
  on \(M \tensor_{\bb{K}} W\) is given by that of \(M\).

  Consequently, theres an isomorphism \(\Der(S) \tensor_S M \iso \Der(S,M)\).
\end{prop}

\begin{proof} Picking a basis \(Y_0,..,Y_n\) for \(W^{*}\), we have
  \(S \iso \bb{K}[Y_0,..,Y_n]\). Let \(\theta \in \Der(S,M)\), let \(g_i
  = \theta(Y_i)\), by linearity and the Leibniz product rule we get
  these \(g_i\) completely determine \(\theta\). It follows that
  \(\theta\) is equal to the derivation \(\sum_{i=0}^n g_i
  \partialfrac{}{Y_i}\).

  Hence if \(W_0,..,W_n\) is a basis of \(W\) dual to 
  \(Y_0,..,Y_n\), meaning \(Y_i(W_j) = \delta_{i,j}\). Then
  \(\theta = \sum_{i=0}^n g_i
  \tensor W_i \in M \tensor_{\bb{K}} W\) establishing the first
  result.

  The second statement follows from the isomorphisms
  \[M \tensor_S \Der(S) \iso M \tensor_S (S \tensor_{\bb{K}} W) \iso M \tensor_{\bb{K}} W\]
\end{proof}

\begin{defn}\label{defn:StanleysRho} Let \(S = \Sym(W^{*})\), and fix
  a basis \(Y_0,Y_1,..,Y_n\) of \(W^{*}\), so \(S \iso
  \bb{K}[Y_0,..,Y_n]\). Also take \(W_0,..,W_n\) to be
  the dual basis of \(W\). Given
  \(\lambda = \sum_{i} f_i \tensor W_i \in S \tensor W\), the preceding proposition
  shows \(\lambda\) defines a derivation \(\theta_\lambda \in
  \Der(S)\). Namely,
  \[\theta_\lambda(g) = \sum_{i} f_i \partialfrac{g}{Y_i}.\]

  Moreover, \(\lambda\) defines a polynomial map
  \( \rho_\lambda:W \to W, \)  or equivalently a rational map \(\bb{P}(W) \to \bb{P}(W)\),
  via
  \[\begin{aligned}\rho_\lambda(w) &= \sum_{i=0}^n f_i(w) W_i \\
      &= (f_0(w):f_1(w):..:f_n(w))
    \end{aligned}
  \]
  
  Finally, it defines a pairing \(\angvec{\;,\,}_\lambda:W \times W^{*}
  \to \bb{K}\), linear only in \(W^{*}\), where for \((s,\ell) \in W
  \times W^{*}\)
  \[ \begin{aligned}
      \angvec{s,\ell}_\lambda &:= \sum_{i=0}^n (f_i(s))(\ell(W_i)) \\
      \text{ or in coordinates; } &\\
       \angvec{(a_0,..,a_n),c_oY_0+..c_nY_n}_\lambda &:= \sum_{i=0}^n 
      f_i(a_0,..,a_n) c_i
    \end{aligned}
  \]
  We extended this definition to a pairing \(\angvec{\;,\,}_\lambda:W \times V \to \bb{K}\) via
  \[ \angvec{s,t}_\lambda := \sum_{i} f_i(s)(B(t,u_i))\]
\end{defn}

The following is immediate from the definitions

\begin{lemma}\label{lem:PairingLemma} For \((s,t) \in W \times
  V\) and \(\lambda \in S \tensor W\)
  \[ \left[\theta_\lambda(\ell_t)\right](s) = \angvec{s,t}_\lambda = \ell_t(\rho_\lambda(s)).\]
\end{lemma}

This proposition is essentially due to
Stanley, though the presentation is our own.

\begin{prop}\label{prop:StanleysCriterion}
  Let \(\lambda \in S \tensor W\), and \(\cl{A} \subseteq
  W\) a hyperplane arrangement with \(Q_{\cl{A}} = \prod_{H}
  \ell_H\). Then the following are equivalent:
  \begin{enumerate}[(i)]
  \item \(\theta_\lambda \in D(\cl{A})\)
  \item \(\theta_\lambda(\ell_H) \subseteq I(H)\) for all \(H \in
    \cl{A}\)
  \item \(\rho_\lambda(H) \subseteq H\) for all \(H \in \cl{A}\)
  \item For all \(H \in \cl{A}\), the restriction of \(\angvec{-,-}_\lambda\) to \(H \times
    H^{\perp} \subseteq W \times V\) is identically \(0\).
  \end{enumerate}
\end{prop}

\begin{proof}

  [\((i) \iff (ii)\)] The implication \((i)  \implies (ii)\) follows from
  the definition. For the converse note that \(I(H)\) is generated by \(\ell_H\), so
  every element \(f \in I(H)\) may be written \(f = g\ell_H\). 
  Applying the Leibniz product rule we get 
  \[\theta_\lambda(f) = \theta_\lambda(a_i) \ell_H + a_i\theta_\lambda(\ell_H).\]
  The first term is necessarily in \(I(H)\), and so if
  \(\theta_\lambda(\ell_H) \in I(H)\) then we conclude
  that the second sum is in \(I(H)\) as well, establishing the result.

  [\((ii) \iff (iii) \iff (iv)\)] \((ii)\) can be rephrased as
  follows:
  ``for all \(\ell \in [I(H)]_1\) and all \(p \in L\), \([\theta_\lambda(\ell)](p)
  = 0\)``.
  
  Now using the fact that \([I(H)]_1\) is naturally isomorphic to
  \(H^{\perp}\) under our isomorphism \(V \iso W^{*}\), we conclude by
  applying \ref{lem:PairingLemma}.
\end{proof}

\begin{defn}
Under the characterization above, the identity map
  on \(W\) corresponds to a derivation known
  as the \emph{Euler Derivation} which we denote \(\theta_e\). In coordinates, if \(S = \bb{K}[Y_0,..,Y_n]\), then
  \[\theta_e = Y_0 \partialfrac{}{Y_0} + Y_1\partialfrac{}{Y_1} +
    \ldots + Y_n \partialfrac{}{Y_n}.\]
\end{defn}

The Euler Derivation can be alternatively characterized as the unique derivation where \(\theta_e(f) = \deg(f) f\) for all
  homogeneous \(f\), an identity originally due to Euler.

\begin{defn}\label{defn:reducedDerivations}(Reduced Module of
  Derivations) 
  Let \(\theta_e\) denote the \emph{Euler Derivation} we define the \emph{Reduced Module of Derivations},
  denoted \(D_0(\cl{A})\), as the quotient
  \[ D_0(\cl{A}) := D(\cl{A}) / (S \theta_e).\]
  By convention, we set \(D(\emptyset) = \Der(S)\) and
  \(D_0(\emptyset) = \Der(S)/(S\theta_e)\). 
\end{defn}

\begin{defn}\label{defn:restrictedDerivations} Let \(\cl{A} \subseteq
  \bb{P}(W)\) be a hyperplane arrangement, then \(D(\cl{A})\) defines
  a reflexive sheaf, \(\widetilde{D(\cl{A})}\) on \(\bb{P}(W)\) of rank \(\dim W\). 

  If \(L \subseteq
  \bb{P}(W)\) is a line we may tensor 
  \(\widetilde{D(\cl{A})}\) with the structure sheaf
  \(\cl{O}_L\). This may equivalently be viewed as a sheaf of
  \(\cl{O}_{\bb{P}(W)}\) modules, or the restriction of
  \(\widetilde{D(\cl{A})}\) to \(L\).
  We let \(D(\cl{A}) \mid_L\) denote the corresponding graded module, that is
  \[[D(\cl{A})\mid_L]_d = H^0(\widetilde{D(\cl{A})}\tensor
    \cl{O}_{L}(-d),L)\]

  We may similarly define \(D_0(\cl{A})\mid_L\).
\end{defn}

The above object, \(D(\cl{A}) \mid_L\) has an equivalent algebraic
  definition, at least if the line \(L\) is general which we state now.

\begin{prop}For a general line \(L \subseteq \bb{P}(W)\), and for
  \(\ell \in S\) let \(\bar{\ell}\) denote the image of \(\ell\) in
  \(S/I(L)\), then 
  \[ D(\cl{A}) \mid_L = \{ \theta \in \Der(S,S/I(L)) \mid
    \theta(\ell) \in (\bar{\ell}) \textrm{ for all \(\ell\)
      dividing \(Q(\cl{A})\)} \}\]
  and similarly for \(D_0(\cl{A})\).
\end{prop}

\begin{proof} First, for any \(f \in S\) we let \(\bar{f}\) denote
  the image of \(f\) in \(S/I(L)\), and similarly if \(\theta =
  \sum_{i=0}^n f_i \partialfrac{}{Y_i}\) we let \(\bar{\theta} =
  \sum_{i=0}^n \bar{f_i} \partialfrac{}{Y_i} \in \Der(S,S/I(L))\).

  Note that \(D(\emptyset) \mid_L \) is
  isomorphic to \(\Der(S,S/I(L))\), and so
  \(D(\cl{A})\mid_L\) is isomorphic to a submodule of
  \(\Der(S,S/I(L))\). 

  Now consider the case that \(\cl{A}\) consists of a
  single hyperplane \(H\). Choosing our coordinates \(Y_0,..,Y_n\)
  so that \(H = (Y_0=0)\),  then \(D(\cl{A})\) is free on generators
  \(\{Y_0 \partialfrac{}{Y_0},\partialfrac{}{Y_1}, \ldots ,
  \partialfrac{}{Y_n}\}\). Then \(D(\cl{A})\mid_L\) is a free
  \(S/I(L)\) module with basis
  \(\{\bar{Y_0} \partialfrac{}{Y_0},\partialfrac{}{Y_1}, \ldots ,
  \partialfrac{}{Y_n}\}\).  Yet these are also precisely the derivations \(\theta
  \in \Der(S,S/I(L))\) where \(\theta(X_0) \in (\bar{X_0})\), so the
  result follows in this case.

  More generally, if \(\cl{A} = \{H_0,H_1,\ldots,H_k\}\) and \(L\)
  is any line not contained in a hyperplane in \(\cl{A}\). Then for all \(i \neq j\) we have that \(L \cap H_i\) and
  \(L \cap H_j\) consist of distinct points. Consequently, letting
  \(U_i\) denote the complement of \(\cl{A} \setminus \{H_i\}\),
  we have that \(\{U_i \cap L\}_{i = 0,..,k}\) is an open cover
  of \(L\). Therefore, for a section \(\sigma \in H^{0}(\widetilde{\Der(S)} \tensor
  \cl{O}_L(-d), L)\), we have that \(\sigma \in D(\cl{A})\mid_L\) if
  and only if
  \[\sigma \mid_{U_i} \in  H^0(\widetilde{D(\cl{A})}\tensor
    \cl{O}_L(-d),U_i \cap L) =  H^0(\widetilde{D(H_i)}\tensor
    \cl{O}_L(-d),U_i \cap L)\] for all \(j=1,..,k\).

 Finally, note that for \(i \neq j\) we have that \(H^0(\widetilde{D(H_i)}\tensor
    \cl{O}_L(-d),U_j \cap L) = \sigma \in H^{0}(\widetilde{\Der(S)} \tensor
  \cl{O}_L(-d), U_j \cap L)\). Therefore, it follows that for \(\sigma
  \in H^{0}(\widetilde{\Der(S)} \tensor
  \cl{O}_L(-d), L)\), that we have the following string of
  equivalences
  \[\begin{aligned}\sigma \in 
      H^0(\widetilde{D(\cl{A})} \tensor \cl{O}_L(-d),L) & \iff \\
      \textrm{ for all } i \in \{1,..,k\}, \sigma \mid_{U_i \cap L} \in H^0(\widetilde{D(H_i)} \tensor
      \cl{O}_L(-d),L \cap U_i) &
      \iff \\
      \textrm{ for all } i \in \{1,..,k\}, \sigma \in H^0(\widetilde{D(H_i)} \tensor
      \cl{O}_L(-d),L) &
\end{aligned}\]
The result now follows from the previous case.
\end{proof}

We can emulate the constructions from \ref{defn:StanleysRho} for
the module \(D_0(\cl{A})\mid_L\), to achieve a characterization similar to \(\ref{prop:StanleysCriterion}\).

\begin{defn}\label{defn:StanleysRestrictedRho}  Let \(L \subseteq
  \bb{P}(W)\) be a line, for \(\gamma =\sum_{j} f_j \tensor w_j \in S/I(L) \tensor W\), we obtain a pairing \(\angvec{\; ,\,}_\gamma:L
  \times V \to \bb{K}\), defined by
  \[ \angvec{p,v}_\gamma := \sum_{j} f_j(p) (\ell_v(w_j)) \]

  Similarly, define the polynomial map \(\rho_\gamma:L \to W\) and
  \(\theta_\gamma \in \Der(S,S/I(L))\).
\end{defn}

\begin{prop}\label{prop:StanleysRestrictedCriterion}
  Let \(L \subseteq \bb{P}(W)\) be a general line, and \(\cl{A}
  \subseteq \bb{P}^n\) a hyperplane arrangement, then for \(\gamma
  \in S/I(L) \tensor_\bb{K} W\), the following are equivalent.
  \begin{enumerate}[(i)]
  \item \(\theta_\gamma \in D_0(\cl{A})\mid_L\)
  \item \(\rho(L \cap H) \subseteq H\) for all \(H \in \cl{A}\)
  \item The restriction of \(\angvec{\;,\,}_\gamma\) to \((H \cap L)
    \times H^{\perp}\) is identically \(0\).
  \end{enumerate}
\end{prop}

\begin{proof} The proof is essentially identical to that of
  \cref{prop:StanleysCriterion} so we omit it. Note that in particular,
  we still have an analogue of \cref{lem:PairingLemma} for \((p,q) \in L \times V\) that
  \[\angvec{p,q}_\gamma = \theta_\gamma(\ell_q)(p) = \ell_q(\rho_\gamma(p))\].
\end{proof}

If \(M\) is a finite reflexive graded module over \(\bb{K}[Y_0,..,Y_n]\) 
defining a reflexive sheaf \(\tilde{M}\) on
\(\bb{P}^n_{\bb{K}}\). Then the restriction of \(\tilde{M}\) to a
general line \(L \subseteq \bb{P}^n_{\bb{K}}\), defines a vector bundle
  over \(L\).   By the well known theorem of Birkhoff and Grothendieck, there exist integers
  \(k_0 \leq k_1 \leq k_2 \leq .. \leq k_m\) where
  \[\tilde{M} \mid_L \iso \oplus_{i=0}^m \cal{O}_L(-k_i).\]

  If \(\cl{A} \subseteq \bb{P}^n_{\bb{K}}\) is a hyperplane
  arrangement, then \(D(\cl{A})\) can be naturally identified with
  the first syzygy module of the ideal \(J = (Q_{\cl{A}},
  \partialfrac{}{Y_0} Q_{\cl{A}}, \partialfrac{}{Y_1} Q_{\cl{A}},
  \ldots, \partialfrac{}{Y_n} Q_{\cl{A}})\). This ensure that
  \(D(\cl{A})\) is reflexive.

  We show that \(D_0(\cl{A})\) is reflexive for any nonempty
  arrangements \(\cl{A} \subseteq \bb{P}^n_{\bb{K}}\). If \(|\cl{A}| \neq 0 \mod
  \chr \bb{K}\), this is well known as \(J = \textrm{Jac}(Q_{\cl{A}})
  = (\partialfrac{}{Y_0} Q_{\cl{A}}, \partialfrac{}{Y_1} Q_{\cl{A}},
  \ldots, \partialfrac{}{Y_n} Q_{\cl{A}})\) and \(D_0(\cl{A})\) can be
  identified with the syzygy module of \(\textrm{Jac}(Q_{\cl{A}})\).

  We establish this more generally, our proof requires the following
  reflexive criterion. We refer to \cite{STK} (see \href{https://stacks.math.columbia.edu/tag/0EB8}{Lemma 15.23.5})  for a proof.

  \begin{prop}[Reflexive Criterion] Suppose
    \[\xymatrix{0 \ar[r] & M \ar[r] & L \ar[r] & K }\]
    is an exact sequence of finite modules, over a
    commutative notherian domain \(R\). Then if \(L\) is reflexive and
    \(K\) is torsion free, then \(M\) is reflexive.
  \end{prop}

  With this criteria we can establish our claim. This is well known
  for arrangements over \(\bb{C}\) and likely in general we include it for completeness.
  
\begin{prop} If \(\cl{A} \subseteq \bb{P}^n_{\bb{K}} = \Proj(\bb{K}[Y_0,..,Y_n])\) is a nonempty
  hyperplane arrangement, then \(D_0(\cl{A})\) is a reflexive module.
\end{prop}

\begin{proof} The proof is by induction on the number of hyperplanes
  in \(\cl{A}\). First we consider the case \(|\cl{A}| = 1\) or
  \(|\cl{A} = 2|\), in these cases we can choose coordinates so that
  \(Q_{\cl{A}} = Y_0\) or \(Q_{\cl{A}} = Y_0Y_1\) respectively. It can
  now be checked by direct computation that \(D_0(\cl{A})\) is free on
  generators \(\{\partialfrac{}{Y_1},\ldots, \partialfrac{}{Y_n}\}\)
  and \(\{Y_1\partialfrac{}{Y_1},\ldots, \partialfrac{}{Y_n}\}\)
  respectively.

  For the general case if \(\cl{A}'\) is a hyperplane arrangement with
  \(k > 2\) hyperplanes pick two distinct hyperplane \(L\) and \(H\) in
  \(\cl{A}\). Let \(\cl{A} = \cl{A}' \setminus \{H\}\) and let \(\cl{B}\) denote the
  hyperplane arrangement \(\{L,H\}\). Then we have the following exact sequence
  \[\xymatrix{0 \ar[r] & D_0(\cl{A'}) \ar[r] & D_0(\cl{A}) \oplus
      D_0(\cl{B}) \ar[r] & D_0(\{L\}) }.\]

  As \(D_0(\{L\})\) is free it is in particular torsion free. Furthermore by inductive hypothesis \(D_0(\cl{A})\) and
  \(D_0(\cl{B})\) are both reflexive so we conclude by applying the preceding proposition.
\end{proof}

\begin{defn}[Splitting Type] \label{defn:SplittingType} 
  If \(\cl{A} \subseteq \bb{P}^n\) is a hyperplane
  arrangement (resp. nonempty hyperplane arrangement), then there exists 
  tuple of integers \((a_0,a_1,..,a_n), \; \left(\textrm{resp.} \;
    (a_1,..,a_n)\right)\) referred to as the \emph{Splitting Type} of
  \(D(\cl{A}) \; \left(\textrm{resp.} \; D_0(\cl{A})\right) \).
  
  This is the unique tuple satisfying \( 0 \leq a_0 \leq a_1 \leq
  ... \leq a_n\), so that if \(L\) is a general line then there's an isomorphism 
  \[\begin{aligned}
      & D(\cl{A}) \mid_L \iso \bigoplus_{i=0}^n S/I(L)(-a_i) &
      \hspace{0.5in} & 
    \left(\textrm{resp.} \; D_0(\cl{A}) \mid_L \iso \bigoplus_{i=1}^n
      S/I(L)(-a_i)\right).\\
  \end{aligned}\]

\end{defn}

\PFEIndent

\section{Derivation Bundle of Hyperplane Arrangements and the Ideals
  of Dual Points}\label{sec:DerivationIdealCorrespondence}

In this section we introduce our duality and establish a relationship between
\(D_0(\cl{A}_Z)\) and \(I(Z)\). We can summarize this relationship as
follows:
Given a set of points \(Z \subseteq \bb{P}(W)\) with
 dual hyperplane arrangement \(\cl{A}_Z \subseteq \bb{P}(V)\), we consider
 a ring \(T = R \tensor_{\bb{K}} \bb{K}[\Gr(n-2,V)]\) that contains naturally isomorphic copies of \(R =
 \Sym(V^*)\) and \(S = \Sym(W^{*})\).  We then show in
 \cref{thm:GlobalIso} that
 \(D_0(\cl{A})\) is isomorphic to an \(S\)-submodule of the extended
ideal \(I(Z) T\). This is analogous to the standard construction used
in \cite{FV}. In \cref{thm:restrictedVectIso} we then give a novel interpretation of the restriction
of this \(S\)-submodule to a general line.

\begin{defn}\label{defn:ExteriorAlg} Let \(\bigwedge^{\bullet} V\), denote the exterior
  algebra of \(V\). This is the graded \(\bb{K}\)-algebra generated in
  degree \(1\) by \(V\), subject to the relation \(v^2 = v \wedge v =
  0\) for all \(v \in V\).
\end{defn}

\begin{defn}\label{defn:Grassmanian}\label{defn:Plucker} Let \(\Gr(k,V)\) denote the \(k\)-th grassmanian of
  \(V\) as a projective subvariety of \(\bb{P}(\bigwedge^{k} V)\). The
  projective coordinate ring of \(\Gr(k,V)\) as a quotient of the
  polynomial ring of the ambient space is the \emph{Pl\"{u}cker
    Algebra}, \(\Pl(k,V)\).
\end{defn}

  Fix a set of coordinates \(X_0,..,X_n\) on \(V\), so that \(\Sym(V^{*}) \iso \bb{K}[X_0,..,X_n]\). Extend these
  to coordinates on \(V^{\oplus k}\) for some \(1 \leq k
  \leq n\), by letting \(A_{i,0},..,A_{i,n}\) denote
  an isomorphic copy of \(X_0,..,X_n\), for each \(i \in \{0,..,k-1\}\).
  We organize these into a \(k \by n+1\) matrix
  \(\mathbf{A}\) with entries \((\mathbf{A})_{i,j} = A_{i,j}\).

  Let \(c(\Gr(k,V))\) denote the affine cone of
  \(\Gr(k,V)\) as a subvariety of \(\bigwedge^{k} V\). Then the
  multiplication map \(\wedge:V^{\oplus k} \to c(\Gr(k,V)) \subseteq
  \bigwedge^{k} V\), identifies the Pl\"{u}cker algebra \(\Pl(k,V)\) with
  the \(\bb{K}\) algebra generated by the maximal \((k \by k)\)
  minors of \(\mathbf{A}\). 

  Restricting to the case where \(k = n\), multiplication in \(\bigwedge V\) gives a non-degenerate pairing
  \(\wedge: V \cross \bigwedge^{n} V \to \bigwedge^{n+1}
  V\). Choosing an isomorphism \(\bigwedge^{n+1} V \iso \bb{K}\)
  gives an isomorphism \(\bigwedge^{n} V \iso V^{*}\), natural up to a
  \(\bb{K}\)-scalar. We fix one of these isomorphisms and let \(\tau\)
  denote the induced isomorphism of polynomial rings \(\tau:
  \Sym(W^{*}) \iso \Sym(\bigwedge^{n} V^{*}) \). As \(n = \dim V -1\),
  then \(\bigwedge^{n} V^{*} = \Gr(n,V)\), and we can identify
  \(\Pl(n,V)\) with \(\Sym(\bigwedge^n V^*) \iso \Sym(W^{*})\).  We further describe
  \(\tau\) in coordinates below.

\begin{defn}\label{defn:TauIso} Taking the definitions of
  \(X_i\) and \(A_{j,\ell}\) from above, further require that
  \(A_{0,i} = X_i\). Define
  \[\begin{aligned}
      \bb{K}[\mathbf{A}] &:= \bb{K}[A_{i,j} \mid 0 \leq i \leq n-1 , 0
      \leq j \leq n] \\
       &:= \bb{K}[X_0,..,X_n][A_{i,j} \mid 1 \leq i \leq n-1, 0 \leq j
       \leq n]
     \end{aligned}\]

  Let \(\Pl(n)\) be the subalgebra of
  \(\bb{K}[\mathbf{A}]\) generated by the determinants
  \[M_i :=  \left| \begin{matrix}0 & \ldots & 0 & 1 & 0 & \ldots & 0\\
        A_{0,0} &  \ldots & A_{0,i-1} & A_{0,i} & A_{0,i+1} & \ldots  & A_{0,n}\\
        A_{1,0} & \ldots & & A_{1,i} &  & \ldots & A_{1,n}\\
        \vdots & & & \vdots & & & \vdots\\
        A_{n-1,0} & \ldots & & A_{n-1,i} &  & \ldots & A_{n-1,n}
      \end{matrix}\right| = \left| \begin{matrix}0 & \ldots & 0 & 1 & 0 & \ldots & 0\\
        X_0 &  \ldots & X_{i-1} & X_i & X_{i+1} & \ldots  & X_n\\
        A_{1,0} & \ldots & & A_{1,i} &  & \ldots & A_{1,n}\\
        \vdots & & & \vdots & & & \vdots\\
        A_{n-1,0} & \ldots & & A_{n-1,i} &  & \ldots & A_{n-1,n}
      \end{matrix}\right|.\]

  Finally, taking \(Y_i\) to be a dual basis of \(X_i\) we define
  \(\tau : \bb{K}[Y_0,..,Y_n] \to \Pl(n)\) via \(\tau(Y_i) = M_i\).
\end{defn}
  
The preceding conversation shows that \(\Pl(n))\) is a polynomial
algebra in the generators \(M_i\). The lemma below shows that our
definition of \(\tau\) above matches the 
construction from the preceding remark.

\begin{lemma}\label{lemma:LinearVanishing} Let \(v \in V =
  \Spec(\bb{K}[X_0,..,X_n])\) and let \(\ell_v = \sum_{i=0}^n c_i Y_i \in W^{*}\)
  be the corresponding linear form. Then as a
  polynomial in \(X_0,..,X_n\) the linear form \(\tau(\ell_v) = \sum_{i} c_i M_i\) vanishes on \(v\).
\end{lemma}

\begin{proof} Following \cref{defn:TauIso} we see that
  \(\tau(\ell_v) = \sum_{i=0}^n c_i M_i\) is the 
  Laplace expansion along the first row of the determinant of the
  matrix \(\begin{bmatrix} \vec{v} \\ \mathbf{A}\end{bmatrix}\), where
   \(\vec{v} = \begin{bmatrix} c_0 & c_1 & ... &
    c_n\end{bmatrix}\). If we then evaluate
  \(X_0,..,X_n\) at \(v\), so that \(X_i \mapsto c_i\), the matrix is
  singular as two rows are identical hence the determinant vanishes.  
\end{proof}

\begin{defn}\label{defn:Flambda} If \(X_0,..,X_n\) form a basis of
  \(V^*\) and \(Y_0,..,Y_n\) are dual coordinates on \(W^*\).  Then for any \(\lambda = \sum_{i=0}^n f_i(Y_0,..,Y_n)
   \ell_{X_i} \in \Sym(W^{*}) \tensor W\), we define a polynomial
  \(F_\lambda \in  \bb{K}[\mathbf{A}]\) via
  \[F_\lambda := \sum_{i=0}^n X_i \tau(f_i) = \sum_{i=0}^n X_i f_i(M_0,..,M_n)\]
\end{defn}

\begin{defn}\label{defn:ShiftedIdealsGlobal}
   Let \(J \subseteq \Sym(V^{*}) = \bb{K}[X_0,\ldots,X_n]\) be any homogeneous
  ideal. We define a graded module denoted \(J^{\gg}\) over
  \(\Pl(n)\), thought of as a polynomial ring in the minors
  \(M_0,..,M_n\).

  First, if \(\mathfrak{m}=(X_0,X_1,\ldots,X_n)\) is the maximal
  ideal we define \(\mathfrak{m}^{\gg}\) as the \(\Pl(n)\)-submodule of
  \(\bb{K}[\bm{A}]\) generated by \((X_0,\ldots,X_n)\). We grade both
  \(\Pl(n)\) and \(\mathfrak{m}^{\gg}\) by the \(X\)-degree, meaning
  \(\deg(M_i) = \deg(X_i) = 1\). Equivalently, the  \(d\)-th graded
  component of \(\mathfrak{m}^{\gg}\) is generated over \(\bb{K}\) by all
  terms of the form 
  \(X_i M_0^{e_0}M_1^{e_1} \ldots M_n^{e_n}\) where \(\sum_{i=0}^n e_i
  = d-1\).

  More generally for any homogeneous ideal \(J\), we set \(J^{\gg} = (J \bb{K}[\bm{A}]) \cap \mathfrak{m}^{\gg}\).
\end{defn}

The following example and proposition that follows are the main
motivation for the definition of \(J^{\gg}\).

\begin{example} If \(P\) is the \(0\)-th coordinate point in
  \(\bb{P}^2 = \Proj(\bb{K}[X_0,X_1,X_2])\), then \(I(P) =
  (X_1,X_2)\). Given \(\sum_{i=0}^2 G_i X_i \in \mathfrak{m}^{\gg}\)
  where \(G_i\) is a polynomial in the maximal minors \(M_0,M_1,M_2\) of the
  matrix
  \[\begin{bmatrix}
      X_0 & X_1 & X_2 \\
      A_0 & A_1 & A_2
    \end{bmatrix}.\]
    
  It's not hard to see
  that \(\sum_{i=0}^n G_i X_i \in I(P)^{\gg}\) if and only if \(G_0
  \in I(P) \bb{K}[\bm{A}]\).

  Treating \(G_0\) as a polynomial in
  \(X_0,..,X_n\) with coefficients in the ring \(\bb{K}[A_0,A_1,A_2]\)
  and consider the evaluation map \(e_P: \bb{K}[\mathbf{A}]
  \to \bb{K}[A_0, A_1,A_2]\) obtained by evaluting at \(P\), we note that
  \(G_0 \in I(P) \bb{K}[\mathbf{A}]\) if and only if
  \(\varepsilon_P(G_0) = 0\). Furthermore, as \(e_P\) 
   sends \(M_0 \mapsto 0\), \(M_1 \mapsto -A_2\) and
  \(M_2 \mapsto A_1\) then \(G_0(P) = 0\) if and only if \(M_0\)
  divides \(G_0\). From this it follows that \(I(P)^{\gg}\) is 
  generated by \(\{X_0 M_0, X_1,X_2\}\).

  In fact this generating set is
  redundant as we have the nontrivial relation \(X_0M_0 + X_1M_1 +
  X_2M_2=0\), and so a minimal generating set for \(I(P)^{\gg}\) is given by \(\{X_1,X_2\}\).
\end{example}

The preceding definition is motivated by the following proposition.

\begin{theorem}\label{thm:GlobalIso} Let \(Z \subseteq \bb{P}(V) = \Proj(\bb{K}[X_0,..,X_n])\)
  be a finite set of points, and
  let \(\cl{A}_Z \subseteq \bb{P}(W)\) denote the dual hyperplane arrangement. Then
  for \(\lambda \in S \tensor W\) the following are
  equivalent:
  \begin{enumerate}[(i)]
  \item \(\theta_\lambda \in D(\cl{A}_Z)\)
  \item \(F_\lambda \in I(Z) \cdot \bb{K}[\mathbf{A}]\)
  \end{enumerate}

  Moreover, \(F_\lambda = 0\) if and only if there exists \(g \in
  S\) so that \(\theta_\lambda = g
  \theta_e\), where \(\theta_e = \sum_{i=0}^n
  Y_i\partialfrac{}{Y_i}\) is the Euler derivation.

  In essence there's an isomorphism \(\eta: \Pl(n,V) \tensor_S D_0(\cl{A}_Z)(-1) \to
  I^{\gg}(Z)\) given by
  \[\eta\left(\sum_{i=0}^n f_i(Y_0,..,Y_n) \partialfrac{}{Y_i}
    \right) =\sum_{i=0}^n f_i(M_0,..,M_n) X_i \]
\end{theorem}

The above theorem is a consequence of the following lemma which is
useful in it's own right.

\begin{lemma}\label{lem:HypFactorization} Fix \(\bm{\alpha} = (\alpha_1,..,\alpha_{n-1})\) a
  tuple of \((n-1)\) linearly independent vectors in \(V\). Letting \(\alpha_i :=
  (\alpha_{i,0},\alpha_{i,1},\ldots,\alpha_{i,n})\) in our chosen set
  of coordinates. We define the partial evaluation map
  \[\varepsilon_{\bm{\alpha}}:\bb{K}[\mathbf{A}] \to \bb{K}[X_0,..,X_n]\]
  via \(\varepsilon_{\bm{\alpha}}(A_{i,j}) = \alpha_{i,j}\) for \(1
  \leq i \leq n-1\).

  Let \(\lambda =\sum_{i=0}^n f_iX_i\in [S \tensor W]_d = \Sym^d(W^{*}) \tensor W\), 
  for any nonzero \(w \in W\) where \(\ell_w \in V^{*}\) vanishes on
  \(\Spn(\bm{\alpha})\), there exists some nonzero linear form \(h\)
  vanishing on \(\Spn(\bm{\alpha})\) so that
  \[\varepsilon_{\bm{\alpha}}(F_\lambda) \equiv h^{d}
    \ell_{\rho_\lambda(w)} =h^d \left(\sum_{i=0}^n X_i f_i(\rho_\lambda(w))\right)
    \mod I(\ell_w)\]
\end{lemma}

\begin{proof}
Take \(w,\bm{\alpha}\) and \(\lambda\) as stated above.
  Let \(\lambda = \sum_{i=0}^n f_i(Y_0,..,Y_n) \tensor X_i \in [S \tensor
  W]_d\), then \(F_\lambda = \sum_{i=0}^n X_i f_i(M_0,..,M_n)\).

  Write \(\ell_w = \sum_{i=0}^n c_iX_i\) and assume without loss of
  generality that \(c_n \neq 0\). Fix some index \(j \in
  \{0,..,n-1\}\) and let \(\ell_u = c_nY_j -
  c_jY_n \in [S]_1 = W^{*}\). Noting that \(\ell_u(w) = \ell_w(u) = 0\),
  we may write \(\varepsilon_\alpha(\tau(\ell_u))\) as the determinant of the matrix
\[ \varepsilon_\alpha(\tau(\ell_u)) = \left| \begin{matrix}0 & \ldots & 0 & c_n & 0 & \ldots & -c_j\\
X_{0} &  \ldots & X_{j-1} & X_{j} & X_{j+1} & \ldots  & X_{n}\\
\alpha_{1,0} & \ldots & & \alpha_{1,j} &  & \ldots & \alpha_{1,n}\\
\vdots & & & \vdots & & & \vdots\\
\alpha_{n-1,0} & \ldots & & \alpha_{n-1,j} &  & \ldots & \alpha_{n-1,n}
\end{matrix}\right|  \]
As \(u \in \ker \ell_w\) and \(\Spn(\alpha) \subseteq \ker w\), then
either \(\varepsilon_{\bm{\alpha}}(\tau(\ell_u))=0\) and \(u \in \Spn(\alpha)\), or
\(\varepsilon_{\bm{\alpha}}(\tau(\ell_u))\neq 0\) and \(\Spn(\alpha,u)
= \ker \ell_w\) which implies theres a scalar \(r \in
\bb{K}\) so
\(\varepsilon_\alpha(\tau(\ell_u)) = r \ell_w\). In
either case we have \(c_n M_i - c_i M_n \equiv 0 \mod (\ell_w)\).
We conclude with the equalities below where here \(M_i =
\varepsilon_{\bm{\alpha}} (M_i)\),
\[\begin{aligned}
    \varepsilon_\alpha(F_\lambda) &= \sum_{i=0}^n
  f_i(M_0,M_1,..,M_n) X_i\\
    &\equiv \sum_{i=0}^n
  f_i\left(\frac{c_0}{c_n}M_n,\frac{c_1}{c_n}M_n,..,\frac{c_{n-1}}{c_n}M_{n},M_n\right)
    X_i & \mod (\ell_w)\\
    &\equiv \left(\frac{M_n}{c_n}\right)^{d}\sum_{i=0} f_i(c_0,..,c_n)
    X_i & \mod (\ell_w)\\
    &\equiv \left(\frac{M_n}{c_n}\right)^{d} \rho_\lambda(w) &\mod (\ell_w)
  \end{aligned}\]

Noting that because \(c_n \neq 0\), we must have that \(E_n =
(0:..:0:1) \not\in \Spn(\bm{\alpha}) \subseteq \ker
\ell_w\). Therefore, \(\varepsilon_{\bm{\alpha}}(M_n) \neq 0\) as it is
the determinant of a non-singular matrix thereby establishing the result.
\end{proof}

\begin{proof}[Proof of \cref{thm:GlobalIso}] Let \(\lambda = \sum_{i}
  f_i \tensor w_i \in S \tensor W\).  We note that since \(D(\cl{A}_Z) = \bigcap_{P
  \in Z} D(H_P)\) and \(I(Z) = \bigcap_{P \in Z} I(P)\), it suffices to
establish the equivalence in consider the case \(Z\) consists of a
single point \(P\). Furthermore, to establish the case 
for a single point it suffices to show that
\(\epsilon_\alpha(F_\lambda)\) is in \(I(P)\) for every (or even for
general) \(\alpha\). This is because \(\theta \in \Der(S)\) is in \(D(\cl{A}_Z)\)
if and only if  the restriction of \(\theta\) to \(L\) is in
\(D(\cl{A})\mid _L \)  for general \(L\), and similarly \(F_\lambda\)
vanishes at \(P\) if and only if
\(\varepsilon_{\bm{\alpha}}(F_\lambda)\) vanishes on \(P\) for general \(\bm{\alpha}\).

Continuing, assume that \(\bm{\alpha}\) is sufficiently general and let \(\ell_Q\) denote the
linear form vanishing on \(\bm{\alpha}\) and \(P\). We consider
\(\varepsilon_{\bm{\alpha}}(F_\lambda)  \mod (\ell_Q)\). By
\cref{lem:HypFactorization}, we get that
\(\varepsilon_{\bm{\alpha}}(F_\lambda) \equiv h^d
\ell_{\rho_\lambda(Q)} = h^d \rho_\lambda(\ell_Q)\). 
Yet for general \(\bm{\alpha}\), we see that \(h(P) \neq 0\) so \(F_\lambda\) vanishes
on \(P\) if and only if \(\rho_\lambda(\ell_Q)\) vanishes on \(P\).

Now for any linear form \(\ell_L\) recall that \(\ell_L(P) = 0\) if
and only if the corresponding \(L \in W\) lies on \(P^{\perp} =
H_P\). Hence, we apply \cref{prop:StanleysCriterion} and conclude the
proof of the first statement with
the following chain of equivalences:

\[\begin{aligned}
    F_\lambda \in I(P) & \iff \text{ for general } \bm{\alpha}, \,
    \varepsilon_{\bm{\alpha}}(F_\lambda) \in I(P) \\
&\iff \text{ for general } \bm{\alpha}, \, \rho_\lambda(\ell_Q) \in I(P)  \text{
  where } \ell_Q \text{ vanishes on }\Spn(\bm{\alpha},P)
\\&\iff \text{ for general } Q \in H_p, \,  \rho_\lambda(Q) \in
H_P \\
&\iff \theta_\lambda \in D(H_P)
\end{aligned}\]

To finish the proof, we must establish the claim about the kernel of \(\eta\).
We see that \(F_\lambda = 0 \)
if and only if for general \(\bm{\alpha}\) and arbitrary \(\ell_H\)
vanishing on \(\bm{\alpha}\) that
\(\varepsilon_{\bm{\alpha}}(F_\lambda) \equiv 0 \mod I(\ell_H)\). By
\cref{lem:HypFactorization} the later condition occurs precisely when \(\rho_\lambda(\ell_H) \in (\ell_H)\) for every linear form
\(\ell_H\). If this occurs we conclude for all \(H \in W\) that
\(\rho_\lambda(H) = r_H H\) for some scalar \(r_H\). It immediately
follows that as a rational map on \(\bb{P}(W)\), \(\rho_\lambda\) can
be extended to the identity, allowing us to conclude that
\(\theta_\lambda = f \theta_e\) where \(\theta_e\) is the Euler derivation.
\end{proof}

We note that the proof above also establishes the following.

\begin{theorem}\label{thm:restrictedVectIso} Let \(Z \subseteq \bb{P}(V)\)
  and let \(\cl{A}_Z  \subseteq \bb{P}(W)\) be the dual hyperplane
  arrangement. Let \(L \subseteq \bb{P}(W)\) be a general line, and
  \(Q = L^{\perp} \subseteq \bb{P}(V)\) the dual linear subspace. Then
  there's an isomorphism of vector spaces  
  \[[I(Z) \cap I(Q)^{m}]_{m+1} \iso [D_0(\cl{A}_Z)\mid_L]_{m}.\]
\end{theorem}

We can in fact, prove a slightly stronger
statement. Namely, the above isomorphism corresponds to an isomorphism of
modules over naturally isomorphic (up to scalar) rings, we give this
proof after \ref{thm:restrictedModIso}. In order to make this stronger
statement and to aid with the exposition for the rest of the paper, we
introduce some new notation.

\begin{defn}
  For \(Q \subseteq \bb{P}^n\) a codimension \(2\) subspace we define
  a ring \(\cl{F}_Q\) via
  \[\cl{F}_Q := \Sym_{\bb{K}}([I(Q)]_1).\]

  We note that if \(L_0,L_1\) are
  linear forms which generate \(I(Q)\), then \(\cl{F}_Q\)
  is a polynomial ring in the generators \(L_0,L_1\).
\end{defn}

Fix \(Q\) and let \(\bm{\alpha}\) be any basis of \(Q\). The following
proposition shows that \(\cl{F}_Q\) can be viewed yet another way, as
the image of the map  \(\varepsilon_{\bm{\alpha}} : \Sym(W^{*}) \to \Sym(V^{*})\).  

\begin{prop}\label{prop:TauIso}\label{prop:epsilonIso} Let \(Q = \Spn(\bm{\alpha})\), and
  \(L = Q^{\perp}\), then the map \(\varepsilon_{\bm{\alpha}} : S=\Sym(W^{*}) \to
  R = \Sym(V^{*})\), induces an isomorphism of \(\bb{K}\)-algebras
  \[\tau_{\bm{\alpha}}: S/I(L) \to \cl{F}_Q.\]

\end{prop}

\begin{proof} First consider the restriction of \(\tau_{\bm{\alpha}}\)
  as a map \([S]_1 \to [R]_1\). By \cref{lemma:LinearVanishing}, \(\varepsilon_{\bm{\alpha}}(\ell)\) must
  vanish on all points of \(Q = \Spn(\bm{\alpha})\), hence
  \(\varepsilon_{\bm{\alpha}}(\ell) \in I(Q)\). In fact, given \(P \in
  \bb{P}(V) \setminus Q\), we see again by \cref{lemma:LinearVanishing}, that
  \(\varepsilon_{\bm{\alpha}}(\ell_P)\) defines the hyperplane
  \(\Spn(Q,P)\). It follows that \(\tau_{\bm{\alpha}}\) induces an 
    isomorphism of vector spaces \([\Sym(W^{*})/I(L)]_1 \iso [\cl{F}_Q]_1\).

  As \(\Sym(W^{*})\) and \(\cl{F}_Q\) are both symmetric
  algebras generated over
  \(\bb{K}\) in degree \(1\). The isomorphism
  \(\tau_{\bm{\alpha}}:\Sym(W^{*}) \to  \cl{F}_Q\) follows.
\end{proof}

\begin{defn}\label{defn:ShiftedIdeals} Let \(Q \subseteq \bb{P}(V)\) be a codimension \(2\)
  subspace, and let \(J \subseteq \Sym(V^{*})\) be any homogeneous
  ideal. We define a graded \(\cl{F}_Q\)-module, \(J^{\gg}_Q\), as
  the \(\cl{F}_Q\)-submodule of \(\Sym(V^{*})\)
  whose \(d\)-th graded component is given by
  \[[I^{\gg}_Q(Z)]_d := [I(Z) \cap I(Q)^{d-1}]_d.\]
\end{defn}

We state the full version of this duality.

\begin{theorem}\label{thm:restrictedModIso} Let \(Z \subseteq
  \bb{P}(V) = \proj{R}\) be a finite set of points and \(\cl{A}_Z \subseteq
  \bb{P}(W) = \proj{S}\) the dual hyperplane arrangement. Let \(L \subseteq
  \bb{P}(W)\) be a general line, then the isomorphism of \(\bb{K}\) algebras
  \(\tau_Q: S/I(L) \iso \cl{F}_Q =
  \Sym([I(Q)]_1)\), induces an isomorphism of graded modules
  \(I^{\gg}_Q(Z)(-1) \iso D_0(\cl{A}_Z) \mid_L \tensor_{S/I(L)} \cl{F}_Q\)
  via the map
  \[\eta_Q:D_0(\cl{A}_Z) \mid_L \tensor
    \cl{F}_Q \iso I^{\gg}_Q(Z)(-1)\]
  \[\eta_Q\left(\sum_{i=0}^n f_i \partialfrac{}{Y_i}\right) = \sum_{i}
    \tau_Q(f_i) X_i \]

  Here \(\{Y_i\}_{i \in [n+1]}\) and \(\{X_i\}_{i \in [n+1]}\) are dual bases of \(W^{*}\)
  and \(V^{*}\) respectively.
\end{theorem}

\begin{proof} This proof is very similar to the proof of
  \cref{thm:GlobalIso}. We make note of some of the differences. Given
  \(\gamma \in S/I(L) \tensor W\), we get both a rational map
  \(\rho_{\gamma}:L \to \bb{P}(W)\) and a derivation 
  \(\theta_\gamma\) of \(\Sym(W^{*})\) into \(\Sym(W^{*})/I(L)\).

  Similarly, we get a polynomial \(F_\gamma \in
  I^{\gg}_Q(\emptyset)\) uniquely determined up to scalar. Now it again follows that \(\theta_{\gamma}
  \in D_0(\cl{A}) \mid_L\) if and only if \(\rho_{\gamma}(H \cap L)
  \subseteq H\) for all \(H \in \cl{A}\). Additionally, we have that for any \(\ell \in [I(Q)]_1\), that
  \(F_\gamma = \ell_Q^{d-1} \rho_\gamma(\ell) \mod (\ell)\).

   The proof now continues as in \cref{thm:GlobalIso}.
\end{proof}

Applying this isomorphism of modules, we note that the splitting type
of of \(\cl{A}_Z\) determines the dimension of \(I^{\gg}_Q(Z)\).

\begin{cor}\label{cor:DimensionFormula} If \(D_0(\cl{A}_Z)\) has
  splitting type \((a_1,a_2,..,a_n)\), then for a general codimension
  \(2\) linear subspace, 
  \[\dim [I_Q^{\gg}(Z)]_d = \sum_{i=1}^n \max\{0,d-a_i\}\] 
\end{cor}

\begin{example}[Ceva Arrangements] Fix \(m \geq 3\) and suppose that \(\bb{K} = \bb{C}\) or
  more generally that \(\chr(\bb{K})\) is coprime from \(m\).
  \(C_m\) denote the Ceva Arrangement of Hyperplanes  in
  \(\bb{P}^n_{\bb{K}} = \bb{P}(W)\). Taking \(\Sym(W^{*}) =
  \bb{C}[Y_0,Y_1,\ldots,Y_n]\) this is the arrangement of
  \(n+1+m\binom{n+1}{2}\) hyperplanes with defining
  polynomial
  \[Q_{C_m} = Y_0Y_1\ldots Y_n \prod_{0 \leq i < j \leq n} (Y_i^m-Y_j^m).\]

  It is shown in \cite{ot1992} that \(D(C_m)\) is free with a basis given by
  \[\left\{Y_0^{\ell m +1}\partialfrac{}{Y_0} +Y_1^{\ell m +
      1}\partialfrac{}{Y_1}+ \ldots Y_n^{\ell m +
      1}\partialfrac{}{Y_n} \mid \textrm{ for } 0 \leq \ell \leq n
  \right\}.\]

Let \(Z_m\) denote the dual set of points inside the dual
projective space  \(\bb{P}(V)\), and supposing \(\Sym(V^{*}) =
\bb{K}[X_0,..,X_n]\) with \(\{X_i\}\) forming a basis dual to
\(\{Y_i\}\). Then by \cref{thm:GlobalIso} we get that \(I^{\gg}(Z_m)\) is free with basis 
\[\left\{M_0^{\ell m +1}X_0 +M_1^{\ell m +
      1}X_1+ \ldots M_n^{\ell m +
      1}X_n \mid \textrm{ for } 1 \leq \ell \leq n \right\},\]
where the \(M_i\) are defined as in \cref{defn:TauIso}.

For both basis we verify that each of these elements lies in their
respective modules to help illustrate  \cref{thm:GlobalIso}. First
given \(\theta_\ell = \sum_{i=0}^n Y_i^{\ell m + 1}
\partialfrac{}{Y_i}\), we note any factor  of \(Q_{C_m}\) is of the form \(Y_i\) or \(Y_j
- \zeta^e Y_k\) where \(\zeta\) is a principle \(m\)-th root of
unity. As \(\theta_\ell(Y_i) = Y_i^{\ell m +1}\) it
follows that \(\theta_\ell \in D(\{Y_i=0\}\). Additionally,
\(\theta_\ell(Y_j - \zeta^e Y_k) = Y_j^{\ell m +1} - \zeta^e Y_k^{\ell
  m + 1}\), yet as \(\zeta^{e \ell m} = 1\) this is the same as \(\theta_\ell(Y_j - \zeta^e Y_k) =
Y_j^{\ell m +1} - (\zeta^e Y_k)^{\ell m +1}\). The identity \(a^d -
b^d = (a-b) \left(\sum_{i=0}^{d-1} a^i b^{d-1-i}\right)\) now
establishes that \(\theta_\ell \in D(C_m)\).

Similarly let \(F_\ell = \sum_{i=0}^n M_i^{m \ell + 1} X_i\) to show
\(F_\ell \in I^{\gg}(Z)\) it suffices to show that \(F_\ell\) is \(0\)
when evaluated at each point in \(Z_m\). For each \(P \in Z_m\) we let
\(\textrm{ev}_p:\bb{K}[\bm{A}] \to \bb{K}[\bm{A}]\)
denote this evaluation map. Every point in \(Z_m\) is represented by
\(E_i\) or \(E_j - \zeta^e E_k\) for \(E_i\) the standard basis
vectors of \(V\). Then \(\textrm{ev}_{E_i}(F_\ell) =
\textrm{ev}_{E_i}\left(M_i^{\ell m +1}\right)\) and \(\textrm{ev}_{E_j -
\zeta^e E_k}(F_\ell) = \textrm{ev}_{E_i}\left(M_j^{\ell m +1} -
\zeta^e M_k^{\ell m +1}\right)\). As similarly \(M_i \mid M_i^{\ell m
+1}\) and \(M_j - \zeta^e M_k \mid M_j^{\ell m + 1} - \zeta^e
M_k^{\ell m +1}\) we can conclude by \cref{lemma:LinearVanishing}.
\end{example}

The previous results about the module \(I^{\gg}(Z)\) and it's
relationship with \(D_0(\cl{A}_Z)\) can be summed up as stating
that the isomorphism of projective spaces \(\bb{P}(W) \iso
\bb{P}(\bigwedge^n V)\) extends to an isomorphism of sheaves
\(\widetilde{D_0(\cl{A}_Z)} \iso \widetilde{I^{\gg}(Z)}(-1)\). This
sheaf \(\widetilde{I^{\gg}(Z)}\) and it's relationship with
\(\widetilde{D_0(\cl{A}_Z)}\) is implicit in \cite{FV}. The
relationship of \(I^{\gg}_Q(Z)\) with a general codimension \(2\)
subspace however was only made explicit in the case \(Z
\subseteq \bb{P}^2\).

As we have committed to working algebraically we state and prove one
more result which is a simple corollary of the fact that previously
stated isomorphisms correspond to an underlying isomorphism of sheaves.

\begin{prop}\label{prop:NaturalDiagram} The following diagram commutes
  for all codimension \(2\) subspaces \(Q\), 
    \[\xymatrix{
      D_0(\cal{A}_Z) \ar[r]^{\res_{Q^{\perp}}} \ar[d]^{\eta} &
      D_0(\cal{A}_Z)\mid_{Q^\perp} \ar[d]^{\eta_Q}\\ 
      I^{\gg}(Z) \ar[r]^{\varepsilon_Q}   & I^{\gg}_Q(Z)  \\
    }\]
  with the sides isomorphisms for general \(Q\).
\end{prop}

\begin{proof} First note, that by proposition  \ref{prop:TauIso} we
  have a commuting diagram of commutative \(\bb{K}\)-algebras 
  \[\xymatrix{\Sym(W^{*}) \ar[r] \ar[d]^\tau & \Sym(W^{*})/I(Q^{\perp}) \ar[d]^{\tau_Q} \\
      \Pl(n) \ar[r]^{\varepsilon_Q} & \cl{F}_Q },\]
  where the top map sends \(f \in \Sym(W^{*})\) to its coset \(\bar{f}
  \in \Sym(W^{*})/I(Q^{\perp})\). 

Working in coordinates given \(\theta = \sum_{i=0}^n F_i
\partialfrac{}{Y_i}\), we have that
\[\begin{aligned}
  \varepsilon_Q\eta(\theta) &= \sum_{i=0}^n \varepsilon_Q(\tau(F_i))
  X_i = \sum_{i=0}^n \tau_Q(\bar{F_i}) X_i = \eta_Q(\res_{Q^{\perp}}(\theta))
\end{aligned}\]
establishing the result.
\end{proof}

\PFEIndent

\section{Unexpected Hypersurfaces}\label{sec:UnexpectedHypersurfaces}
In \cite{CHMN}, the authors gave a characerization of the degrees \(d\),
in which a finite set of points \(Z \subseteq \bb{P}^2_{\bb{C}}\) admits unexpected curves in the
specific case when \(m = d-1\). In this section we introduce the concept of very unexpected
  hypersurfaces (\cref{defn:VeryUnexpectedHypersurface}) and study them using the duality of
\cref{sec:DerivationIdealCorrespondence}. Namely in \cref{thm:CHMNGeneralization}, we achieve a
higher dimensional generalization of the main result of
\cite{CHMN} which we recall below.

\CHMNthm

The most striking part of this characterization is that it does not
depend directly on \(\dim [I(Z)]_d\) or \(\dim [I(Z) \cap I(Q)]_d\), a
feature also present in our generalization. As some papers have
already introduced a notion of unexpected hypersurface we recall this definition below, before
discussing why it is inadequate for our needs.

\begin{defn}\label{defn:OstensiblyUnexpected}
  If \(Z \subseteq \bb{P}^n=\Proj(R)\) is a set of points and \(Q\) is
some general linear subspace, we say \(Z\) admits 
  unexpected \(mQ\)-hypersurfaces in degree \(d\) if
\[\begin{aligned} \dim [I(Z) \cap I(Q)^m]_d &> \max\left\{0,\binom{n+d}{n} - \dim [R/I(Z)]_d - \dim[R/I(Q)^m]_d\right\} \\
     &> \max\left\{0,\dim[I(Q)^m]_d - \dim [R/I(Z)]_d\right\}
  \end{aligned}
\]
If \(Z \subseteq \bb{P}^2\) we instead say that \(Z\) admits
unexpected curves.
\end{defn}

If \(Q\), \(m\) or \(d\) are obvious from context, we may avoid these
qualifiers and simply specify that \(Z\) admits unexpected
hypersurfaces.

If we filter \([R/I(Q)^m]\) by \(I(Q)\) we get that each graded
component, \([I(Q)^{i-1}/I(Q)^i]\), of the corresponding filtered module
is a free module over \([R/I(Q)]\) 
generated by \([I(Q)^{i-1}/I(Q)^{i}]_{i-1}\). It follows that
\(\displaystyle \dim
[I(Q)^{i-1}/I(Q)^i]_d = \binom{\dim Q +d - i}{\dim Q} \binom{\codim
  Q -1+i}{\codim Q-1}\) and consequently we conclude that
\[\displaystyle \dim [R/I(Q)^m]_d = \sum_{i=0}^{m-1} \binom{\dim Q +d
  - i}{\dim Q} \binom{\codim Q -1+i}{\codim Q-1}.\]  This result in combination with the
Chu-Vandermonde identity, 
\(\displaystyle \sum_{j=0}^{m-b} \binom{a+j}{a}\binom{m-j}{b} =
\binom{a+m+1}{a+b+1}\),
allows us to conclude the following 
\begin{prop}\label{prop:NumericsOfUnexpectedness} \(Z \subseteq \bb{P}^n\) admits unexpected
  \(mQ\)-hypersurfaces in degree \(d\) if and only if 
\[\begin{aligned}
    \dim [I(Z) \cap I(Q)^m]_d
    & > \max\left\{0,\binom{n+d}{n} - \dim [R/I(Z)]_d
      -\sum_{i=0}^{m-1} \binom{\dim Q+d-i}{\dim Q}\binom{\codim Q
        -1+i}{\codim Q - 1}\right\} \\
    \text{ or equivalently,  }& \\
\dim[I(Z) \cap I(Q)^m]_d    & >\max\left\{0, - \dim [R/I(Z)]_d+\sum_{j=m}^d\binom{\dim Q+d-j}{\dim Q}\binom{\codim Q
        -1+j}{\codim Q - 1} \right\}
  \end{aligned}\]
In particular, if \(m = d-1\) and \(\dim Q = n-2\) the inequality becomes
\[\dim[I(Z) \cap I(Q)^{d-1}]_d > \max\{ 0, nd+1 - \dim[R/I(Z)]_d\}\]
\end{prop}

Despite the ease in which the above definition can be stated, it has a
few shortcomings. The first shortcoming is of a semantic nature, namely there are sets of points
which by definition admit unexpected 
\(mQ\)-hypersurfaces, but where we believe the difference in dimension
is unsurprising. The
second issue is somewhat larger if we hope to generalize
\cref{thm:UnexpectedCurveChar}, namely it can not be determined from
\(D_0(\cl{A}_Z)\) whether or not \(Z\) admits ostensibly unexpected
\((d-1)Q\)-hypersurfaces in degree \(d\).

Both of these issues are illustrated by the following example.

\begin{example}\label{example:TotallyExpected}
  Let \(H \subseteq \bb{P}^3\) be any plane, and let \(U\) consist of
  \(10\) of points on \(H\). Now take two general points \(P_0\) and \(P_1\)
  not on \(H\). Let \(Z = P_0+P_1+U\), and \(Q \subseteq \bb{P}^2\) a
  generic line, if we let \(\ell_H\) be a linear form defining \(H\), and
  \(\ell_0,\ell_1\) be linear forms defining \(\Spn(Q,P_0)\) and \(\Spn(Q,P_1)\)
  respectively. Then taking \(f = \ell_H\ell_0\ell_1\), we get that
  \(f\) lies in \([I(Z + 2Q)]_3\). If the points in \(U\) are general
  points on \(H\),  then \(h_U(3) = \min\left\{\binom{2+3}{3},|U|\right\}=10\), \(h_Z(3) =
  12\), and \(h_{2Q}(3) = 4+(2)(3) = 10\), in which case \(Z\) admits an
  unexpected hypersurface in degree \(3\).

  However, taking \(U'\) to be \(10\) points lying on a smooth conic
  in \(H\) and letting \(Z' = P_0+P_1+U'\), then \(h_{U'}(3) =
  7\) and \(h_Z'(3)=9\) so \(Z'\) does not admit unexpected
  hypersurfaces in degree \(3\).

  Note though that there is an isomorphism of intersection lattices \(L_{\cl{A}_Z}
  \iso L_{\cl{A}_Z'}\), and that both \(D_0(\cl{A}_Z)\) and
  \(D_0(\cl{A}_{Z'})\) have splitting type
  \((2,4,5)\).
  \end{example}

  In the above example, the ``unexpectedness'' is explained
  by the fact that most of the points of \(Z\) lie on the plane
  \(H\). This gives us a lower bound on \(\dim[I(Z+2Q)]_d\) since,
  \[\dim [I(Z+2Q)]_3 > \dim[I(H +2Q) \cap I(P_0+P_1)]_3 \geq \dim
  [I(H+2Q)]_3 - 2.\]

Furthermore, there is no reason to expect equality in the inequality \[\dim [I(H) \cap
  I(Q)^2]_3 \leq \max\left\{0,\dim [I(Q)^2]_3 - \dim [R/I(H)]_3\right\}\]
  since \(Q\) and \(H\) have nonempty intersection. 
  This situation is elaborated on further by the following proposition which
  computes the dimension of \([I^{\gg}_Q(H)]_d = [I(H) \cap I(Q)^{d-1}]_d\)
  can impose on \(I(Q)^{d-1}\).
  
\begin{prop}\label{lemma:LinearConditions} Let \(H,Q \subseteq \bb{P}(V)\) be nonempty linear subspaces,
  with \(Q\) general of codimension \(2\). Then
  \[\dim [I^{\gg}_Q(H)]_d = \dim [I(H) \cap I(Q)^{d-1}]_d = d (\codim H).\]

  As a consequence, if \(Z \subseteq H\), then \(\dim [I^{\gg}_Q(Z)]_d
  \geq d (\codim H)\).
\end{prop}

\begin{proof}   Let \(h = \dim H\). We may choose a basis
  \(\{X_0,..,X_n\}\) of \(V^{*}\) so that  \(I(H) =
  (X_{h+1},..,X_n)\).   Moreover, let \(\ell_i := \varepsilon_Q(M_i) \in
  [\cl{F}_Q]_1\), denote the linear form vanishing on
  \(Q\) and the \(i\)-th coordinate point.

  We proceed by induction on \(h\), establishing that \(I^{\gg}(H)\)
  is free with basis \(\{X_{d+1},..,X_n\}\). First consider the case
  \(h=0\), so that \(H\) is the \(0\)-th coordinate point. For each \(f \in [I^{\gg}(H)_Q]_d\), we may write \(f = \sum_{i=0}^n
  f_i X_i\) with each \(f_i \in [\cl{F}_Q]_{d-1}\). Evaluating \(f\)
  at \(H\) shows that \(f_0(H) = 0\). As \(\cl{F}_Q\) is a polynomial
  ring in two variables, we conclude that \(\ell_0\) divides \(f_0\). Using the identity \(\sum_{i=0}^n X_i\ell_i = 0\), and letting
  \(g_i = f_i - \ell_if_0/\ell_0\) we get \(f = \sum_{i=1}^n g_iX_i\). It
  follows that \(I^{\gg}(Q)\) is a free \(\cl{F}_Q\)-module with basis
  \(X_1,..,X_n\).

  Now when \(h \geq 1\), let \(H_0 \subseteq H\) be the coordinate
  subspace, with defining ideal \(I(H_0) = (X_{h},..,X_n) \supset
  I(H)\). We get by inductive hypothesis that every element
  \(f \in I^{\gg}_Q(H_0)\) may be written in the form \(f = \sum_{i=h}^n f_iX_i\).
  As \(X_j \in I(H)\) for \(j > h\), we see that \(f \in I(H)\) if and only if \(f_h \in I(H) \cap
  \cl{F}_Q\). However, as \(\bb{K}\) is infinite and \(h > 0\) we have for general \(Q\) that
  there is no finite collection of hyperplanes through \(Q\) which vanish on \(H\), and consequently
  we must have \(I(H) \cap \cl{F}_Q = 0\). Hence \(\sum_{i=h}^n f_iX_i
  \in I^{\gg}_Q(H)\) if and only if \(f_h = 0\), and so \(I^{\gg}_Q(H)\) is free with
  basis \(X_{h+1},..,X_n\) as claimed.

  Noting that \(\dim [\cl{F}_Q]_{t-1} = t\), it follows that \(\dim
  [I^{\gg}_Q(H)]_d = (\codim H)(\dim [\cl{F}_Q]_{d-1})= d(\codim H)\) as desired.
\end{proof}

\begin{example}
  In view of the preceding lemma, we see that \cref{example:TotallyExpected} can be generalized. Namely,
for \(n > 2\) we let \(H \subseteq \bb{P}^n\) be a proper linear
subspace of dimension \(d > 1\). Fix a degree \(t > 1\) and let \(Z\)
consist of \(\binom{t+d}{t}\) general points on \(H\), so that \(\dim
[R/I(Z)]_s = \min\{\binom{s+d}{s},|Z|\}\). Then the prior lemma shows
\(\dim [I^{\gg}_Q(Z)]_s = \max\{(n-d)s,ns+1-|Z|\}\), and hence that
\(Z\) admits unexpected \(Q\)-hypersurfaces in all degrees \(2
\leq s \leq t\).
\end{example}

With this discussion in mind we introduce our definition of
very unexpected hypersurface. 

\begin{defn}\label{defn:VeryUnexpectedHypersurface} Let \(Z \subseteq
  \bb{P}(V)\) be a finite set of points and \(R = \Sym(V^{*})\) the
  projective coordinate ring. For \(Q\) a generic linear subspace,  we say that \(Z\) admits
  \emph{very unexpected \(mQ\)-hypersurfaces} in degree \(d\), if there
  is a subset \(W \subseteq Z\) satisfying the following conditions:
  \begin{itemize}
  \item[(I)] \([I(Z) \cap I(Q)^{m}]_d =  [I(W) \cap I(Q)^{m}]_d\)
  \item[(II)] For all irreducible subvarieties
    \(X \subseteq \bb{P}(V)\),
    \[|W \cap X| \leq \dim [I(Q)^{m}/(I(X) \cap I(Q)^{m})]_d\]
  \item[(III)] \(W\) imposes less condition on \([I(Q)^{m}]_d\) than on
    \([R]_d\), that is
    \[\dim [R/I(W)]_d > \dim [I(Q)^{m}/(I(W) \cap  I(Q)^{m})]_d\]
  \end{itemize}

\end{defn}

\begin{remark} We note that condition \((II)\) only needs to be
  checked on \emph{positive dimensional} irreducible subvarieties.
\end{remark}

\begin{example} To illustrate the difference between
  \cref{defn:VeryUnexpectedHypersurface} and
  \cref{defn:OstensiblyUnexpected}, we revisit
  \cref{example:TotallyExpected}. Letting \(Z\) and \(Z'\) denote the
  points sets from that example, we show neither set admits very
  unexpected hypersurfaces where \(m=d-1\).

  Note that \(Z\) does not admit very unexpected
  hypersurfaces in particular because it does not admit unexpected
  hypersurfaces. Namely for any subset \(W \subseteq Z\), we have
  \[\dim [I(Q)^{d-1}/(I(W) \cap  I(Q)^{d-1})]_d \geq\dim [R/I(Z)]_d
    \geq  \dim [R/I(W)]_d,\]
  and so \(W\) does not satisfy \((III)\).

  For \(Z'\), if \(W \subseteq Z'\) is a subset satisfying \((II)\)
  for some \(d > 0\) then \(|W \cap H| \leq \dim [I(Q)^{d-1}/(I(H) \cap I(Q)^{d-1})]_d
  = 2 d +1\) in particular if \(d =3\) we have \(|W \cap H| \leq
  7\). If \(W\) also satisfies \((I)\) then \(\dim [I^{\gg}_Q(W)]_3 =
  1\) and so in particular \(\dim [I(Q)^{m}/(I(W) \cap  I(Q)^{m})]_d =
  9\)  as \(\dim [R/I(W)]_d \leq |W| \leq |W \cap H| + 2 = 9\) we see
  any such \(W\) cannot satisfy condition \((III)\).
 \end{example}

\begin{remark} It's possible that there are other definitions that are
  preferable in some ways. One change that might be useful is to
  require condition (ii) in the case where \(X\) is not necessarily
  irreducible, or if we allow \(Z\) to be nonreduced perhaps take \(X\) to be a
  positive dimensional subscheme. We use the above definition for now as it
  is strong enough for our purposes while still being relatively
  easy to check.

  In this paper we will be focusing on the case where \(\codim Q = 2\)
  and \(m = d-1\). We introduce this definition in general because we
  think it is a natural and potentially useful modification given our
  discussion in \cref{example:TotallyExpected}.
\end{remark}

\begin{remark}
Despite the fact the above definition is strictly stronger than
\cref{defn:OstensiblyUnexpected}, the two definitions agree in
\(\bb{P}^2\). This is a consequence of the fact that the only positive dimensional
subvarieties that are needed to check in condition \((ii)\) are
hypersurfaces. More generally, if \(Z \subseteq \bb{P}^n\) is a finite
set of points contained in a hypersurface defined by \((f=0)\) and  \(Q \in
\bb{P}^n\) is the generic point. Then applying the dimension count from
\ref{prop:NumericsOfUnexpectedness}, that 
\[\begin{aligned}
    \dim [I(Q)^m/(f) \cap I(Q)^m]_d &= \dim [I(Q)^m/(f I(Q)^m)]_d =
    \dim[I(Q)^m]_d - \dim  [I(Q)^m]_{d-\deg f} \\
    &= \max\left\{0,\binom{n+d}{n} - \binom{n+d-f}{n}\right\}
\end{aligned}\]
It follows for all \(m\) and \(d\) that 
\[\dim [R/I(Z)]_d \leq \dim [R/(f)]_d = \binom{n+d}{n} -
  \binom{n+d-\deg(f)}{n} \leq \dim [I(Q)^m/((f) \cap I(Q)^m)]_d\]
Establishing that condition \((iii)\) could never be satisfied under
these conditions.

A similar argument shows that \(Z\) can never admit very unexpected
\(mQ\)-hypersurfaces if \(Q\) is a hyperplane.
\end{remark}

One potential issue with \ref{defn:VeryUnexpectedHypersurface} is that condition \((II)\) seems
difficult to verify, given that naively there is a
potentially infinite number of irreducible varieties we must check. However, we make a few observations
showing that it is easier to verify than it may seem, and can be
reduced to a finite number of subvarieties.

Suppose that \(Z \subseteq \bb{P}^n\) admits unexpected
\(mQ\)-hypersurfaces in degree \(d\) and furthermore, that there's no \(P \in
Z\) where \( [I(Z) \cap
I(Q)^m]_d \subsetneq [I(Z - P) \cap I(Q)^m]_d\). This is a relatively
harmless assumption since if such a \(P\) does exist, then \(Z \setminus P\) still
admits unexpected hypersurfaces in degree \(d\).

Now if there is some positive dimensional variety \(X_1 \subseteq \bb{P}^n\) so that \(|Z \cap X_1| > \dim
[I(Q)^m/I(X_1) \cap I(Q)^m]_d\). Then \(Z \cap X_1\) imposes
less than \(|Z \cap X_1|\) conditions on \(I(Q)^m\) and so we may can find a
subset \(U_1 \subseteq Z \cap X_1\) with \(|U_1| = \dim
[I(Q)^m/(I(X_1) \cap I(Q)^m)]_d\) and \([I(Q)^{m} \cap I(U_1)]_d = [I(Q)^{m}
\cap I(X \cap Z)]_d\). Setting \(Z_1 = (Z \setminus X) \cup U_1\) we make
two observations both of which follow readily:
\begin{enumerate}[(A)]
\item \(\displaystyle [I(Z) \cap I(Q)^m]_d = [I(Z \setminus X) \cap
  I(X \cap Z) \cap I(Q)^m]_d = [I(Z_1) \cap I(Q)^m]_d\)
\item If there's a strict containment \(\displaystyle [I(X_1) \cap I(Q)^m]_d \subsetneq [I(U_1) \cap I(Q)]_d\),
  then \(Z_1\) admits unexpected hypersurfaces if and only if
  \(Z\) does.
\end{enumerate}
We may continue in this way stopping when we find a subset \(Z_k \subseteq Z_{k-1}
\subseteq ... \subseteq  Z\), where either
\begin{enumerate}[1.]
  \item \(Z_k\) does not admit unexpected hypersurfaces; or 
\item \(W = Z_k\) satisfies the conditions \((I),(II)\) and \((III)\) of
\cref{defn:VeryUnexpectedHypersurface}.
\end{enumerate}

If \(Z_k\) does not admit unexpected hypersurfaces then by
observation \((B)\), we must have \([I(X_k) \cap I(Q)^m]_d = [I(U_k) \cap I(Q)^m]_d\). Then
\[[I(Z) \cap I(Q)^m]_d \subseteq [I(U_k) \cap I(Q)^m]_d = [I(X_k) \cap
  I(Q)^m]_d.\]
Hence, the polynomials in \([I(Z) \cap I(Q)^m]\) vanish on the positive dimensional variety \(X_k\).

From the preceding discussion we can conclude the following proposition.

\begin{prop}\label{prop:GoodDefnUHyp} Let \(Z \subseteq \bb{P}^n\) and suppose \(Z\) admits
  unexpected \(mQ\)-hypersurfaces in degree \(d\). Then
  there exists \(W \subseteq Z\), so that \(W\) satisfies conditions
  \(I\) and \(II\) of \cref{defn:VeryUnexpectedHypersurface} and \(Z\)
  admits very unexpected hypersurfaces if and only if \(W\) admits
  unexpected hypersurfaces.
\end{prop}

With this discussion in mind we introduce the following definition.

\begin{defn}\label{defn:BaseLocus}
  Fix positive integers \(m,n,c\) and \(d\). 
  If \(Z \subseteq \bb{P}^n\) is a finite set of points we
  set
  \[\Bloc_d(Z,m,c) := \bigcap_{Q} V([I(Z) \cap
    I(Q)^{m}]_d).\]
  Where \(Q\) is over all linear subspaces of dimension
  \(c\). Moreover, we set \(\Bloc_d(Z) := \Bloc_d(Z,d-1,n-2)\) as this is the case
  we will focus on.
  
  If \(m=d-1\) and \(c=n-2\), we also define \(\Bloc_d(M)\) for a
  submodule \(M \subseteq \mathfrak{m}^{\gg}\) via
  \[\Bloc_d(M) = \bigcap_{F \in [M]_d}\bigcap_{Q \in \Gr(n-2,n)}
    V(\varepsilon_Q(F)).\]
  That is \(\Bloc(F_\sigma)\) is the intersection of all the
  hypersurfaces defined by \(\varepsilon_Q(F_\delta)\) as \(Q\)
  varies.
\end{defn}

From the discussion proceeding this definition, we may conclude the
following

\begin{prop}\label{prop:BetterDefUHyp} Fix \(m,n,c\) and \(d\) as
  above.  For \(Z \subseteq \bb{P}^n = \proj{R}\), and \(Q\) the generic
  \(c\)-dimensional linear subspace, we have \(Z\) admits very unexpected
  \(mQ\)-hypersurfaces if and only if there's a subset \(W \subseteq Z\) satisfying

  \begin{enumerate}
  \item[(I)] \([I(Z) \cap I(Q)^{m}]_d =  [I(W) \cap I(Q)^{m}]_d\)
  \item[(II')] For all irreducible components \(X\) of \(\Bloc_d(Z,m,c)\)
    \[|W \cap X| \leq \dim [I(Q)^{m}/(I(X) \cap I(Q)^{m})]_d\]
  \item[(III)] \(W\) imposes less condition on \([I(Q)^{m}]_d\) than on
    \([R]_d\), that is
    \[\dim [R/I(W)]_d > \dim [I(Q)^{m}/(I(W) \cap  I(Q)^{m})]_d\]
  \end{enumerate}

  Consequently, if \(\dim \Bloc_d(Z,m,c) = 0\) then \(Z\) admits
  very unexpected hypersurfaces if and only if \(Z\) admits unexpected hypersurfaces.
\end{prop}

\begin{remark} Note that \(W \subseteq Z\) may satisfy \((II')\)
  without satisfying \((II)\). For instance the points \(Z = C_{5}\) dual to the
  Ceva Arrangement \(\cl{A}_{C_{5}} \subseteq \bb{P}^2_{\bb{C}}\)
  consist of \(18\) points which admit unexpected curves in all degrees
  \(d\) with \(6 < d < 11\). Taking \(d = 7\) we note that \(W = Z\)
  does not satisfy condition \((II)\)  of
  \cref{defn:VeryUnexpectedHypersurface}, since taking \(X =
  \bb{P}^2\) we see that \(|W \cap X| = 18 > 15 = \dim [ I(Q)^6/(0)]_7\). 

  More generally, if \(H\) is a \(2\)-dimensional linear subspace in \(\bb{P}^3\) then taking
  \(Z \subseteq H\) it follows from \cref{prop:UnexForSumOfPoints}
  that \(W = Z\) does not satisfy condition \((II)\) in degree \(d =
  7\). However, in either case \(W = Z\) satisfies condition \((II')\) above.
\end{remark}

\begin{example} It should be noted here that \(\Bloc_d(Z)\) and
  \(\Bloc_d(I^{\gg}(Z))\) are not necessarily the same. For instance, if
  \(Z \subseteq \bb{P}^2_\bb{C}\) is \(5\) general points, then a
  computation shows that 
  \([I^{\gg}(Z)]_3 = 0\), and so \(\Bloc_3(I^{\gg}(Z)) = \bb{P}^2\).
  Yet a direct computation shows that \(\dim
  [I(Z) \cap I(Q)^{2}]_3 = 2\) and 
  \(\Bloc_3(Z) = Z\). It is true, however, that \(\Bloc_d(Z) \subseteq \Bloc_d(I^{\gg}(Z))\).
\end{example}

\begin{remark}From here on we restrict the view of the paper, to the case where
\(c=n-2\) and \(m = d-1\) that is we study \([I(Z) \cap
I(Q)^{d-1}]_d\).
\end{remark}

The following proposition provides a classification of those varieties
that can appear in \(\Bloc_d(Z)\).

\begin{prop}\label{prop:BlocIsLinear}
  For any submodule \(M \subseteq \mathfrak{m}^{\gg}\), (resp. \(Z
  \subseteq \bb{P}^n\)) the base locus
  \(\Bloc_d(M)\) (resp. \(\Bloc_d(Z)\)) is a union of linear subspaces.
\end{prop}

\begin{proof} We prove both statements in parallel, let \(B = \Bloc_d(M)\) or \(B =
  \Bloc_d(Z)\).

  Let \(C\) be a positive dimensional irreducible
  subvariety which is contained in \(B\) and not a linear subspace.  We establish that \(\Spn(C)
  \subseteq B\) from which the result follows.

  First we show for a general hyperplane \(H\),
  that \(\Spn(H \cap C) = H \cap \Spn(C)\). Note that \(\Spn(H \cap C)
  \subseteq H \cap \Spn(C)\) and so it suffices to show they have the
  same dimension. To do this take \(c_1,..,c_t
  \in C\) to be \(t = \dim
  \Spn(C)\) linearly independent points, and let \(L\) be any hyperplane containing
  \(\Spn(c_1,..,c_t)\), but with \(C \not\subseteq L\). Then
  \[\dim  \Spn(L \cap C) = \dim \Spn(C)-1 = \dim(L \cap \Spn(C)).\]
  It now follows that
  \(\dim \Spn(H \cap C) \geq -1 + \dim \Spn(C)\) for a general
  hyperplane \(H \subseteq \Spn(C)\), since among hyperplanes \(H\)
  which properly intersect \(C\), the quantity \(\dim \Spn(H \cap C)\) is lower
  semi-continuous. So in particular,
  \[\dim \Spn(C)-1 = \dim \Spn(L \cap C) \leq \dim \Spn(H \cap C)
    \leq \dim H \cap \Spn(C) = \dim \Spn(C)-1.\]
  Thus establishing the claim.

  Proceeding let \(Q \subseteq \bb{P}^n\) be a general codimension \(2\) subspace
  , and let \(\ell\) a general linear form vanishing on \(Q\). 
  As \(Q\) is a hypersurface considered as a subvariety of \((\ell =
  0)\), we get for any \(f \in \varepsilon_Q([M]_d)\) (resp. any \(f \in
  [I_Q^{\gg}(Z)]_d\)) that there exist linear forms \(r \in [R]_1\) and \(\ell_Q \in I(Q)\) so that
  \[f = (\ell_Q)^{d-1} r \mod (\ell).\]
  Note that as \(\ell\) is general, we may assume that \(f \neq 0 \mod (\ell)\).
  Since \(Q\) is general we can assume that  for every positive
  dimensional component \(C\) of \(B\), that \(C \not\subseteq Q\) and
  furthermore that \(Q\) contains no component of \(C \cap (\ell = 0)\).
  As \(r\) is linear, it vanishes on \(\Spn((\ell = 0) \cap B) = \Spn(B) \cap (\ell = 0)\). It
  follows that for any component \(C\) of \(B\) that \(f\) vanishes on a general hyperplane section of
  \(\Spn(C)\). As \(\Spn(C)\) is irreducible we conclude that if
  \(\dim(C) > 0\) then  \(f\) vanishes on \(\Spn(C)\) as desired.
\end{proof}

\begin{remark}
  Given \(B \subseteq \bb{P}^n\) as in the proof above, we note the proof shows that in
  fact for a general hyperplane \(H\), that  \(\Spn(B \cap H) = n-2\).
  This places further constraints on \(B\). With some more work it can be shown
  that \(B = \Bloc_d(I^{\gg}(Z))\) for some \(Z\) and some \(d\) if
  and only if \(B\) is a union of linear subspaces \(H_0,..,H_s\) which
  satisfy 
  \[\sum_{i \in J} \dim H_i < \dim \Spn\left(\bigcup_{i \in J}
      H_i\right) \;\textrm{ for all } \; J \subseteq \{0,..,s\},
    \textrm{ with } |J| > 2.\]

  In particular it follows the \(H_i\) are disjoint.
\end{remark}

Combining the 2 preceding proposition with
\cref{prop:NumericsOfUnexpectedness}, it follows that conditions
\((I)\) and \((II)\) of \cref{prop:BetterDefUHyp} can be checked by looking at the
combinatorics of linear subspaces spanned by subsets of \(Z\). With
this in mind we state one of the main theorems of this section though we
postpone the proof until after \cref{thm:SegreBound}.

\begin{defn}\label{defn:BestDefUHyp} Given a finite set of points \(Z
  \subseteq \bb{P}^n\) and a real number \(d \in \bb{R}\) we define
  the \emph{modified expected number of conditions}, as the integer
   \(\Exc(Z,d)\), which is the solution to
  optimization problem
  \[\Exc(Z,d) = \min\left\{\sum_{i=0}^s (d\dim(H_i) + 1) \mid
  \{H_0,...,H_s\} \text{ are nonempty linear subspaces with } Z \subseteq
  \bigcup_{i=0}^s H_i \right\}.\]

\end{defn}

\begin{theorem}\label{thm:BestDefUHyp}
 Let \(Z \subseteq \bb{P}(V)\) be a finite set of points.
Then \(Z\) admits unexpected \((d-1)Q\)-hypersurfaces in degree \(d\)
if and only if
\[\dim [I(Q)^{d-1}]_d - \dim [I(Q)^{d-1} \cap I(Z)]_d < \Exc(Z,d)\]
or in the notation of \cref{sec:DerivationIdealCorrespondence},
\[\dim [I(Q)^{d-1}]_d - \dim [I^{\gg}_Q(Z)]_d = \dim[I(Q)^{d-1}
  /I^{\gg}_Q(Z)]_d< \Exc(Z,d)\]

\end{theorem}

It turns out that the linear program defined in \cref{defn:BestDefUHyp}
can be studied via combinatorial objects known as Matroids. We recall
the definition of a Matroid
now for convenience.

\begin{defn}\label{defn:Matroid} A \emph{matroid} is a finite set \(M\) along with a rank
  function \(\rk_M:2^{M} \to \bb{Z}\), which satisfies the following
  \(3\) conditions

  \begin{enumerate}[(\bf RK 1)]
  \item \(0 \leq \rk_M(A) \leq |A|\)
  \item If \(A \supseteq B\), then \(\rk_M(A) \leq \rk_M(B)\)
  \item \(\rk_M(A) + \rk_M(B) \geq \rk_M(A \cup B) + \rk_M(A \cap B)\)
  \end{enumerate}
  We note a function satisfying only {\bf (RK 3)} is a
  \emph{submodular function}.

  Equivalently, it may be defined as a nonempty collection of subsets
  \(\cl{I}\) of \(M\), which satisfy
  \begin{enumerate}[\bf ({IND} 1)]
  \item If \(A \in \cl{I}\) and \(B \subseteq A\) then \(B \in \cl{I}\)
  \item If \(A,B \in \cl{I}\) and \(|A| < |B|\), then there exists
    some \(b \in B\) so that \(A \cup \{b\} \in \cl{I}\).
  \end{enumerate}
\end{defn}

We now recall a few more pieces of related terminology. We refer to
\cite{Oxley} for definitions.

\begin{itemize}
\item  A subset \(I \subseteq M\) is \emph{Independent} if and only if \(|I| =
  \rk_M(I)\). Conversely for \(A \subseteq M\), \(\rk_M(A)\) is equal
  to the largest size of an independent \(I \subseteq A\).
\item A maximal independent set is a \emph{Basis} of \(M\). Every
  basis has the same size namely \(\rk_M(M)\), and every independent
  subset is contained in some basis.  
\item A \emph{flat} of rank \(r\) is a subset \(F \subseteq M\), which is maximal
  among subsets of \(M\) with rank \(r\). Every subset \(A\) of \(M\) is
  contained in a unique flat, \(F\), with \(\rk_M(A) = \rk_M(F)\),
  this flat \(F\) is often denoted \(\Cl_M(A)\).
\item Let \(f:M \to \bb{R}\) be any increasing submodular function (meaning it
  satisfies only conditions \(ii\) and \(iii\) of \cref{defn:Matroid}) defines a matroid
  \(M_f\). Its independent sets are those \(I \subseteq M\), which
  satisfy \(|A| < f(A)\) for all nonempty \(A \subseteq I\). 
\end{itemize}

\begin{example} Every finite set of points \(Z \subseteq \bb{P}(V)\)
  defines a matroid, \(M(Z)\) namely for every nonempty \(A \subseteq Z\), we set
  \[\rk_{M(Z)}(A) = 1 + \dim \Spn(A)\]
  A flat of this matroid, is the intersection of \(Z\) with a linear
  subspace \(L \subseteq \bb{P}^n\).

  Matroids of this type are referred to as \emph{representable
    matroids}, and are in some sense the prototypical example of a matroid.
\end{example}

One result which we will need later on relates information about the matroid
determined by \(Z\) to information about the ideal \(I(Z)\). The
result below is shown in section 4 of \cite{NT}, where a further
generalization is given to schemes of fat points. 

\begin{theorem}[\cite{NT}]\label{thm:SegreBound} Let \(Z \subseteq
  \bb{P}^n\) be a nonempty finite set
  of points satisfying
  \[|Z \cap L| \leq d(\dim L) + 1\]
  for all nonempty linear subspaces \(L \subseteq
  \bb{P}^{n}\). Then \(Z\) imposes independent conditions on \(d\)
  forms meaning
  \[\dim [\Sym(V^{*})/I(Z)]_d = |Z|\]
\end{theorem}

We are now ready to state and prove the main result of this
section. 

\begin{proof}[Proof of \Cref{thm:BestDefUHyp}]
 Fix a integer \(d\), note that it follows from the definition that 
  \[\Exc(Z,d) =  \min\left\{\sum_{i=0}^s (d\dim \Spn(A_i)
      +1 ) \mid
  \{A_0,...,A_s\} \text{ are nonempty subsets of \(Z\) with } Z =
  \bigcup_{i=0}^s A_i \right\}.\]

Before proving either direction of the equivalence. We establish the
claim below.

\begin{claim}
\(\Exc(Z,d)\) is equal to the largest size of a subset \(B \subseteq
Z\) which satisfies the following \(3\) conditions
\begin{enumerate}[(C1)]
\item \([I_Q^{\gg}(B)]_d = [I^{\gg}_Q(Z)]_d\)
\item For all linear subspaces \(L\), \(|B \cap L | \leq
  \dim[I(Q)^{d-1}/I^{\gg}_Q(L)]_d = d(\dim B) + 1\)
\item \(B\) imposes independent conditions on \(d\) forms.
\end{enumerate}
\end{claim}

\begin{proof}[proof of claim]
Applying results from \cite{Eureka} (see theorem (8) and comment
(16)), we may define a matroid \(M_d\) on the set \(Z\) whose independent sets
are precisely those \(I \subseteq Z\) where \(|A| \leq d \dim (\Spn A)  - 1\) for all nonempty \(A \subseteq I\). The linear
programming duality given in \cite{Eureka}, now states that
\[\rk(M_d) =\Exc(Z,d).\]
From this we can conclude that \(\Exc(Z,d)\) is equal to the largest
size of a subset which satisfies condition \((C2)\), namely any basis
of \(M_d\) works. To finish the proof of the claim we find a basis of
\(M_d\) satisfying \((C1)\) and \((C3)\).

By \cref{prop:GoodDefnUHyp}, there is some \(W \subseteq Z\) so that \(W\) satisfies conditions
\((C1)\) and \((C2)\). As \(W\) satisfies \((C2)\) it is independent in
\(M_d\) and we can therefore
extend it to a basis \(W \subseteq B\) of \(M_d\). Now as \(W \subseteq B \subseteq
Z\) we have that \([I^{\gg}(B)]_d = [I^{\gg}(Z)]_d\), and therefore
\(B\) satisfies \((C1)\). 

Lastly, we note that \cref{thm:SegreBound} ensures that \(B\) since
\(B\) satisfies \((C2)\) it necessarily imposes
independent conditions on \(d\) forms, thereby establishing condition
\((C3)\) and the claim.
\end{proof}

Now continuing with the proof of the equivalence. If
\(\dim [I(Q)^{d-1}]_d - \dim [I^{\gg}_Q(Z)]_d < \Exc(Z,d)\)
we can find some \(B \subseteq Z\) so that \(|B| = \Exc(Z,d)\) and
\(B\) satisfies conditions \((C1)\), \((C2)\) and \((C3)\).
Then we have 
\[\dim [I(Q)^{d-1}]_d - \dim [I(Q)^{d-1} \cap I(Z)]_d < \Exc(Z,d) =
  |B| = \dim [\Sym(V^{*})/I(B)]_d.\]
Letting \(W = B\), we see that \(W\) satisfies the necessary criteria of \cref{defn:VeryUnexpectedHypersurface},
and so \(Z\) admits very unexpected hypersurfaces.

Conversely, suppose that \(Z\) admits very unexpected
hypersurfaces. Then by definition
there exists \(U \subseteq Z\) so that for general \(Q\),
\begin{enumerate}[(I)]
\item \([I_Q^{\gg}(U)]_d = [I^{\gg}_Q(Z)]_d\)
\item For all linear subspaces \(L\), we have \(|U \cap L | \leq
  \dim[I(Q)^{d-1}/I^{\gg}_Q(L)]_d = d(\dim L) + 1\)
\item \(\dim [R/I(U)]_d > \dim[ I(Q)^{d-1}/I^{\gg}_Q(U)]_d\)
\end{enumerate}

Finding a subset \(W \subseteq U\) so that \([I(U)]_d = [I(W)]_d\) and
\(W\) imposes independent conditions on \(d\) forms. We get by the
claim above that \(|W| \leq \Exc(Z,d)\) and so
\[\dim[ I(Q)^{d-1}/I^{\gg}_Q(U)]_d < \dim [R/I(U)]_d = |W| <
  \Exc(Z,d).\]

\end{proof}

\begin{remark} Let \(L \subseteq \bb{P}^n\) be a nonempty linear
  subspace. We note that the above proof relies on a somewhat
  remarkable agreement between the dimension \(\dim [I(Q)^{d-1}/(I(L) \cap
  I(Q)^{d-1})]_d\) and the quantity \(d \dim L +
  1\) appearing in the inequality from \cref{thm:SegreBound}. This
  is even more remarkable considering that the proof of
  \cref{thm:SegreBound} is almost entirely combinatorial relying on a
  generalization of Edmonds Matroid Partition Theorem.
\end{remark}

Combining the above result with \cref{thm:restrictedModIso}, we obtain
the following as a corollary.

\BestDefUHyp

\begin{remark}
   Note one consequence of this is if \(Z\) admits very unexpected hypersurfaces in degree
  \(d\), then \(a_1 < d < a_n\).
\end{remark}

\begin{proof}
Let \(H_1,..,H_s\) be any collection of linear subspaces covering
\(Z\). Note that \(I^{\gg}_Q(Z) \supseteq  \bigcap_{i=1}^s
I^{\gg}_Q(H_i)\)
and that \(\bigcap_{i=1}^s I^{\gg}_Q(H_i)\)
is the kernel of the canonical map \([I(Q)^{d-1}]_d \to
\oplus_{i=1}^s [I(Q)^{d-1}/I^{\gg}_Q(H_i)]_d\). We have by dimension counting
that for a fixed \(d\)
\[\dim [I^{\gg}_Q(Z)]_d \geq \dim \bigcap_{i=1}^s
[I^{\gg}_Q(H_i)]_d \geq nd+1 - \left(\sum_{i=1}^s d \dim(H_i) + 1\right).\]
  Taking \(H_1,..,H_s\) so \(\sum_{i=1}^s d \dim(H_i) + 1 =
  \Exc(Z,d)\), the rest follows directly from theorem \ref{thm:BestDefUHyp}
  and corollary \ref{cor:DimensionFormula}.

  The final consequence follows since if \(d \leq a_1\) then
  \(I^{\gg}_Q(Z) = 0\), if \(d \geq a_n\) then note that \(\Exc(Z,d)
  \leq |Z|\), and so
  \[nd-(|Z|-1) = \sum_{i=1}^n \max\{0,d-a_i\}
      \geq nd+1 -\Exc(Z,d) \geq nd+1 - |Z|\]

  Establishing that \(\Exc(Z,d) = |Z|\) and that the middle inequality
  is an equality.
\end{proof}

The following lemma, shows that the inequality in the preceding
corollary above may be replaced by
\[\sum_{i=1}^n \max\{0,a_i-d\} \geq  |Z| - \Exc(Z,d) \geq 0.\]

\begin{lemma}\label{lemma:InequalityNumerology} Let \(Z \subseteq \bb{P}^n\) be a finite set of points
  and suppose that \((a_1,..,a_n)\) is the splitting type of
  \(D_0(\cl{A}_Z)\). Then for all real numbers \(c\) and \(d\) 
  \[\begin{aligned}
      &\sum_{i=1}^n\max\{0,d-a_i\} \geq& nd+1 - c\\
      \iff& \sum_{i=1}^n \max\{0,a_i-d\} \geq & |Z| -c 
    \end{aligned}
  \] 
\end{lemma}

\begin{proof} Using that \(\sum_{i=1}^n a_i = |Z|-1\) we obtain
  \[
  \begin{aligned}
    \sum_{i=1}^n\max\{0,d-a_i\} &\geq& nd+1 - c& \iff\\
    \left(\sum_{i=1}^nd-a_i\right) - \left(\sum_{j; a_j \geq d} d - a_j\right)  &\geq& nd+1 - c&\iff\\
    nd - (|Z|-1) + \sum_{j; a_j \geq d} \left(a_j-d\right) & \geq &nd+1 - c&\iff\\
    \sum_{i=1}^n \max\{0,a_i-d\} &\geq & |Z| - c &
  \end{aligned}
\]
  
\end{proof}

We now conclude this section by discussing a few conditions on \(Z\)
which makes it easier to determine if \(Z\) has very unexpected
hypersurfaces in some degree \(d\).
The first is a consequence of the preceding lemma and \cref{thm:BestDefUHyp}.

\begin{cor}\label{cor:CHMNGeneralization} Let \(Z \subseteq
  \bb{P}^n\) be a finite set of points, with \((a_1,a_2,..,a_n)\) the
  splitting type of \(D_0(\cl{A}_Z)\). Suppose we have for a fixed integer \(d \geq 0\) that
  \[\Exc(Z,d) = \min\{|Z|,nd+1\}.\]
  Then the following are equivalent:
  \begin{enumerate}[(a)]
  \item \(Z\) admits very unexpected hypersurfaces in degree \(d\)
  \item \(Z\) admits unexpected hypersurfaces in degree \(d\)
  \item \(a_1 < d < a_n\)
  \end{enumerate}
\end{cor}

\begin{proof} Note first note that by definition \(\Exc(Z,d) \leq
  \min\{nd+1,|Z|\}\), since \(\Exc(Z,d) \leq |Z|\) and \(\Exc(Z,d)
  \leq d \dim(\bb{P}^n) + 1\).

  \(\bf [(a) \iff (c)]\) First, as mentioned after
  \cref{cor:BestDefnUHyp} we have that \((a) \implies (c)\). For the reverse direction assume that \(a_1 < d < a_n\). First in the
  case that \(\Exc(Z,d) = nd+1\) we see that \(Z\) admits unexpected
  hypersurfaces in degree \(d\) as \(d > a_1\), and so the inequality
  in \cref{cor:BestDefnUHyp} is strict.
  For the case when \(\Exc(Z,d) = |Z|\), we similarly conclude by
  applying \cref{lemma:InequalityNumerology} and using that \(d < a_n\).

  \(\bf [(a) \iff (b)]\) The forward direction is by
  definition. For the reverse we use the equivalence of (a) and (c),
  and note it suffices to show that \(Z\) cannot admit unexpected
  hypersurfaces in degree \(d\) if \(d \leq a_1\) or \(d \geq
  a_n\). If \(d \leq a_1\), we note this is impossible as
  \([I^{\gg}_Q(Z)]_d = 0\). If \(d \geq a_n\), then
  \[\dim[I^{\gg}(Z)]_d = \sum_{i=1}^n \max\{0,d-a_n\} = nd -
    \sum_{i=1}^n = nd -(|Z| - 1) = nd+1 - |Z|.\]
  As \(\dim [I(Q)^{d-1}]_d = nd+1\) we conclude that \(Z\) imposes
  independent conditions on \([I(Z)]_d \), and so \(Z\) cannot admit
  unexpected hypersurfaces.
\end{proof}

In the case that the points of \(Z\) are not too concentrated on one
or more proper subspaces, it turns out that \(\Exc(Z,d) =
\max\{nd+1,|Z|\}\) holds for all \(d\) and we obtain the following result.

\CHMNGeneralization

\begin{proof}
  By \cref{cor:CHMNGeneralization}, it suffices to
  show that \(\Exc(Z,d) = \min\{nd+1,|Z|\}\). Let \(\cl{H} = \{H_1,..,H_s\}\) be a
  collection of positive dimensional linear
  subspaces, so that setting \(W = Z \setminus \bigcup_{i=1}^s H_i\)
  we have
  \[|W| + \sum_{i=1}^s d \dim(H_i) + 1 = \Exc(Z,d).\]
  As \(|W| + \sum_{i=1}^s d \dim(H_i)+1\) is at a minimum, we make the
  following observations:
  \begin{enumerate}[\bf ({Ob.} 1)]
  \item \(d \dim(H_i) + 1 \leq |H_j \cap Z|\).
  \item For all \(J \subseteq
    \cl{H}\) we have \(\sum_{H_j \in J} d \dim(H_j) + 1 \leq d
    \dim\Spn\left(\bigcup_{H_j \in J} H_j\right) + 1\).
  \item \(\sum_{i=1}^s \dim (H_i) \leq \dim\Spn\left(\bigcup_{i=1}^s
      H_i\right) < n\).
  \end{enumerate}

  {\bf (Ob. 1)} and {\bf (Ob. 2)} must hold since otherwise we could find a set of points \(W'\)
  and a collection of subspaces \(\cl{H}' = \{H'_1,...,
  H'_k\}\)  with \(Z \subseteq W'
  \cup \bigcup_{H'_i \in \cl{H}'} H'_i\) and \(\sum_{H'_i \in \cl{H}'} d \dim
  H'_i + 1 < \Exc(Z,d)\). For instance, in {\bf (Ob. 1)} we would
  consider \(W' = W \cup (Z \cap H_i)\) and \(\cl{H}' = \cl{H}
  \setminus \{H_i\}.\) {\bf (Ob. 3)} is a consequence of {\bf (Ob. 2)}.

  Note that {\bf (Ob. 1)} implies that \(\Exc(Z,d) = |Z|\) for all \(d \geq
  \frac{|Z|-1}{n} \geq \frac{|Z \cap H_i| - 1}{\dim H_i}\), so suppose that \(nd+1 < |Z|\). Let \(g_i =
  |Z \cap H_i| - \left(d \dim (H_i)+1\right) \geq 0\), and note that
  by hypothesis \(\frac{g_i}{\dim H_i} = \frac{|Z \cap H_i|
    -1}{\dim H_i} - d \leq \frac{|Z| - nd - 1}{n}\). Combining this
  with our formula for \(\Exc(Z,d)\), we obtain the following
  \[\begin{aligned}
      \Exc(Z,d)&=|W| + \sum_{i=1}^s \left(d \dim(H_i) + 1\right)\\
      &= |Z| - \sum_{i=1}^sg_i \geq |Z| - \sum_{i=1}^s (\dim H_i) \left(\frac{|Z| -nd-1}{n}\right)\\ 
      &\geq |Z| - \left(|Z| -nd-1\right)\left(\frac{\sum_{i=1}^s \dim
        H_i}{n}\right) \\
    \end{aligned}\]
  Now as \(\sum \dim H_i \leq n\) by {\bf (Ob. 3)} we obtain
  \(\Exc(Z,d) \geq |Z| - \left(|Z| -nd-1\right) = nd+1\).
  As it's always true that \(\Exc(Z,d) \leq nd+1\), the result now follows.
\end{proof}

\begin{remark}To close we  spell out the connection between
  \cref{thm:CHMNGeneralization} and the original
  \cref{thm:UnexpectedCurveChar} from \cite{CHMN}
\end{remark}

\begin{proof}[proof of Theorem \ref{thm:UnexpectedCurveChar}]  Let \(Z \subseteq \bb{P}^2\) and let
  \((a_1,a_2)\) be the splitting type of \(D_0(\cl{A}_Z)\). First
  consider the case where there exists some \(L \subseteq \bb{P}^2\) so that \(|L
  \cap Z| > a_1 + 1\). Let \(Q \in \bb{P}^2\) be a general point, and take \(f \in [I^{\gg}_Q(Z)]_{a_1+1}\) be a
  minimal generator of \(I^{\gg}(Z)\), then applying Bezout's Theorem
  we see that \(L\) must be a component of the variety \(f=0\). Hence,
  \(f\) factors as \(f = \ell g\) where \(\ell\) is the linear form
  defining \(L\) and \(g \in [\cal{F}_Q]_{a_1} =
  [I(Q)^{a_1}]_{a_1}\) is a product of linear forms. As \(f\) is a minimal generator and \(Q\) is
  general it follows each linear form in \(g\) vanishes at precisely
  one point of \(|Z \setminus L|\). Therefore, as \(a_1 = \deg g = |Z
  \setminus L|\), we conclude that \(|Z \cap L| = |Z| - |Z \setminus
  L| = a_1 + 1\).

  Now noting that if \(\Exc(Z,d) = \sum_{i=1}^k d \dim H_i +1\) for
  linear subspaces, \(H_i\) that we must have \(d \dim (H_1+H_2) + 1
  \geq d \dim H_1 + d \dim H_2 + 2\). It follows that \(\dim (H_1 +
  H_2) > \dim H_1 + \dim H_2\), in the case that \(Z \subseteq
  \bb{P}^2\), this implies that there is at most one line or plane among
  the \(H_i\). Therefore, we conclude that \(\Exc(Z,d) =
  \min\{2d+1,(d \dim L +1)+|Z\setminus L|,|Z|\}\) or equivalently
  \[
    \Exc(Z,d) = \begin{cases}
      2d+1 & \text{ If } d \leq a_1\\
      d+1+a_1 & \text{ If } a_1 \leq d \leq a_2\\
      a_1+a_2+1 & \text{ If } a_2 \leq d
    \end{cases}\]
  Applying \cref{thm:restrictedModIso} and a direct comparison now
  shows that \(Z\) admits no unexpected curves, establishing this case.

  For the other case we have \(|L \cap Z| \leq a_1+1\) for all \(L
  \subseteq \bb{P}^2\).  Then for all lines \(L \subseteq \bb{P}^2\)
  we have the inequality
  \(|Z \cap L| \leq a_1+1 \leq \left\lfloor\frac{|Z|-1}{2}\right\rfloor+1 \).
  Subtracting through by \(1\) gives 
  \[|Z \cap L| -1 \leq a_1 \leq \frac{|Z|-1}{2}\]
  allowing us to conclude by \cref{thm:CHMNGeneralization}.

\end{proof}

\subsection{Computations and Examples of Unexpected Hypersurfaces}
Combinatorial Optimization problems similar to the linear program,
\(\Exc(Z,d)\), from in \cref{defn:BestDefUHyp} have been studied before.
One notable instance of this is in the paper \cite{Nara}. In
\cite{Nara} the author fixed a submodular function \(\mu:S \to
\bb{R}\) and a real parameter \(\lambda\), and studied the optimization
problem
\[\min\left\{ \sum_{i=0}^t \mu(S_i) - \lambda \;\middle| \,\{S_0,..,S_t\} \textrm{ is
      a partition of } S\right\}.\]

It was shown in section 3 of \cite{Nara} that for a fixed \(\mu\) and \(\lambda\)
that there is a unique \emph{finest} and a unique \emph{coarsest} partition of \(S\) achieving this
minimum. Here we say a partition \(\pi\) if \emph{finer} than the
partition \(\tau\) (or equivalently that \(\tau\) is \emph{coarser}
than \(\pi\))
and write \(\pi \leq \tau\), if every block of \(\pi\) is contained in
a block of \(\tau\).
In section 4 of \cite{Nara} an algorithm was given which solves this problem for a fixed
\(\mu\). It was shown in particular that minimum is a piecewise linear function of \(\lambda\).

We note that \(\Exc(Z,d)\) is equivalent to 
\[d \min\left \{\sum_{i=0}^t \left(\rk_{M(Z)}(A_i) - \frac{d-1}{d}\right) \; \middle|
    \,\{A_0,..,A_t\} \textrm{ is a partition of } Z\right\}\]
and so the algorithm given in \cite{Nara}  can be used to solve \(\Exc(Z,d)\). 

\begin{defn}\label{defn:ExpectedBaseLocus} Let \(Z \subseteq
  \bb{P}^n\) be a finite set of points, for each \(d \geq 0\), we
  define the \emph{modified expected base locus}, which we denote \(\Exbl(Z,d)\) to be the coarsest partition in the partition
  order which satisfies
  \[\sum_{B \in \Exbl(Z,d)} \left(d \dim\Spn(B) + 1\right) =
    \Exc(Z,d).\]

  Meaning that if \(\Pi\) is any other partition with \(\sum_{P \in
    \Pi} \left(d \dim\Spn(P) + 1\right) =
    \Exc(Z,d),\) then for every \(P \in \Pi\) there is some \(B \in
    \Exbl(Z,d)\) so that \(P \subseteq B\).
\end{defn}

Section 3 of \cite{Nara} establishes not only that \(\Exbl(Z,d)\)
exists, but also that in the partition order \(\Exbl(Z,d)  \geq
\Exbl(Z,d+1)\). We now make a few
observations about \(\Exbl(Z,d)\) and \(\Exc(Z,d)\) in order to
compute \(\Exc(Z,d)\) more easily.  These results
are heavily influenced by the results and techniques in \cite{Nara}.
However, our results are stronger in some cases as we can take advantage of the fact that
\(\rk_{M(Z)}\) is the rank function of a matroid, and not merely
a submodular function.

\begin{lemma}\label{lem:BlockLemma} For any real number \(d > 0\), and
  for distinct blocks \(B_1,..,B_k \in \Exbl(Z,d)\) we have
  \[\sum_{i=1}^k (d \dim\Spn(B_i) + 1) < d
    \dim\Spn\left(\bigcup_{i=1}^k B_i\right)+1\]
  In particular, for each pair of distinct blocks \(B_1\) and \(B_2\),
  \(\Spn(B_1)\) and \(\Spn(B_2)\) are disjoint subspaces.

  Similarly, if \(C_1 \sqcup ... \sqcup C_\ell\) is a partition of a
  block \(B_\ell \in \Exbl(Z,d)\) into nonempty subsets, then
  \[\sum_{j=1}^\ell (d \dim\Spn(C_j) + 1) \geq d
    \dim\Spn(B)+1\]
\end{lemma}

\begin{proof}
  If \(\Exbl(Z,d) = \{B_1,..,B_m\}\) then let
   \(A = \left\{\bigcup_{i=1}^k B_i,B_{k+1},...,B_m\right\}\).
  As \(A\) is coarser than \(\Exbl(Z,d)\), we get  \[\sum_{a \in A} d \dim \Spn(a) + 1 > \sum_{b \in \Exbl(Z,d)} d \dim
    \Spn(b) + 1.\]
  Subtracting away the shared terms now gives the desired inequality.

  The proof of the second claim follows similarly, since we must have
  \[ \left(\sum_{j=1}^\ell d \dim\Spn(C_j) + 1\middle) +
      \middle(\sum_{B \in \Exbl(Z,d); B \neq B_\ell} d \dim\Spn(B) +1\right) \geq \sum_{B \in \Exbl(Z,d)} d \dim\Spn(B) +1\]
\end{proof}

It can be somewhat laborious to determine if a given set of points
satisfies the combinatorial condition in \cref{thm:CHMNGeneralization}.
Furthermore, most of the observed configurations of points \(Z\) which admit
unexpected curves possess certain kinds of symmetry, namely their dual
arrangements \(\cl{A}_Z\) are reflection arrangements. We designed
this next proposition with these examples in mind.

\begin{defn} A \emph{psuedoreflection} is a matrix \(R \in
  \bb{GL}(n,\bb{K})\) so that \(R^k = I_n\) for some \(k > 1\) and the
  set of points in \(\bb{K}^n\), which are fixed by \(R\), denoted
  \(\textrm{Fix}_R\), form a hyperplane.
  A \emph{reflection group} is a
  subgroup, \(G\), of \(\bb{GL}(n,\bb{K})\), which is generated by
  \emph{psuedoreflections}. \(G\) is an \emph{irreducible reflection
    group} if there no nontrivial \(G\)-invariant subspace of \(\bb{K}^n\).
\end{defn}

\begin{prop}\label{prop:ReflectIsGood} If \(Z \subseteq
  \bb{P}^n_{\bb{K}}\) is a finite set of points, and there is an
  irreducible reflection group \(G \subseteq \bb{PGL}(\bb{K},n)\)
  acting on \(Z\). Then for all positive dimensional linear subspaces
  \(H \subseteq \bb{P}^n\) we have 
  \[\frac{|Z \cap H| -1 }{\dim H} \leq \frac{|Z| - 1}{n}.\]
  Consequently by \cref{cor:CHMNGeneralization}, \(Z\) admits very
  unexpected hypersurfaces in degree \(d\) for precisely those \(d\)
  with \(a_1 < d < a_n\), where \((a_1,..,a_n)\) is the splitting type
  of \(D_0(\cl{A}_Z)\).
\end{prop}

 We first note a useful criterion, which is used in the proof of the
 above proposition.

\begin{claim} Let \(Z \subseteq \bb{P}^n\), then the following are
  equivalent:
  \begin{enumerate}
  \item For all positive dimensional subspaces \(H \subseteq \bb{P}^n\),
    \[\displaystyle \frac{|Z \cap H| -1 }{\dim H} \leq \frac{|Z| - 1}{n}.\]
  \item \(\Exc(Z,q) = \min\{qn+1,|Z|\}\) for all \(q \in \bb{Q}\).
  \end{enumerate}
\end{claim}

\begin{proof}[proof of claim]
  The forward direction is established in
  \cref{thm:CHMNGeneralization}. For the reverse direction, we prove
  the contrapositive. Namely, suppose that there is some \(H \subseteq
  \bb{P}^n\) with \(\frac{|Z \cap H|-1}{\dim H} >
  \frac{|Z|-1}{n}\). Then choose any \(q\) with
  \[\frac{|Z \cap H|-1}{\dim H} > q >
    \frac{|Z|-1}{n}.\]
  Note then that \(|Z \cap H| > q \dim H + 1\) and that \(q n + 1 >
  |Z|\), hence we have that
  \[\Exc(Z,d) \leq q \dim H + 1 + |Z \setminus H| < |Z| < qn+1\]
  establishing the result.
\end{proof}
  
\begin{proof}[Proof of \cref{prop:ReflectIsGood}]

  First, note that if \(G\) is any group acting on \(Z\) then this
  action extends to the lattice of partitions of \(Z\). Furthermore,
  if \(\Pi\) is any partition of \(Z\), then for any \(g \in G\) we have \(\sum_{P \in \Pi} d
  \dim \Spn(P) + 1 = \sum_{P \in \Pi} d
  \dim \Spn(gP) + 1\). From this it follows that \(\Exbl(Z,d)\) is
  fixed by the \(G\) action, in the sense that blocks of
  \(\Exbl(Z,d)\) are taken to other blocks of \(\Exbl(Z,d)\).

  Now we continue to establishing the proposition. By the preceding
  claim it suffices to show that for rational \(q\), \(\Exbl(Z,q)\) is either the
  discrete or the indiscrete partition. Suppose that \(B \in
  \Exbl(Z,d)\) is a block with \(|B| \geq 2\), let \(r \in G\) be a
  psuedoreflection and \(H_r = \textrm{Fix}_r\) the hyperplane of the points
  fixed by \(r\). As  \(\dim \Spn(B) \geq 1\) then consequently \(H_r
  \cap \Spn(B)\) and hence \(\Spn(r B) \cap \Spn(B)\) are both
  nonempty. Applying \cref{lem:BlockLemma}, we see that we must have
  \(rB = B\) and so \(r \Spn(B) = \Spn(B)\). Therefore, \(\Spn(B)\) is
  a nonzero \(G\)-invariant subspace of \(\bb{P}^n\). As \(G\) is an
  irreducible reflection group we must have that \(\Spn(B) =
  \bb{P}^n\) and so \(B = Z\) by \cref{lem:BlockLemma}. 
\end{proof}

For a set of points \(Z \subseteq \bb{P}^n\), if the splitting type of \(D_0(\cl{A}_Z)\) is
known, then determining when \(Z\) admits very unexpected
hypersurfaces comes down to computing \(\Exc(Z,d)\). The following two
propositions can be useful in determining \(\Exc(Z,d)\). The first
places bounds on how \(\Exc(Z,d)\) can change between degrees.

\begin{lemma} For \(Z \subseteq \bb{P}(V)\),  the sequence of forward differences \(\delta_d =
  \Exc(Z,d+1) - \Exc(Z,d)\) is nonincreasing.
  Furthermore, we have
  \[ \sum_{A \in \Exbl(Z,d)} \dim(A) \geq \delta_d \geq \sum_{B \in \Exbl(Z,d+1)} \dim(B).\]
\end{lemma}

\begin{proof} 
Consider the following inequalities
  \[\begin{aligned}
      \sum_{A \in \Exbl(Z,d)} \dim \Spn(A) &= \sum_{A \in \Exbl(Z,d)} \left[ (d+1)
        \dim\Spn(A) + 1\right] - \sum_{A \in \Exbl(Z,d)}\left[d \dim \Spn(A) + 1\right] \\
      &\geq \sum_{A \in \Exbl(Z,d)} \left[ (d+1)\dim\Spn(A) + 1\right]
      - \sum_{B \in \Exbl(Z,d+1)} \left[d \dim\Spn(B) + 1\right]\\
      &\geq \sum_{B \in \Exbl(Z,d+1)} \left[ (d+1)\dim\Spn(B) + 1\right]
      - \sum_{B \in \Exbl(Z,d+1)} \left[d \dim\Spn(B) + 1\right]\\
      &= \sum_{B \in \Exbl(Z,d+1)} \dim \Spn(B).
    \end{aligned}\]
  Now noting that \(\delta_d = \sum_{A \in \Exbl(Z,d+2)} \left[ (d+1)\dim\Spn(A) + 1\right]
  - \sum_{B \in \Exbl(Z,d)} \left[d \dim\Spn(B) + 1\right]\)
  establishes the result.
\end{proof}

The following proposition shows that if the splitting type
\((a_1,..,a_n)\) is known it suffices to check if \(Z\) admits very unexpected hypersurfaces by
only looking around the degrees in the splitting type.

\begin{prop}Let \(Z \subseteq \bb{P}^n\) and let \((a_1,a_2,..,a_n)\)
  denote the splitting type of  \(D_0(\cl{A}_Z)\). If \(Z\) does not
  admit very unexpected hypersurfaces in degree \(d\), but does admit
  them in either degree \(d-1\) or degree \(d+1\), then \(d = a_i\)
  for some \(i\).
\end{prop}

\begin{proof}
  Let \(d\) be an index satisfying the hypothesis. Define indexes
  \(j\) and \(\ell\) so that \(a_k < d\) for all \(k
  \leq j\), and  \(a_k < d+1\) for all
  \(k \leq \ell\). The proposition is established if we show
  \(\ell > j\).

  Applying the inequality from \cref{cor:BestDefnUHyp} in degrees
  \(d-1,d\) and \(d+1\) we obtain the following three equations
\begin{enumerate}[\bf (Eq. 1)]
\item \(  \Exc(Z,d-1) + \sum_{k=1}^j (d- 1 - a_k) \geq n (d-1) + 1 \)
\item \(\Exc(Z,d) + \sum_{k=1}^\ell (d - a_k) = nd + 1\); and
\item \(\Exc(Z,d+1) + \sum_{k=1}^\ell (d+1 - a_k) \geq n(d+1) + 1\).
\end{enumerate}

Subtracting {\bf (Eq. 1)} from {\bf (Eq. 2)}  and {\bf (Eq. 2)} from {\bf (Eq. 3)}
gives   {\bf (Eq. 4)} and  {\bf (Eq. 5)} below.

\begin{enumerate}[\bf (Eq. 1)]
  \setcounter{enumi}{3}
\item \(\delta_{d-1} +j = \Exc(Z,d) - \Exc(Z,d-1) + j \leq n\)
\item \(\delta_d +\ell = \Exc(Z,d+1) - \Exc(Z,d) + \ell \geq  n \)
\end{enumerate}

By the preceding lemma \(\delta_d \leq
\delta_{d-1}\), with this and {\bf (Eqs. 4 \& 5)} we have
\[\delta_d + j \leq \delta_{d-1} + j \leq n \leq \delta_{d} + \ell.\]
We may conclude that \(j < \ell\), if either {\bf (Eq. 4)} or
{\bf (Eq. 5)} is strict. Yet this happens precisely when \(Z\) admits very
unexpected hypersurfaces in degree \(d-1\) or \(d+1\).
\end{proof}

\begin{example}\label{example:FiniteProjectiveSpaces} Let \(\bb{F}_q\)
  be the finite field with \(q = p^e\) elements, and \(\bb{K}\) an
  infinite field containing \(\bb{F}_q\). Let \(
  \bb{P}^n_{\bb{F}_q} \subseteq \bb{P}^n_{\bb{K}}\) consist of those
  points which in homogeneous coordinates can be written as
  \((\alpha_0:\alpha_1:..:\alpha_n)\) with \(\alpha_i \in
  \bb{F}_q\). It is well known that \(|\bb{P}^n_{\bb{F}_q}| = \frac{q^{n+1}-1}{q-1} =
  q^{n} + q^{n-1} + \ldots + q + 1\), and that \(\cl{A}_{\bb{P}^n_{\bb{F}_q}}\) is free
  with exponents \((1,q,q^2,...,q^n)\). The generator in degree
  \(q^i\) is of the form
  \[\sum_{j=0}^n Y_j^{q^i} \partialfrac{}{Y_i}\]
  and so the corresponding generator of \(I^{\gg}(\bb{P}^n_{\bb{F}_q})\) is
  \[\sum_{j=0}^n M_j^q X_i = \left|\begin{matrix}X_0 & X_1 & \ldots &
      X_n \\
      X_0^{q^i} & X_1^{q^i} & \ldots & X_n^{q^i}\\
      A_{1,0} ^{q^i} &A_{1,1}^{q^i} &\ldots &A_{1,n}^{q^i}\\
      \vdots & & \ddots & \vdots \\
    A_{n-1,0}^{q^i} &\ldots & \ldots &
    A_{n-1,n}^{q^i}\end{matrix}\right|\]

Furthermore, note that \(\bb{P}^n_{\bb{F}_q}\) is acted on by the
group \(\bb{GL}(n,\bb{F}_q)\). This in particular contains the
irreducible reflection group consisting of the permutation matrices, so by \cref{prop:ReflectIsGood} \(\bb{P}^n_{\bb{F}_q}\)
admits very unexpected hypersurfaces in all degrees \(d\) with \(q < d
< q^n\).
\end{example}

\begin{example}\label{example:FermatTypeArrangements} Fix some primitive \(m\)-th
  root of unity \(\zeta \in \bb{C}\), for \(m \geq 2\). Define a configuration of
  points \(F_m \subseteq \bb{P}^n_{\bb{C}}\) as consisting of the \(m \binom{n+1}{2}\)
  points whose \(i\)-th coordinate is \(-1\) and \(j\)-th coordinate
  is \(\zeta^k\) for all \(0 \leq k \leq d-1\) and all pairs \(0 \leq
  i < j \leq n\). Let \(C_m = F_m \cup \{E_0,E_1,..,E_m\}\) here \(E_i\) is the
  \(i\)-th coordinate point.

  Then \(\cl{A}_{C_m}\) is an Extended Ceva Arrangement, it is a
  reflection arrangement corresponding to the reflection group
  \(G(m,1,n+1) \subseteq \bb{PGL}(\bb{C},n)\). The splitting type of
  \(D_0(\cl{A}_{C_m})\) is \((m+1,2m+1,...,nm+1)\) (see \cite{ot1992}
  for details). As \(\cl{A}_{C_m}\) is a reflection arrangement, we again apply
  \cref{prop:ReflectIsGood} to conclude that \(\cl{A}_{C_m}\) admits
  very unexpected hypersurfaces in all degrees \(d\) with \(m+1 < d < nm+1\).
\end{example}

Both of our classes of examples come from reflection arrangements,
more generally \cref{prop:ReflectIsGood} gives a good criterion for
determining if the points dual to a given reflection arrangement admit
unexpected hypersurfaces. We note that reflection arrangements have
been classified and that their exponents and hence their splitting
type can be found in the appendix of \cite{ot1992}.

Our final example shows that the degrees in which a set of points
\(Z\) admits very unexpected hypersurfaces do not need to be
consecutive. This is in contrast with the situation in the plane as
shown in \cref{thm:UnexpectedCurveChar}. Before
outlining the example we state a useful proposition and definition.

\begin{defn}Let \(V_1\) and \(V_2\) be finite dimensional
  \(\bb{K}\)-vector spaces, and suppose we have finite sets of points \(Z_1 \subseteq \bb{P}(V_1)\), \(Z_2
  \subseteq \bb{P}(V_2)\). There are inclusion maps \(\iota_i:\bb{P}(V_i) \to  \bb{P}(V_1 \oplus V_2)\) for
  \(i=1,2\). We then define \(Z_1  \oplus Z_2 \subseteq \bb{P}(V_1 \oplus
  V_2)\) as the set of points
  \[Z_1 \oplus Z_2 := \iota_1(Z_1) \cup \iota_2(Z_2).\]

\end{defn}

\begin{prop}\label{prop:UnexForSumOfPoints} \(Z_1 \oplus Z_2\) admits very unexpected curves in
  degree \(d \geq 1 \) if and
  only if \(Z_1\) or \(Z_2\) admits unexpected curves in degree \(d\).
\end{prop}

\begin{proof} First, note that for hyperplane arrangements
  \(\cl{A}_1 \subseteq \bb{P}(W_1)\) and \(\cl{A}_2\ \subseteq
  \bb{P}(W_2)\) there is an arrangement \(\cl{A}_1 \cross
  \cl{A}_2 \subseteq \bb{P}(W_1 \oplus W_2)\) induced by the
  projections \(p_i:\bb{P}(W_1 \oplus
  W_2) \to \bb{P}(W_i)\). Namely, \(\cl{A}_1 \cross \cl{A}_2\) is
  formed by taking all hyperplanes of the form \(\pi_i^{-1}(H)\) for
  \(H \in \cl{A}_i\).

  We now note two facts:
  \begin{enumerate}[\bf ({Fact} 1)]
  \item \(\cl{A}_{Z_1 \oplus Z_2} = \cl{A}_{Z_1} \cross
    \cl{A}_{Z_2}\),
  \item If \(S\) is the projective
    coordinate ring of \(\bb{P}(W_1 \oplus W_2)\) there is an isomorphism
    of \(S\)-modules,
    \(D(\cl{A}_{Z_1} \cross \cl{A}_{Z_2}) \iso  (S \otimes
    D(\cl{A}_{Z_2})) \oplus (S \otimes D(\cl{A}_{Z_1})).\)
  \end{enumerate}
  
  The first can be seen by following each the constructions through the
  duality. We omit a proof of the second referring to \cite{ot1992} for details.

  One consequence of fact 2 is that if \(D_0(\cl{A}_1)\) has splitting
  type \((a_1,..,a_n)\) and \(D_0(\cl{A}_2)\) has splitting type
  \((b_1,..,b_m)\), then \(D_0(\cl{A}_1 \cross \cl{A}_2)\) has a
  splitting type (up to reordering) of \((1,a_1,...,a_n,b_1,..,b_m)\).
  Applying \cref{cor:BestDefnUHyp}, now yields the inequalities valid
  for any \(d \geq 1\). Each inequality strict if and only if the
  corresponding set of points admits very unexpected hypersurfaces
\begin{enumerate}[\bf ({Ineq.} 1)]
\item \( \sum_{i=1}^n \max\{0,d-a_i\} \geq nd+1 - \Exc(Z_1,d) \)
\item \( \sum_{j=1}^m \max\{0,d-b_j\} \geq md+1 - \Exc(Z_2,d) \)
\item \( d-1 + \left(\sum_{i=1}^n \max\{0,d-a_i\}\middle)+\middle(\sum_{j=1}^m \max\{0,d-b_j\}\right) \geq (n+m+1)d+1 -
  \Exc(Z_1\oplus Z_2,d) \)
\end{enumerate}
We now claim that \(\Exc(Z_1 \oplus Z_2,d) = \Exc(Z_1,d) +
\Exc(Z_2,d)\). First note that if we assume this claim and subtract \(d-1\)
from both sides of {\bf (Ineq. 3)}, then the resulting inequality may be
written as the sum of {\bf (Ineq. 1)} and {\bf (Ineq. 1)}. From this it follows that {\bf (Ineq. 3)} is
strict if and only if either {\bf (Ineq. 3)} or {\bf (Ineq. 3)} is
strict and the proposition follows.

Continuing to the proof of our claim, we first note that if \(d =1\)
then for any set of points \(\Exbl(Z,1) = \{Z\}\).
A direct computation establishes the claim in this case.

Now we may assume \(d \geq 2\). Take a block \(B \in \Exbl(Z_1 \oplus
Z_2,d)\) and define \(B_1 = B \cap Z_1\) and \(B_2 = B \cap Z_2\). We
note that if \(B_1\) and \(B_2\) are
nonempty, then \cref{lem:BlockLemma} states
\[\displaystyle d (\dim \Spn B - \dim \Spn B_1 - \dim\Spn
  B_2) \leq 1.\]
Yet as \(B_1\) and \(B_2\) are contained in disjoint subspaces, \(\dim
\Spn(B) = \dim \Spn(B_1) + \dim \Spn(B_2) + 1\) and the inequality
becomes \(d \leq 1\) giving a contradiction. Therefore, for each block \(B\) we have \(B \subseteq Z_1\) or
\(B \subseteq Z_2\), and consequently \(\Exbl(Z_1 \oplus Z_2,d) = \Pi_1 \cup
\Pi_2\) for some partitions \(\Pi_1\) and \(\Pi_2\) of \(Z_1\) and \(Z_2\)
respectively. From this it readily follows from the definition that
\(\Exbl(Z_1 \oplus Z_2,d) = \Exbl(Z_1,d) \cup \Exbl(Z_2,d)\). This
establishes the claim that \(\Exc(Z_1 \oplus Z_2,d) = \Exc(Z_1,d) +
\Exc(Z_2,d)\) and completes our proof. 
\end{proof}

\begin{example}\label{example:OneMoreThing} If \(C_2,C_7 \subseteq \bb{P}^2_{\bb{C}}\) are the
  configurations of points described in
  \cref{example:FermatTypeArrangements}, then \(C_2 \oplus C_7\) is a
  configuration of \(33\) points in \(\bb{P}^5_{\bb{C}}\). The
  module of derivations, \(D_0(\cl{A}_{C_2} \times
  \cl{A}_{C_7})\), has splitting type \((1,3,5,8,15)\). Using the
  computation from
  \cref{example:FermatTypeArrangements} along with
  \cref{prop:UnexForSumOfPoints}, it follows that \(C_2 \oplus
  C_7\) admits very unexpected hypersurfaces in degree \(d\) if and only if
  \(d =4\) or \(8 < d < 15\).

  \end{example}

\PFEIndent

\section{A Lifting Criterion and the Structure of Unexpected Curves
  in \texorpdfstring{\(\bb{P}^2_{\mathbb{C}} \)}{the Complex Plane}}\label{sec:LiftingAndCurves} 

One feature of the \cref{thm:GlobalIso,thm:restrictedModIso}, is they
allow us to view elements of reduced Module of Derivations as explicit
polynomials. This permits us to use techniques such as unique factorization and
polynomial division that are not as well developed for
general modules. In this section we give a few applications of this
view point. First, we state a lifting criterion in
\cref{prop:DivisibilityProp}, this allows us under certain conditions to lift an element of the
restricted module \(D_0(\cl{A}_Z)\mid_L \) to the module
\(D_0(\cl{A}_Z)\). This criterion has especially strong implications
in \(\bb{P}^2_{\bb{C}}\), such as in \cref{thm:UnexpectednessIsGlobal}, where we show that for \(Z\subseteq \bb{P}^2\) every
polynomial defining an unexpected curve in \(I^{\gg}_Q(Z)\) can be
lifted to an element of \(I^{\gg}(Z)\).

This result ends up putting very strong conditions on the
combinatorics of sets of points \(Z\) which admit unexpected
curves, which we explore in the next section.

\begin{prop}\label{lemma:DivisibilityLemma}Let \(Z \subseteq
  \bb{P}(V) \iso \bb{P}^n\). Consider \(G \in I^{\gg}(Z) \subseteq
  \bb{K}[\mathbf{A}]\). If there is some \(F \in \bb{K}[\mathbf{A}]\) so
  that for general \(\bm{\alpha} \in \textrm{Gr}(n-1,V)\) we have that
  \(\varepsilon_{\bm{\alpha}}(F) \in I^{\gg}_{\bm{\alpha}}(Z)\) and 
  \(\varepsilon_{\bm{\alpha}}(F) \mid \varepsilon_{\bm{\alpha}}(G)\). Then \(F\) and
  \(G\) have a common divisor \(H \in I^{\gg}(Z)\).
\end{prop}

\begin{proof} For any prime ideal, \(I\), we set \(\nu_I(F)\) as the
  valuation \(\nu_I(F) := \sup\{m \geq  0 \mid F \in I^{(m)}\}\). Now we
  define two ideals of \(\bb{K}[\bold A]\),
  \(X\) is the ideal \((X_0,..,X_n)\) and we let \(M\) denote the
  ideal generated by the maximal minors of the matrix \(\mathbf{A}\).
  Lastly, for \(\bm{\alpha} \in \Gr(n-1,V)\), \(I(\bm{\alpha})\) is the ideal of
  \(\bb{K}[X_0,..,X_n] \subseteq \bb{K}[\bold A]\) defined by the
  subspace \(\bm{\alpha}\).

  Before continuing we note a few facts:
  \begin{enumerate}[\bf {Fact} 1]
   \item For each of the \(3\) ideals, \(X\), \(I(\alpha)\) and \(M\), that we have defined
     we have \(I^k = I^{(k)}\).
   \item For any \(f \in \bb{K}[\mathbf{A}]\) we have  \(\nu_X(f)
     \geq \nu_M(f)\) and
     \(\nu_{I(\bm{\alpha})}(\varepsilon_{\bm{\alpha}}(f)) = \nu_M(f)\)
     for general \(\bm{\alpha}\).
   \item For any \(f \in \bb{K}[\mathbf{A}]\), we have the inequality
     \[(\nu_X(f) - \nu_M(f)) + n^2 \nu_M(f) \leq \deg(f).\]
   \item If \(\nu_X(f) = \nu_M(f)+1\), then equality occurs in {\bf Fact 3} if and only if \(f
     \in \mathfrak{m}^{\gg}\).
   \end{enumerate}
    The first fact  follows for \(X\) and \(I(\bm{\alpha})\) since
    both are complete intersections, for \(M\) we refer to section 2.2 of \cite{Hoch}.
    The second fact follows since \(M\) is essentially \(I(\alpha)\) for
    \(\alpha\) the generic point. The third is a consequence of the first and
    that \(\deg(M_i) = n^2\). Lastly, the fourth fact follows since if
    \(\nu_X(f) = \nu_M(f)\) then setting \(d = \nu_X(f)\) we have  \(f
    \in [(X_0,..,X_n) M^{d-1}]_{n^2(d-1) + 1}\), but this is precisely
    \([\mathfrak{m}^{\gg}]_d\).

    First we claim for general \(\bm{\alpha}\), that
    any \(f  \in I^{\gg}_{\bm{\alpha}}(Z)\) factors into irreducible components as \(f = f_0
    \prod_{i=0}^k \ell_i\) where each \(\ell_i\) is an element of the
    special fibre ring \(\cl{F}_{\bm{\alpha}} =
    \Sym([I({\bm{\alpha}})]_1)\). It suffices to show that if \(f = 
    pq\), then either \(p\) or \(q\) is in
    \(\cl{F}_{\bm{\alpha}}\). Noting that \(\nu_X(f) = 1 +
    \nu_{I({\bm{\alpha}})}(f)\), and using the additive property of
    valuations, we obtain
    \[\nu_X(p) + \nu_X(q) = 1 + \nu_{I({\bm{\alpha}})}(p) +
      \nu_{I({\bm{\alpha}})}(q) \leq 1+\nu_X(p) + \nu_X(q).\]
    
    Since all numbers above are integers, and
  \(\nu_{I({\bm{\alpha}})}(h) \leq \nu_X(h)\) for every \(h \in \bb{K}[X_0,..,X_n]\), we may
  assume without loss of generality that \(\nu_{I({\bm{\alpha}})}(p) = \nu_X(p)\) and
  \(\nu_{I({\bm{\alpha}})}(q) = \nu_X(q)+1\). It now follows that \(p
  \in [I(\bm{\alpha})^{\nu_X(p)}]_{\nu_X(p)} \subseteq\cal{F}_{\bm{\alpha}}\), which establishes our claim.

   Continuing with the proof of the proposition,
   we let \(\varepsilon_{\mathfrak{g}}\) denote the generic
   evaluation. In other words \(\varepsilon_{\mathfrak{g}}\) is the
   inclusion \(\varepsilon_{\mathfrak{g}}:\bb{K}[\mathbf{A}] \to \bb{F}[X_0,..,X_n]\),
   for \(\bb{F}\) the function field \(\bb{F} := \bb{K}(A_{i,j} \mid (1,0) \leq
   (i,j) \leq (n,n+1))\). By assumption \(\varepsilon_{\mathfrak{g}}(F)\) divides
   \(\varepsilon_{\mathfrak{g}}(G)\), so there exists \(h
   \in \bb{K}[\mathbf{A}]\) and \(k \in \bb{K}[A_{i,j} \mid (1,0) \leq
   (i,j) \leq (n,n+1)]\) with \(h\) and \(k\) coprime so that
   \(\frac{h}{k} \varepsilon_{\mathfrak{g}}(F) = \varepsilon_{\mathfrak{g}}(G)
   \iff hF = kG\) where this last equality is in \(\bb{K}[\mathbf{A}]\). Now by unique factorization in the
   polynomial ring \(\bb{K}[\mathbf{A}]\), we get \(k \mid F\).
   Setting \(\tilde{F}=\frac{F}{k} \in \bb{K}[\mathbf{A}]\), we have \(h \tilde{F} =
   G\) and \(k\tilde{F} = F\).  Moreover, since \(\varepsilon_{\mathfrak{q}}(F) \not\in
   \cal{F}_{\mathfrak{q}}\) it follows that \(\tilde{F} \not\in
   \cal{F}_{\mathfrak{q}}\) and so \(h \in
   \cal{F}_{\mathfrak{q}}\). We finish the proof by establishing that
   \(\tilde{F} \in I^{\gg}(Z)\).

   Since \(F\) differs from \(\tilde{F}\) only by \(\bb{F}\) scalar,
   and \(F \in I(Z)\) we have \(\tilde{F} \in I(Z)\) and so it
   suffices to show that \(\tilde{F} \in \mathfrak{m}^{\gg}\). Since
   \(\varepsilon_{\mathfrak{q}}(\tilde{F}) \in I^{\gg}_{\mathfrak{g}}(Z)\) and
   \(\varepsilon_{\mathfrak{q}}(h) \in \cal{F}_{\mathfrak{q}}\) we have the
   inequalities
   \[\begin{aligned}
       \nu_M(h) &\leq \nu_{I(\mathfrak{q})}(h) =
     \nu_X(\varepsilon_{\mathfrak{q}}(h)) = \nu_X(h); \text{ and } \\
       \nu_M(\tilde{F}) &\leq \nu_{I(\mathfrak{q})}(\tilde{F}) =
     \nu_X(\varepsilon_{\mathfrak{q}}(\tilde{F})) - 1 =
     \nu_X(\tilde{F})-1.
   \end{aligned}\]
 As \(h \tilde{F} \in I^{\gg}(Z)\), we have that \(\nu_M(h) + \nu_M(F)
 = \nu_X(h) + \nu_X(F) -1\) and so the above inequalities must be equality.
   Similarly, using the inequalities \(1 + n^2 \nu_M(\tilde{F}) \leq \deg(\tilde{F})\),
   \(n^2\nu_M(h) \leq \deg(h)\) and  \(\tilde{F} h = G \in I^{\gg}(Z)\) we have that
   \[1 + n^2 \nu_m(\tilde{F}) + n^2 \nu_M(h) \leq \deg\tilde{F} +
     \deg h = \deg G = 1 + n^2 \nu_M(\tilde{F} h).\]

   Allowing us to conclude that \(1 + n^2 \nu_M(\tilde{F}) =
   \deg(\tilde{F})\) and \(n^2\nu_M(h) = \deg(h)\) which completes the proof.
\end{proof}

The preceding lemma when combined with the results of \cref{sec:DerivationIdealCorrespondence} allows us under certain
circumstances to lift elements
of \(D_0(\cl{A})\mid_L\) to elements of \(D_0(\cl{A})\). One example of this is illustrated in the following
proposition.

\begin{prop}\label{prop:DivisibilityProp}  Let \(\cl{A} \subseteq
  \bb{P}^n_{\bb{K}}\) and suppose that \(D_0(\cl{A})\) has splitting type
  \((a_1,a_2,..,a_n)\) with \(a_1 < a_2 \leq a_3 \leq ... \leq
  a_n\). If \(\theta_\lambda \in D_0(\cl{A})\) is a nonzero element of degree
  \(< a_2\), then \(D_0(\cl{A})\) has a minimal generator in degree \(a_1\). 
\end{prop}

\begin{proof} Using the translation given by \cref{thm:GlobalIso},
  there's a nonzero \(F_\lambda \in [I^{\gg}(Z)]_d\) where \(d <
  a_2+1\). If \(Q\) is the generic codimension \(2\) linear subspace,
  then by \cref{thm:restrictedModIso} \(I^{\gg}_Q(Z)\) is free on generators \(f_1,..,f_n\) with
  \(\deg f_i = a_i + 1\). Hence, \(\varepsilon_Q(F_\lambda) =
  \sum_{i=1}^n g_i f_i\). Yet as \(\deg f_j > \deg
  \varepsilon_Q(F_\lambda)\) for all \(j \geq 2\), we must have
  \(\varepsilon_Q(F_\lambda) = g_1 f_1\). After clearing denominators we may 
  lift \(f_1\) to an element \(\tilde{f}_1\) of
  \(\bb{K}[\mathbf{A}]\).

  Now as \(\varepsilon_Q(\tilde{f}_1)\) divides
  \(\varepsilon_Q(F_\lambda)\) we see by the previous lemma that there exists \(F_1
  \in I^{\gg}_Q(Z)\) which divides both \(\tilde{f}_1\) and
  \(F_\lambda\). As \(F_1\) divides \(\tilde{f}_1\) we must have
  \(F_1 \in [I^{\gg}(Z)]_{a_1+1}\) and so by
  \cref{thm:GlobalIso} there's a nonzero \(\theta_1 \in [D_0(\cl{A}_Z)]_{a_1}\).
\end{proof}

The previous two propositions will prove to be especially useful when our points
(or line arrangements) are in the plane \(\bb{P}^2\). We will
establish this using some results on vector bundles on
\(\bb{P}^2_{\bb{C}}\) which we recall now.

\begin{defn} We say a vector bundle \(M\) on \(\bb{P}^n_{\bb{C}}\) is semistable, if for all
  proper subbundles \(N \subsetneq M\), we have
  \[\frac{c_1(N)}{\rank N} \leq \frac{c_1(M)}{\rank M}\]
  where here \(c_1\) is the first Chern class, and \(\rank M\) is the
  dimension of a fibre.
\end{defn}

If \(\rk(\cal{M}) = 2\), semistability has a simpler characterization
originally due to Hartshorne (see lemma \(3.1\) of \cite{Hart}).

\begin{lemma}\label{lemma:HartshorneLemma} Let \(\cal{M}\) be a rank \(2\) bundle on \(\bb{P}^n_{\bb{C}}\), then
  if \(\cal{M}\) is semistable if and only if letting \(c_1 = c_1(\cal{M})\) 
  \[H^0\left(\cal{M}\middle(\left\lfloor\frac{-c_1-1}{2}\right\rfloor\middle)\right)
    = 0.\]
\end{lemma}

In the case that our bundle \(\cal{M}\) is the
Derivation Bundle, \(\widetilde{D_0(\cl{A})}\) of a hyperplane arrangement, it was shown by Terao that \(c_1(M) = 1-|\cl{A}| \),
where here \(|\cl{A}|\) is the number of hyperplanes in
\(\cl{A}\). Using this together with the previous lemma now allows us to characterize semistability of
\(\widetilde{D_0(\cl{A})}\) for \(\cal{A}\) a line
arrangement in \(\bb{P}_\bb{C}^2\). 

\begin{prop}\label{prop:SemistableCharac} For \(\cl{A} \subseteq \bb{P}^2_\bb{C}\) define \(d =
  \left\lfloor\frac{|A|-2}{2}\right\rfloor\). Then the derivation
  bundle \(\cl{D}_0(\cl{A})\) is semistable if and only if
  \([D_0(\cl{A})]_{d} = 0\). In particular, if \(\cl{D}_0(\cl{A})\) is not semistable, then
  \(D_0(\cl{A})\) contains a nonzero derivation in degree  \(\left\lfloor\frac{|A|}{2}\right\rfloor-1\).
\end{prop}

One property of semistable bundles is the celebrated theorem of
Grauert and M\"{u}lich, which characterizes the splitting type of
semistable bundles. In light of \cref{thm:UnexpectedCurveChar}, we see
that if \(D_0(\cl{A}_Z)\) is semistable then \(Z
\subseteq \bb{P}^2_{\bb{C}}\) admits no unexpected curves.

\begin{theorem}[Grauert-M\"{u}lich] If \(\cl{B}\) is a semistable bundle on
  \(\bb{P}^n_{\bb{C}}\), with splitting type \(a_1 \leq a_2 \leq
  ... \leq a_k\), then for all \(1 \leq i < k\), we have \(0 \leq
  a_{i+1} - a_i \leq 1\).
\end{theorem}

\begin{theorem}\cite{CHMN}\label{thm:SSImpliesNoUnex} For \(\bb{Z} \subseteq \bb{P}(V) \iso
  \bb{P}^2_\bb{C}\), if \(\cl{D}_0(\cl{A}_Z)\) is semistable then
  \(Z\) admits no unexpected curves.
\end{theorem}

This theorem in conjunction with \ref{prop:DivisibilityProp} allows us to say
that every unexpected curve in \(\bb{P}^2_{\bb{C}}\) comes from a global
section of the derivation bundle. More precisely a polynomial
defining a degree \(d\) unexpected curve corresponds via the duality
of \cref{thm:GlobalIso} to an element of \([D_0(\cl{A})]_{d-1}\).

\begin{theorem}\label{thm:UnexpectednessIsGlobal} Let \(Z \subseteq
  \bb{P}^2_{\bb{C}}\) be a finite set of points. If \(D_0(\cl{A}_Z)\)
  has splitting type \((a_1,a_2)\) with \(a_2 - a_1 \geq 2\) (in
  particular if \(Z\) admits unexpected curves), then \(D_0(\cl{A}_Z)\) has a generator in degree \(a_1\); Equivalently, \(I^{\gg}(Z)\) has a
  generator in degree \(a_1+1\).
\end{theorem}

\begin{proof} As \(a_2 - a_1 > 1\), \(\widetilde{D_0(\cl{A}_Z)}\) is
  not semistable by the Grauert-M\"{u}lich
  theorem. Hence by \cref{prop:SemistableCharac}, there's a nonzero \(\theta
  \in D_0(\cl{A})\) with \(\deg \theta \leq \left\lfloor \frac{|Z|-2}{2}
  \right\rfloor < a_2\). Applying \cref{prop:DivisibilityProp} now yields
  the required generator of degree \(a_1\).
\end{proof}

Combining this with the Grauert-M\"{u}lich Theorem, we obtain the
following result.

\SplittingTypeChar

Translating the above statement via the duality of
theorem \ref{thm:GlobalIso}, we obtain the
corollary below. 

\begin{cor}\label{cor:PolynomialStructure} For a finite set of points \(Z \subseteq
  \Proj(\bb{C}[X_0,X_1,X_2])\), suppose \(Z\) admits unexpected
  curves. If \(Q = (A_0:A_1:A_2)\) is the generic point, then
  \(I^{\gg}_Q(Z)\) is a free \(\cl{F}_Q\)-module on generators \(f\) and \(g\) with \(\deg
  f < \deg g-1\) and \(f\)
  can be lifted to an element \(F\)  of \(I^\gg(Z)\) with
  \(\epsilon_Q(F) = f\).

  In particular, \(f\) can be written as 
  \[f = X_1f_1 + X_2 f_2 + X_3f_3,\]
  where \(f_i\) is a polynomial of degree \((\deg f)-1\) in the maximal
  minors of \(\begin{bmatrix} A_0 & A_1 & A_2 \\ X_0 & X_1 & X_2 \end{bmatrix}\)
\end{cor}

The above corollary states that the polynomials defining unexpected
curves are ``as simple as possible'', in the sense that they have
the minimal possible degree as a polynomial in the coordinates of our
general point \(Q\). This stands in stark contrast to most other sets of points
where this is not the case. As an illustration, taking \(8\) randomly chosen points \(Z \subseteq
\bb{P}^2_{\bb{C}}\), a computation with Macaulay2 showed that 
\(I^{\gg}_Q(Z)\) has generators of \(X\)-degree \(4\) and \(5\). The first generator
had an \(A\)-degree of \(12\) giving a total degree of \(16\), showing that
the above result is far from expected. Similar computations with
\(6\) points and \(10\) points gave minimal polynomials with \(A\)-degrees
of \(6\) and \(20\), respectively.

Below we present a simpler example illustrating a similar point.

\begin{example}\label{example:PolynomialStructure}If \(Z \subseteq
  \bb{P}^2_{\bb{C}}\) consists of the \(3\) coordinate points and
  \((1:1:1).\)  Then for generic \(Q\), \(I^{\gg}_Q(Z)\) has
  generators \(f_1\) and \(f_2\) of  degrees \(2\) and \(3\)
  as polynomials in \(X\). Many different \(f_2\) are
  possible. On the other hand, if we
  require \(f_1\) to be a polynomial of minimal degree in \(\bb{K}[\bold
  A]\), it is unique up to \(\bb{C}\) scalar.  The corresponding polynomial formula is
  \[f_1 = (A_0-A_1)M_1X_1 + (A_0-A_2)M_2X_2=(A_0-A_1)A_2 X_0X_1 -
    (A_0-A_2)A_1X_0X_2 + (A_1-A_2)A_0X_1X_2 \]
  where here \(M_i\) is the minor of the \(2 \by 3\) matrix from the
  matrix above and \cref{prop:TauIso}. \(f_1\) is irreducible, and
  defines the unique smooth conic through \(Z\) and \(Q\). As the
  \(A\) degree and \(X\) degree of \(f_1\) are the same, we can see \(f_1\)
  cannot be written in the form from \cref{cor:PolynomialStructure}.
\end{example}

\PFEIndent

\section{Combinatorial Constraints on Points Admitting Unexpected
  Curves}\label{sec:CombinatoricsInP2}

In this section we explore combinatorial constraints necessarily
satisfied by sets of points admitting unexpected curves. Most of these constraints apply only
when this unexpected curve is irreducible. Yet this turns out to be a
fairly weak assumption, since if \(Z\) admits a unique unexpected curve in degree \(d\)
there is always a subset \(W \subseteq Z\) so \(|Z \setminus W| = k\)
and \(W\) admits a unique irreducible unexpected curve in degree
\(d-k\).

We start by
exploring the consequences of  \cref{cor:PolynomialStructure}. In the
case the curve of degree \(d\) is irreducible we show in \cref{prop:LpBound} that
\cref{cor:PolynomialStructure} gives a bound on the number of distinct lines through a point \(P \in Z\) and
the remaining points of \(Z\), showing that there are at most \(d\) lines. As \(|Z| \geq 2d+1\) this is a very strong
combinatorial condition which states that on average each line
through a fixed point \(P\) contains \(3\) or more points of \(Z\). We are
able to use this in \cref{thm:PointCardBound} to give a sharp bound on the number of points
in \(Z\), this bound is achieved by the Ceva type point configurations \(C_d\)
from \cref{example:FermatTypeArrangements}.

Furthermore, in \cref{prop:LineCardBound} we also give and upper bound on the
number of lines spanned by points of \(Z\). We then close the section
by applying a theorem of Terao to state a combinatorial condition that guarantees that \(Z\) will
admit an unexpected curve.

We note that throughout this section, we often state theorems with the
assumption that \emph{``there's some nonzero (possibly irreducible) \(f \in
[I^{\gg}(Z)]_d\)''}. By \cref{thm:UnexpectednessIsGlobal} perhaps the
prototypical example for us are points 
\(Z \subseteq \bb{P}^2_{\bb{C}}\) admitting unexpected curves in
degree \(d\). However, this also holds in other contexts for instance if
\(\cl{A}_Z \subseteq \bb{P}^2_{\bb{K}}\) is free.

\begin{lemma}\label{lemma:VanishingEquivalence} For \(Z \subseteq \bb{P}(V) \iso \bb{P}^2_{\bb{K}}\), consider
  \(F(X_0,X_1,X_2;A_0,A_1,A_2) \in [I^{\gg}(Z)]_{d}\) then for every \(P = (P_0:P_1:P_2) \in
  Z\),
  \[ \varepsilon_P(F) = F(X_0,X_1,X_2;P_0,P_1,P_2) \in I(P)^{d}.\]
  
  Moreover, \(\varepsilon_P(F) = 0\) if and only if the
  linear form \(\ell_P = P_0M_0 + P_1M_1+P_2M_2\) divides \(F\).
\end{lemma}

\begin{proof}
  We may choose coordinates so \(P = (1:0:0)\) and
  \(\ell_P = M_0\). Writing \(F\) as \(F
  = F_0 X_0 + F_1 X_1 + F_2 X_2\), then as \(F \in I(P)=(X_1,X_2)\)  we have
  \[F(P_0,P_1,P_2,A_0,A_1,A_2) =
    F_0(1,0,0,A_0,A_1,A_2) =  0.\]
  As each \(M_i\) is antisymmetric in \(A\) and \(X\), it follows that
  \(\varepsilon_P(F_0)  = 0\) and so \(F_0 \in
  I(P^{\perp}) = M_0\) by \cref{prop:epsilonIso}. Then applying the identity \(M_0X_0+M_1X_1+M_2X_2 = 0\),  we
  may write \(F =  f_1X_1+f_2X_2\) where \(f_i =
  \left(F_i-\frac{M_i}{M_0}F_0\right)\). Noting for arbitrary \(Q \in \bb{P}^2\), that
  \(\varepsilon_Q(f_i) \in I(Q)^{d-1}\). It follows that
  \[\varepsilon_P(F)= \varepsilon_P(f_1)X_1+\varepsilon_P(f_2)X_2 \in (X_1,X_2)I(P)^{d-1} =
  I(P)^d,\] which  establishes the first statement.

  Continuing with the proof of the second statement,  assume that \(\varepsilon_P(F) =
  \varepsilon_P(f_1)X_1+\varepsilon_P(f_2)X_2 = 0\). Noting \(\varepsilon_P(M_1) = X_2\)
  and \(\varepsilon_P(M_2) = -X_1\), it follows that
  \(\varepsilon_P(f_1M_2 - f_2M_1) = - \varepsilon_P(f_1)X_1-\varepsilon_P(f_2)X_2 =
  0\), so \(f_1M_2 - f_2M_1 \in (M_0) = \ker \varepsilon_P\). Let
  \(\tilde{f_1},\tilde{f_2} \in \bb{K}[M_1,M_2]\) so that \(f_i =
  \tilde{f_i} \mod (M_0)\) for each \(i\in\{1,2\}\). Then \(\tilde{f_1}M_2 -
  \tilde{f_2}M_1 \in (M_0) \cap \bb{K}[M_1,M_2] = 0\), so
  \(\tilde{f_1}M_2 = \tilde{f_2}M_1\) and we get by unique factorization that there exists some \(g \in
  \bb{K}[M_1,M_2]\) with \(g = \dfrac{\tilde{f_2}}{M_2} =
  \dfrac{\tilde{f_1}}{M_1}.\)

  Finally, applying the identity \(X_0M_0+X_1M_1+X_2M_2 = 0\) again we
  can write 
  \[F = f_1 X_1 + f_2 X_2 - g(M_0X_0+M_1X_1+M_2X_2) = (-gM_0)X_0 +
    (f_1-gM_1)X_1+(f_2-gM_2)X_2.\]
  Noting that \(f_i - gM_i = f_i - \tilde{f_i}
  \equiv 0 \mod (M_0)\), we conclude that \(M_0\) divides \(F\).
  This establishes the forward direction, the reverse direction follows as \(\varepsilon_P(\ell_P) = 0\).
\end{proof}

As we will see, the preceeding lemmas imposes a very strong combinatorial condition on
the configurations of points which can admit unexpected curves. Before
we state the first of these conditions we introduce a new piece of notation.

\begin{defn}\label{defn:LpDefn} Let \(Z \subseteq \bb{P}^2_{\bb{K}}\) be a finite
  configuration of points, with \(|Z| \geq 2\). For each \(P \in \bb{P}^2_{\bb{K}}\),  define a
  set of lines,  \(L_P(Z)\), as follows
  \[L_p(Z) := \{ \Spn(Q_i,P) \mid Q_i \in Z \setminus \{P\}\}.\]
\end{defn}

\begin{remark}
  Note \(|L_P(Z)| \leq |Z \setminus \{P\}|\) with equality if and
  only if for distinct \(Q,Q' \in |Z \setminus \{P\}|\) we have \(\Spn(Q,P)
  \neq \Spn(Q',P)\).
\end{remark}

The number \(|L_p(Z)|\) defined above has an equivalent purely algebraic definition.

\begin{lemma}\label{lemma:LpEquivalence} Let \(Z \subseteq \bb{P}^2_{\bb{K}}\) be a finite set of at
  least \(2\) points, then for any \(P \in \bb{P}^2_{\bb{K}}\) we have
  \[|L_P(Z)| = \min\{d \mid [I(Z) \cap I(P)^{d}]_d \neq 0\}\]
\end{lemma}

\begin{proof} Let \(m = \min\{d \mid [I(Z) \cap I(P)^{d}]_d \neq
  0\}\). For any \(Q \in Z \setminus \{P\}\), we get by Bezout's
  Theorem that the line \(\Spn(P,Q)\) must be
  a component of the base locus of \([I(Z) \cap I(P)^d]_d\). Hence,
  letting \(G_p\) denote the product of the linear forms defining
  the elements of \(L_p(Z)\), we have \(\deg G_p = |L_p(Z)|\) and \([I(Z) \cap I(P)^d]_d =
  [I(L_p(Z))]_d = [(G_p)]_d\) which completes the proof.
\end{proof}

We introduce the first combinatorial constraint below, it occurs
whenever \(I^{\gg}(Z)\) contains an irreducible element. As we will
see this simple constraint ends up having a number of strong consequences.

\begin{lemma}\label{prop:LpBound}Let \(Z \subseteq \bb{P}^2_{\bb{K}}\), with \(|Z| \geq 2\) and suppose \(F \in
  [I^{\gg}(Z)]_d\) is an irreducible polynomial. Then for all \(P \in
  Z\) we have \(|L_P(Z)| \leq d\).

  Consequently, if \(Z \subseteq \bb{P}^2_{\bb{C}}\) admits an
  irreducible unexpected curve in degree \(d\), then \(|L_p(Z)| \leq
  d\) for all \(P \in Z\).
\end{lemma}

\begin{proof}As \(F\) is irreducible, we know by
  \cref{lemma:VanishingEquivalence}  that
  \(\varepsilon_P(F) \neq 0\) and \(\varepsilon_P(F) \in [I(P)^{d}
  \cap I(Z)]_d\) for all \(P \in Z\). Hence, by \cref{lemma:LpEquivalence} we get \(|L_P(Z)| \leq
  d\).

  The second statement follows from the first in light of \cref{cor:PolynomialStructure}.
\end{proof}

G. Dirac conjectured (see \cite{Dirac}) that for any set \(Z\) of noncollinear points in
\(\bb{R}^2\), there always exists some \(P \in Z\) with \(|L_P(Z)|
\geq \lfloor \frac{|Z|}{2}\rfloor\). This turned out to be
false. However, since then alternative conjectures have been proposed,
one version of the conjecture was established in \cite{Han} for points
in \(\bb{C}^2\). This result allows us \cref{thm:PointCardBound}
below. We explore possible further consequences of the conjectures in \cref{sec:ApplicationsToTeraos}.

\PointCardBound

\begin{proof} This follows from \ref{prop:LpBound} and
  Han's improvement of the Dirac Conjecture \cite{Han}, which states for a finite
  set of points \(Z \subseteq \bb{P}^2_{\bb{C}}\) which span
  \(\bb{P}^2\) there always exists some \(P \in Z\) so \(|L_p(Z)| \geq
  \frac{|Z|}{3} + 1\).

  Namely, suppose \(Z\) admits an irreducible curve in degree \(d\),
  then for all \(P \in Z\), \(|L_P(Z)| \leq d\). Now applying Han's result, there exists \(P \in
  Z\) so
  \[ d \geq |L_p(Z)| \geq \dfrac{|Z|}{3} + 1.\]
  Solving for \(|Z|\) now yields \(3 (d-1) \geq |Z|\), the desired inequality.

\end{proof}

\begin{remark}
  It should be noted that the paper \cite{Han}, is rather vague and
  states the result only for points in ``the plane''. However, the
  proof works for complex line arrangements, as the main nonelementary
  tool is a Hirzeburch type inequality for complex line arrangements first proved in \cite{Boj.}

  Equivalently, Langer's Inequality \cite{Lang.}, could replace and or
  rederive \cite{Han}'s result. Langer's Inequality states that
  letting \(\displaystyle \ell_r = \left|\{ L \subseteq \bb{P}^2_{\bb{C}} \mid |L \cap
  Z|=r\} \right|\) we have that if \(\ell_r = 0\) for \(r > \frac{2}{3} |Z|\) that 
  \[ \sum_{P \in Z} |L_P(Z)| = \sum_{r \geq 2} r \ell_r \geq
    \left\lceil\frac{|Z|^2+3|Z|}{3}\right\rceil.\]
  For further discussion see the survey article \cite{Pok}, where the
  author first learned of these results.
\end{remark}

In \cite{CHMN} it was shown that any set of points \(Z \subseteq \bb{P}^2_{\bb{C}}\) in linearly general position
can never admit unexpected curves. This proposition provides a
strengthening of that result, and extends it to an arbitrary field.

\begin{prop}\label{prop:IrrLinearlySpecial} Let \(Z \subseteq
  \bb{P}^2_{\bb{K}}\) suppose there's a nonzero \(F \in
  [I^{\gg}(Z)]_d\) for \(1 \leq
  d \leq  |Z| -2\). Then no subset \(W \subseteq Z\) with \(|W| > d+1\) is in linearly
  general position. If \(F\) is irreducible and \(d\) is even this can be improved to say no subset \(W
  \subseteq Z\) with \(|W| > d\) is in linearly general position.
\end{prop}

\begin{proof}  We first proceed in the special case that
   \(|L_p(Z)| \leq d\) for all \(P \in Z\), note that by
   \cref{prop:LpBound} this includes the case that \(F\) is irreducible. If \(W \subseteq Z\) is in linearly general position, then
   for all \(P \in W\) and all \(L \in L_P(W)\), we have that \(|L \cap
   (W \setminus \{P\})| \leq 1\). Therefore
   \( |W| - 1 = |L_P(W)| \leq |L_p(Z)| \leq d\) implying
   \[|W| \leq |L_P(W)| +1 \leq
   d+1.\]

   If furthermore \(d\) is even, then suppose by contradiction that \(W \subseteq Z\)
   is in linearly general position with \(|W| = d+1\). As \(|W| =
   d+1\) we get that \(|L_P(W)| = |L_p(Z)| = d\) for all \(P \in W\).
   Now fix some \(Q \in Z \setminus W\) and define a partition \(\Pi_Q\) of \(W\), where \(P \in W\) is
   contained in the block \(\Spn(Q,P) \cap W\). Now as \(\Spn(Q,P) \in
   L_P(Z) = L_P(W)\), we get \(|\Spn(Q,P) \cap W| = 2\), therefore
   \(\Pi_Q\) is a partition where each block has size \(2\) contradicting
   the fact that \(|W| = d+1\) is odd.

   Now continuing with the general case, let \(F\) be a nonzero possibly
   reducible polynomial. Let \(Z' \subseteq Z\) be the subset \(Z' = \{
   P \in Z \mid \varepsilon_P(F) \neq 0\}\), and let \(T = Z\setminus
   Z'\). Then by \cref{lemma:VanishingEquivalence}, we see that \(F\) factors as \(F
   = G \prod_{P \in T} \ell_P\). Furthermore, \(G \in I^{\gg}(Z')\) and
   \(\varepsilon_P(G) \neq 0\) for all \(P \in Z'\), so by the proof
   of \cref{prop:LpBound} we have \(|L_P(Z')| \leq \deg(G) = d - |T|\) for
   all \(P \in Z'\).  If \(W \subseteq Z\) is in linearly
   general position, then so is \(W' = W \cap Z'\). Applying the
   result from our first
   case we see  
   \[ |W| \leq |W'| + |T| \leq d + 1\]
   establishing the result.
 \end{proof}

 We immediately obtain the following corollary by applying
 \cref{cor:PolynomialStructure}.
 
 \begin{theorem}\label{thm:LinearlySpecial} If \(Z \subseteq
  \bb{P}^2_{\bb{C}}\) admits an unexpected curve in degree \(d\), then
  every subset \(W \subseteq Z\) of points in linearly general
  position has \[|W| \leq d+1.\] Furthermore, we have \(|W| \leq d\)
  if \(d\) is even and the unexpected curve is irreducible.
\end{theorem}

\begin{remark}

  The author suspects the bound of \cref{thm:LinearlySpecial} can be somewhat
  improved over \(\bb{C}\). Namely, given \(Z \subseteq \bb{P}^2_{\bb{C}}\), which
  admits an unexpected curve in degree \(d\), then every subset \(W
  \subseteq Z\) in linearly general position must have \(|W| \leq
  d\). It should be noted, however, that the bound given in
  \cref{prop:IrrLinearlySpecial} is
  sharp in positive characteristic at least if \(d-1\) is a prime power.

  Namely, let \(\bb{K}\) be a field of characteristic \(p > 0\). Let
  \(q = p^e\) and take \(Z = \bb{P}^2_{\bb{F}_{q}} \subseteq
  \bb{P}^2_{\bb{K}}\). Then as shown in \cref{example:FiniteProjectiveSpaces}
  \(Z\) will have an unexpected curve in degree \(q+1\). A
smooth conic such as \(X_1^2 = X_0X_2\) will contain exactly \(q+1\)
points of \(Z\) which form a subset in linearly general position. This
achieves the bound from \cref{prop:IrrLinearlySpecial} if \(q+1\) is even.

If \(q+1\) is odd, then \(\chr(\bb{K}) = 2\) and for each smooth conic \(C \subseteq
\bb{P}^2_{\bb{K}}\) there is a point \(N_C \in \bb{P}^2_{\bb{K}}
\setminus C\) which is contained in every
tangent line of \(C\). This is often referred to as the \emph{nucleus}
of \(C\). As an example, we can verify that \(C :
X_1^2 = X_0X_2\) has nucleus \(N = (0:1:0)\). In this case taking
\(W\) to be \((Z \cap C) \cup \{N\}\), gives a linearly general subset
of size \(q+2\).
\end{remark}

\begin{remark} Proposition \ref{prop:LpBound} can be applied to generalize
  an inductive technique, stated as Lemma 6.5 in \cite{CHMN},
  restricted to the case the bundle  \(\cl{D}_0(\cl{A}_Z)\) is semistable. 
\end{remark}

\begin{prop}\label{prop:UnexpectedInduction} Let \(Z \subseteq
  \bb{P}^2_{\bb{C}}\) and \(P \in \bb{P}^2\) and suppose \(I^{\gg}(Z)\)
 has splitting type \((a,b)\) with \(a \leq b\). If \(L_P(Z) > a\),
 then \(I^{\gg}(Z+P)\)  has splitting type \((a+1,b)\) (or \((a,a+1)\)
 if \(a=b\)).

 In particular, if \(Z\) does not admit unexpected curves, and \(L_P(Z)
 > \left\lceil\dfrac{|Z|}{2}\right\rceil\), then \(Z+P\) does not
 admit unexpected curves.
\end{prop}

\begin{proof}
  First suppose that \([I^{\gg}(Z)]_a = 0\), then by
  \cref{thm:SplittingTypeChar} we have \((a,b) = \left(\left\lfloor
    \frac{|Z|+1}{2} \right\rfloor,\left\lceil
    \frac{|Z|+1}{2} \right\rceil\right)\), and can conclude that
\(|Z|\) admits no unexpected curves. Now for every \(P \in \bb{P}^2
\setminus Z\) we have that \(Z+P\) does not admit unexpected curves.
Since if it did we would have by \cref{thm:UnexpectednessIsGlobal},
that \([I^{\gg}_Q(Z+P)]_d \neq 0\) for some \(d \leq a = \left\lfloor
  \frac{|Z|+1}{2}\right\rfloor\).

  Now suppose that \([I^{\gg}(Z)]_a \neq 0\) and that \(|L_P(Z)| > a\). Note that \(\ell_p F \in [I^{\gg}(Z+P)]_{a+1}\) and it suffices to show that \([I^{\gg}(Z+P)]_a = 0\), since
  then by \cref{thm:SplittingTypeChar} \(I^{\gg}(Z+P)\) has splitting
  type \((\alpha,|Z| - \alpha)\) where \(\alpha =
  \min\left\{a+1,\middle\lfloor\frac{|Z|+2}{2}\middle\rfloor\right\}\).

  For any \(F \in [I^{\gg}(Z+P)]_d\), we have that \(\varepsilon_P(F)
  \in [I(P)^d \cap I(Z)]_d\). Yet as \(|L_P(Z)| > a\) we have by
  applying \cref{lemma:LpEquivalence} that \(\varepsilon_P(F)\) must
  be \(0\). By \cref{lemma:VanishingEquivalence} this means that
  \(\ell_p\) divides \(F\), however then \(F/\ell_P\) is a nonzero
  element of  \([I^{\gg}(Z)]_{a-1}\) giving a contradiction.
\end{proof}

\begin{prop}\label{prop:LineCardBound} Let \(Z \subseteq
  \bb{P}^2_{\bb{C}}\), define
  \[\cl{L} = \{ \Spn(P,Q) \mid P, Q \in Z \text{ are distinct points }\}\]
  and suppose that \(Z\) admits an irreducible
  unexpected curve in degree \(d\), then
  \[|\cl{L}| \leq d^2-d+1\].
\end{prop}

\begin{proof} We prove the theorem under the slightly weaker
  assumption that there is some irreducible \(F_\lambda \in
  [I^{\gg}(Z)]_d\). Without loss of generality assume that \(E_0 =
  (1:0:0) \in Z\), so we may write \(F_\lambda = f_1 X_1 + f_2X_2\) with \(f_i(M_0,M_1,M_2)
  \in \bb{K}[M_0,M_1,M_2]\). Recall the map \(\rho_\lambda:\bb{P}^2
  \to \bb{P}^2\) from \cref{defn:StanleysRho}, in coordinates
  \[\rho_\lambda(a_0,a_1,a_2) =
    (0:f_1(a_0,a_1,a_2):f_2(a_0,a_1,a_2)).\]
  By \cref{thm:GlobalIso} we get
  that \(\rho_\lambda(H) \subseteq H\) for all \(H \in \cl{A}_Z\),
  namely all \(H\) of the form \(H = P^{\perp}\) for \(P \in
  Z\). Define \(\cl{P} = \{ L \cap H \mid L \text{ and } H \text{ are distinct
    lines in } \cl{A}_Z\}\). By projective duality we have that
    \(\cl{L}^{\perp} = \cl{P}\) and so in particular \(|\cl{L}| =
    |\cl{P}|\). Now as \(\rho_\lambda(H) \subseteq H\) for all \(H \in
    \cl{A}_Z\), we get for any \(H \cap L = Q  \in \cl{P}\) that 
    \[\rho_{\lambda}(Q) = \rho_\lambda(H \cap L) \subseteq
      \rho_\lambda(H) \cap \rho_\lambda(L) \subseteq H \cap L = Q.\]
    Hence \(\cl{P}\) is contained in the vanishing locus of
    the minors of
    \[\begin{bmatrix} Y_0 & Y_1 & Y_2 \\ 0 & F_1 & F_2\end{bmatrix}.\]
    So \(\cl{P}\) is contained in the solutions of the polynomial system
      \begin{equation}\begin{aligned}
          Y_1 F_2 - Y_2 F_1 & = 0\\
          Y_0 F_2 &= 0 \\
          Y_0 F_1 &= 0
      \end{aligned}
    \end{equation}
   To count solutions, let \(V\) denote the variety defined by this system, we look at solutions on the line \(Y_0 = 0\)
   and solutions on the subset \(Y_0 \neq 0\). On \(Y_0 = 0\) we get
   that the system (1), reduces to
   \[\begin{aligned}
       Y_1 F_2 - Y_2 F_1 &= 0 \\
       Y_0 &=0 
     \end{aligned}\]
   from which we get by Bezout's theorem that the number of solutions
   is at most \(\deg(Y_0) \deg(Y_1 F_2 - Y_2 F_1) = d\). On the subset
   \(Y_0 \neq 0\), the first equation in the system is redundant and
   the system reduces to
   \[\begin{aligned}
       F_1 &= 0 \\
       F_2 &=0.
     \end{aligned}\]
   As \(F\) is irreducible \(F_1\) and \(F_2\) have no shared
   component so by Bezout's Theorem this system has at most \(\deg(F_1) \deg(F_2) =
   (d-1)^2\) solutions. Combining both results, we can conclude that
   \[|\cl{L}| = |\cl{P}| \leq |V \cap
   (Y_0=0)| + |V \cap (Y_0 \neq 0)| \leq d+(d-1)^2 = d^2-d+1.\]
 \end{proof}

 \begin{example}\label{example:LineCardBound} We note that the above
   bound is sharp in every degree. Namely, the point configuration
   \(C_m \subseteq \bb{P}^2_{\bb{C}}\) of
   \cref{example:FermatTypeArrangements} achieves the bound. To see
   this write \(C_m = \{E_0,E_1,E_2\} \cup F_m\), where \(F_m = \{-E_i
   + \zeta^k E_j \mid 0 \leq i < j \leq 2 \}\). Then the points in \(F_m\) generate
   \(m^2+3\) lines, the \(3\) coordinate lines (\(X_i=0\)), and all lines of the
   form \(\Spn(-E_0+\zeta^j E_1,-E_1+\zeta^k E_2,-E_0+\zeta^{j+k}E_2)\) with defining
   equation \(\zeta^{j+k} X_1 + \zeta^k X_1 +X_2=0\). The only lines
   unaccounted for in \(\cal{L}_{C_m}\) are those of the form
   \(\Spn(E_i,-E_j+\zeta^t E_k)\) with \(\{i,j,k\}=\{0,1,2\}\) of
   which there are \(3m\).

   Then \(|\cal{L}_{C_m}| = m^2+3m+3\), and \(C_m\) admits
   a unique unexpected curve in degree \(m+2\). Noting that \(m^2+3m+3
   = (m+2)^2-(m+2)+1\) we conclude that the above bound is sharp for
   all \(d \geq 4\).
 \end{example}
 
  We close this section with a previously unnoticed combinatorial
  condition which guarantees the existence of unexpected curves.

  \begin{prop}\label{prop:SufficientUnexpectedness} Let \(Z \subseteq \bb{P}^2_{\bb{C}}\) be a finite set of points.
    Further suppose that no line \(L \subseteq \bb{P}^2_{\bb{C}}\) has \(|Z \cap L| \geq
    \frac{|Z|-1}{2}\). Define \(\cl{L}\) as in \cref{prop:LineCardBound}.
    Then for any integer \(d \leq \frac{|Z|+1}{2}\), if
    \[\left(\sum_{p \in Z}|L_P(Z)| \right)-|Z|-|\cl{L}| + 1 < (d-1)(|Z| - d)\]
    then \(|Z|\) admits unexpected curves in degree \(d\).
  \end{prop}

  \begin{proof} Let \(c_t(D_0(\cl{A}_Z))\) denote the Chern polynomial of \(D_0(\cl{A}_Z)\). Then by theorem 2.5 of \cite{schenck},
    \[c_t (D_0(\cl{A}_Z)) = 1-(|Z|-1)t + \left(\sum_{L \in \cl{L}} \left(|L \cap Z| - 1\right) -|Z| + 1\right)t^2.\]
    Noting that
    \[\sum_{L \in Z} |L \cap Z| = \sum_{L \in Z} \sum_{P \in (L \cap Z)} 1 = \sum_{P \in Z} \sum_{L \in L_p(Z)} 1= \sum_{P\in Z} |L_P(Z)|.\]
    We get the following formula for \(c_2(D_0(\cl{A}_Z))\),
    \[\begin{aligned}
        c_2 (D_0(\cl{A}_Z)) &= \sum_{L \in \cl{L}} \left(|L \cap Z| - 1\right) -|Z| + 1\\
        &= \left(\sum_{P \in Z} |L_P(Z)|\right) - |\cl{L}| - |Z| + 1.\\
      \end{aligned}\]
    In particular, we see that our hypothesized inequality is
    equivalent to \(c_2(D_0(\cl{A}_Z)) < (d-1)(|Z| - d)\).

    Now theorem
    B of \cite{BR} states that if \((a,b)\) denotes the splitting type of \(D_0(\cl{A}_Z)\) then
    \(ab \leq c_2 (D_0(\cl{A}_Z))\). So if we satisfy
    \(c_2(D_0(\cl{A}_Z)) < (d-1)(|Z| - d)\), then \(ab < (d-1)(|Z| - d)\).
    Letting \(k = \frac{|Z|-1}{2}\), \(g_1 = k - (d-1)\) and \(g_2 =
    k-a\), then this inequality becomes 
    \[k^2-g_2^2 = (k-g_2)(k+g_2) =ab<(d-1)(|Z|-d)= (k-g_1)(k+g_1) = k^2-g_1^2.\]
    Therefore, we may conclude that \(g_1 < g_2\) and so \(a <
    d-1\). Applying theorem \ref{thm:UnexpectedCurveChar} now establishes the result.
  \end{proof}

\PFEIndent

\section{Regularity Bounds}\label{sec:Regularity}

 In remark 3.8 of \cite{CHMN}, it is claimed that the definition of
 unexpected curves
 
 \say{  \emph{.. leaves open the possibility that the points of
 \(Z\) do not impose independent conditions on curves of some degree
 \(j+1\), and ... a general fat point \(jP\) fails to impose the
 expected number of conditions on the linear system defined by
 \([I_Z]_{j+1}\). Theorem 3.7 gives the surprising result that this is
 impossible.}}

However, it appears to the author that theorem 3.7 of \cite{CHMN} is a
weaker statement than the above quotation claims. Rather it
establishes that \(Z\) imposes independent conditions on a specific
degree \(t_z \geq j+1\), if \(Z\) admits an unexpected curve in degree
\(j+1\).

In this section, we establish the full
claim for points \(Z \subseteq \bb{P}^2_{\bb{C}}\). In fact we prove a
stronger claim. Namely if \(Z \subseteq \bb{P}^2_{\bb{C}}\) 
admits unexpected curves in degree \(d+1\), 
then \(Z\) imposes independent conditions on forms in degree \(d\). This claim is
false in general in positive characteristic, though it does hold for
certain values of \(d\) (see \cref{prop:CharHypothesis}).

Before proceeding we recall the definition of Castelnuovo-Mumford
Regularity, this number determines when \(Z\) imposes independent
conditions on \(d\) forms.

\begin{defn}\label{defn:Regularity} Given a finite set of points \(Z
  \subseteq \bb{P}(V)\) the Castelnuovo-Mumford Regularity of \(Z\),
  denoted \(\reg(Z)\), is the integer
  \[\reg(Z) := 1+\min\{ r \mid \dim_{\bb{K}} [\Sym(V^{*})/I(Z)]_{r} = |Z|\}.\]
\end{defn}

It should be noted that the above definition is highly nonstandard,
and applies only to this specific situation. We refer to Exercise \(4E.3\) and theorem
\(4.2\) of \cite{GSyz}, for proofs that the definition given is
equivalent the standard definitions for graded modules.

This result has some applications to Terao's conjecture as well, which
we explore in the last section.

\begin{prop}\label{prop:CharHypothesis} Let \(\bb{K}\) be an infinite
  field, and let \(A = \bb{K}[s,t]\)
  be a standard graded polynomial ring on \(2\) variables, let
  \(\pow_d:[A]_1 \to [A]_d\) be the \(d\)-th power map, that is the
  map \(\ell \mapsto \ell^d\).

  Then the image of \(\pow_d\) spans \([A]_d\) over \(\bb{K}\) if and
  only if the pair \((\chr \bb{K}, d)\) satisfies one of the following
  \begin{center} \textbf{Characteristic Hypothesis} \end{center}
  \begin{enumerate}
  \item \((\chr \bb{K},d) = (0,d)\); or 
  \item \((\chr \bb{K},d) = (p,q(p^e)-1)\) for some \( e \geq 0\) and
    \(q\), with \( 1 \leq q \leq p\).
  \end{enumerate}
\end{prop} 

\begin{proof} This result is likely well known and consists of
  standard techniques so we only give a brief sketch.

  Let \(L\) be the
  \((d+1) \by (d+1)\) matrix whose \(i\)-th row is \((s+a_it)^d\)
  in  the standard monomial basis of \([A]_d\), also suppose that
  \(a_i \neq a_j\) for \(i \neq j\). Then \(L\) can be seen to
  be a Vandermonde Matrix whose \(j\) column has been scaled by
  \(\binom{d}{j}\). Using the well known Vandermonde Determinant formula the
  matrix is nonsingular and hence the rows span \([A]_d\) if and only
  if \(\prod_{i=0}^d \binom{d}{i}\) is nonzero as an element of \(\bb{K}\). In
  particular, we may conclude if \(\chr (\bb{K}) = 0\).

  If \(\chr(\bb{K}) = p > 0\), we recall Lucas's
  theorem on Binomial coefficients which states \(\displaystyle \binom{d}{i} \not\equiv
  0 \mod p\) if and only if each digit of \(i\) written in base
  \(p\) does not exceed the corresponding digit of \(d\). In base
  \(p\), the only numbers \(d\) where this criterion holds for all
  \(0 \leq i \leq d\) are those \(d\) where the non-leading digits are all
  \(p-1\). This happens precisely when \(d = qp^e-1\) for \(1 \leq q
  \leq p\).
\end{proof}

\begin{prop}\label{prop:RegularityBound} Let \(Z \subseteq \bb{P}^2_{\bb{K}}\) be a finite set of
  points, suppose that \(\dim \Bloc_{d+1}(I^{\gg}(Z)) = 0\) 
  and that the pair \((\chr \bb{K},d)\) satisfies the
  characteristic hypothesis of \cref{prop:CharHypothesis}. Then
  \(\reg(Z) \leq d+1\).
\end{prop}

\begin{proof} Let \(\bb{P}^2 = \Proj(R)\), where \(R =
  \bb{K}[X_0,X_1,X_2]\). Take \(\ell \in [R]_1\) to be a general
  linear form and consider the short exact sequence
  \[\xymatrix{0 \ar[r] & [R/I(Z)]_{t-1} \ar[r]^{\ell} & [R/I(Z)]_{t}
      \ar[r] & [R/(I(Z)+(\ell))]_t \ar[r] & 0}\]

  Letting \(h_Z(t) := \dim [R/I(Z)]_t\) we conclude that for integers
  \(t\) that
  \[h_z(t) - h_z(t-1) = \dim [R/(I(Z) +  \ell)]_t.\] Furthermore, as \(R/(I(Z)+\ell)\) is principally
  generated we can conclude that \(h_z(t) = h_z(t-1)\) if and only if
  \(h_z(t-1) = |Z|\). From this it follows from 
  \cref{defn:Regularity} that \(\reg(Z) = \min\{ r \mid [R/(I(Z) +
  (\ell))]_r = 0\}\) and it suffices to prove that
  \([R/(I(Z)+(\ell))]_{d+1} = 0\).

  Fix \(F_\lambda \in [I^{\gg}(Z)]_{d+1}\) with \(\dim
  \Bloc(F_\lambda) = 0\). For all points \(Q\) on the line \(\ell =0\), we have by
  \cref{lem:HypFactorization} that
  \[\varepsilon_Q(F_\lambda) = h_Q^{d}\rho_\lambda(\ell) \mod (\ell) \]
  for some linear form \(h_Q\) vanishing on \(Q\).  Noting
  \(\varepsilon_Q(F_\lambda) \in I(Z)\), we get an inclusion
  of \(\bb{K}\)-vector spaces
  \[[(I(Z) + (\ell))/(\ell)]_{d+1} \supseteq \Spn\{h_Q^{d}
    \rho_\lambda(\ell) +(\ell) \mid Q \in \bb{P}^2 \text{ and } \ell(Q) = 0\}.\]
  By \cref{prop:CharHypothesis} the set \(\{h_Q^d \mid Q
    \in (\ell = 0)\}\) spans \([R/(\ell)]_d\) and so
   \[[I(Z) + (\ell)/(\ell)]_{d+1} \supseteq \rho_\lambda(\ell)  [R/(\ell)]_{d}. \]
   Hence, \([I(Z) + (\ell)]_{d+1} \supseteq
   [(\ell,\rho_\lambda(\ell))]_{d+1}\). Let \(P \in \bb{P}^2\) denote
   the point defined by the ideal \((\ell,\rho_\lambda(\ell))\), as
   \(\Bloc(F_\lambda) \cap (\ell = 0) = \emptyset\) we can find some \(H \in \bb{P}^2\) so
   that \(\varepsilon_H(F_\lambda) \in I(Z)\) and
   \(\varepsilon_H(F_\lambda) \not\in I(P) =
   (\ell,\rho_\lambda(\ell))\). Hence,
   \[[I(Z) + (\ell)]_{d+1} \supseteq
   [(\ell,\rho_\lambda(\ell))]_{d+1} +
   [\varepsilon_H(F_\lambda)]_{d+1} = [R]_{d+1}\]
 allowing us to conclude that \([R/(I(Z) + (\ell))]_{d+1} = 0\) as desired.
\end{proof}

\begin{remark}It should be noted that the assumptions of the above
  theorem, may be relaxed in various ways to give slightly different
  bounds, which also require different proofs, and possibly stronger
  assumptions. We have choosen to give 
  only the proof above for the sake of brevity, but will briefly
  comment on two of the possible changes now.

  \begin{enumerate}
  \item The condition \(0 = \dim \Bloc_{d+1}(I^{\gg}(Z)) \) may be
    replaced with the weaker condition that \(0 = \dim \Bloc_{d+1} (Z)\). However the bound then becomes
  \(\reg(Z) \leq d+2\). The proof is similar to above, but \(\ell\) is replaced
  with a line through a general point \(Q\), which also vanishes on a
  point in \(Z\). It can then only be concluded that \(\dim
  [R/(I(Z)+(\ell))]_{d+1} \leq 1\). This worse bound of \(\reg(Z) \leq
  d+2\) is in fact sharp, as can be illustrated by taking \(Z\) to be \(2d+3\) points on
  a smooth conic.

  Interestingly, this technique can be used to give a completely geometric proof of
  \cref{thm:SegreBound} for points in \(\bb{P}^2_{\bb{C}}\), in contrast to the
  combinatorial proof given in \cite{NT}.

  \item A generalization to \(\bb{P}^n\), at least in characteristic
    \(0\), is possible however the
    proof becomes more involved and/or additional assumptions on
    \([I^{\gg}(Z)]_{d+1}\) are necessary. The main technique is still
    roughly the same except now induction is needed. After
    showing \([I(Z)+(\ell)]_{d+1} \supseteq
    (\ell,\rho_\lambda(\ell))\) one proceeds as before showing
    \[[I(Z) + (\ell)]_{d+1} = [I(Z)+(\ell,\rho_\lambda(\ell))]_{d+1}
    \supseteq
    [(\ell,\rho_\lambda(\ell),\rho_\lambda^2(\ell))]_{d+1}\supseteq
    ...\]
  If \((\ell,\rho_\lambda(\ell),...,\rho_\lambda^{n-1}(\ell))\) is the
  ideal of a point we then proceed as in the proposition.
  \end{enumerate}
\end{remark}

Combining the above from some results from earlier sections, we may
obtain the following result
\begin{theorem}\label{thm:UnexpecteCurveRegularity} If \(Z \subseteq \bb{P}^2_{\bb{C}}\) has an unexpected
  curve in degree \(d\), then \(\reg(Z) \leq d\). In particular, \(Z\)
  imposes independent conditions on forms of degree \(d-1\).
\end{theorem}

\begin{proof}
 Suppose that \(Z\) admits an unexpected curve in degree \(d\), without
  loss of generality we assume that \(Z\) does not admit an unexpected
  curve in degree \(d-1\). Then by theorem
  \ref{thm:UnexpectedCurveChar}, we note that \(|Z| \geq 2d+1\) and
  that no line \(L\) contains more that \(d+1\) points of \(Z\).
  Additionally by theorem \ref{thm:UnexpectednessIsGlobal},
  there exists \(F \in [I^{\gg}(Z)]_d\) defining the curve. We claim
  that \(\dim \Bloc(F) = 0\), which in light of
  \cref{prop:RegularityBound} establishes the claim.

  Proceeding by contradiction assume that  \(\dim \Bloc(F) =1\), applying \cref{prop:BlocIsLinear}, we
  get \(\Bloc(F)\) has a component which is a line. If \(\ell \in \bb{C}[X_0,X_1,X_2]\) is the
  linear form defining this line, \(L\), then viewing it as a polynomial in
  \(\bb{C}[X_0,X_1,X_2,A_0,A_1,A_2]\) and applying lemma
  \ref{lemma:DivisibilityLemma}, we get \(F = \ell h\) with \(h \in
  \bb{C}[M_0,M_1,M_2]\). Then as \(\varepsilon_Q(h) \in [I(Q)^{\deg
    h}]_h\) for all \(Q \in \bb{P}^2\) we get that the variety
  \(V(\varepsilon_Q(F))\) is a union of \(L\) and at most \(\deg h\)
  lines through \(Q\). For a general \(Q \in
  \bb{P}^2_{\bb{C}}\), each line in \(V(\varepsilon_Q(h))\)
  contains at most one point of \(Z\). This forces us to conclude that
  \(|L \cap Z| \geq |Z| - \deg h  \geq d+2\) giving us our desired contradiction.
\end{proof}

We note that the \cref{thm:UnexpecteCurveRegularity}, gives a decent
criteria purely algebraic criteria for establishing that a set of points \(Z\) does not admit
unexpected curves in a given degree. We explore this a bit in the next
section in the context of Terao's Conjecture. 

\PFEIndent
\section{Applications to Terao's Conjecture in
  \texorpdfstring{\(\bb{P}^2\)}{the Plane}}\label{sec:ApplicationsToTeraos}

A much studied problem in the theory of line arrangements is the
freeness of the module of derivations \(D(\cl{A}_Z)\). One reason for
this in particular is that if \(D(\cl{A}_Z)\) is free then many of the invariants of
\(D(\cl{A}_Z)\) can be determined from combinatorics of the
intersection lattice \(L(\cl{A}_Z)\) (or equivalently the
matroid \(M(Z)\)). A major open problem in the study of Hyperplane
arrangements is Terao's Freeness Conjecture

\begin{conj}[Terao's Freeness Conjecture]\label{conj:Teraos} Over
  \(\bb{C}\) freeness of \(D_0(\cl{A})\) can be determined by the
  intersection lattice \(L(\cl{A}_Z)\).
\end{conj}

\begin{remark} The above conjecture is usually stated for
  \(D(\cl{A})\). However the two versions are equivalent because over \(\bb{C}\), \(D(\cl{A})\)
  splits as \((\Sym(V^{*}))\theta_e \oplus D_0(\cl{A})\).
\end{remark}

One natural question to ask given theorem \ref{thm:GlobalIso}, is what
freeness of \(D_0(\cl{A}_Z)\) says about \(I(Z)\). Namely, can
\(D_0(\cl{A}_Z)\) be characterized in terms of \(Z\)?
The following proposition (which is well known to experts) is helpful
in addressing this question.

\begin{prop}\label{prop:FreeBySurjection} \(D_0(\cl{A})\) is free if and only if for a general line
  \(L \subseteq \bb{P}(W)\), the restriction map \(D_0(\cl{A}) \to
  D_0(\cl{A})\mid_L\) is surjective.
\end{prop}

\begin{proof} The forward implication is clear. For the reverse
  implication, we apply a corollary of Saito's Criterion which can be
  found as theorem 4.23 of \cite{ot1992}. Namely \(D_0(\cl{A})\) is
  free if there exists \(\theta_1,..,\theta_n \in D_0(\cl{A}_Z)\)
  which are linearly independent over the
  projective coordinate ring, \(S\), of \(\bb{P}(W)\), and where
  \(\sum_{i=1}^n \deg(\theta_i) = |\cl{A}| - 1\).

  So suppose that \(\res_L: D_0(\cl{A}) \to D_0(\cl{A}_Z)\mid_L\) is
  surjective for a general line \(L\). Let
  \(\bar{\theta_1},...,\bar{\theta_n}\) be a \(S/I(L)\)-basis of
  \(D_0(\cl{A}_Z)\mid_L\), then for each \(i\) we can find \(\theta_i
  \in D_0(\cl{A})\) so that \(\res_L(\theta_i) =
  \bar{\theta_i}\). As \(\sum_{i=1}^n \deg (\theta_i) =
  \sum_{i=1}^n \deg (\bar{\theta_i}) = |\cl{A}| - 1\) it suffices to
  show that the \(\theta_i\) are linearly independent over \(S\). Yet
  if \(\sum_{i=1}^n s_i \theta_i = 0\) for some \(s_i \in S\) and some index \(j\), then  \(\sum_{i=1}^n
  \bar{s_i} \bar{\theta_i} = 0\) in \(D_0(\cl{A}_z)\mid_L\). As \(L\)
  is general if \(s_j \neq 0\) for some index \(j\) then we can assume that \(s_j \not\in
  I(L)\) which gives a non-trivial relation among
  \(\bar{\theta_1},..,\bar{\theta_n}\) and a contradiction.
\end{proof}

\begin{cor}\(\cl{A}_Z\) is a free arrangement if and only if the
  evaluation map \(\varepsilon_Q:I^{\gg}(Z) \to I^{\gg}_Q(Z)\) is
  surjective for general \(Q\).
\end{cor}

Theorem \ref{thm:UnexpectednessIsGlobal} has some applications to
Terao's conjecture, namely we give a new criterion for determining
Freeness.

\begin{prop}\label{prop:2ndGeneratorCriterion} Let \(\cl{A} \subseteq \bb{P}^2_\bb{C} = \Proj(S)\) with splitting
  type \((a_1,a_2)\). If \(a_2 - a_1 \geq 2\), then  \(D_0(\cl{A})\) is
  free if and only if it has a minimal generator in degree \(a_2\). 
\end{prop}

\begin{proof} We use the criterion from \cref{prop:FreeBySurjection}. As  \(a_2 - a_1 > 2\), then we may apply
  \cref{thm:UnexpectednessIsGlobal} to see 
  there nonzero \(\theta_1 \in [D_0(\cl{A})]_{a_1}\) and so the image of the restriction map contains the generator of
  \(D_0(\cl{A})\mid_L\) in degree \(a_1\). If \(D_0(\cl{A})\) has a
  minimal generator \(\theta_2\) in degree \(a_2\), then \(\theta_2
  \neq f \theta_1\) for any \(f \in S\). For a general line \(L\), we still
  have \(\res_L(\theta_2)  \not \in S/I(L) \res_L(\theta_1)\) and
  conclude that \(\{\res_L(\theta_1),\res_L(\theta_2)\}\) is a
  generating set for \(D_0(\cl{A})\mid_L\).
\end{proof}

Additionally, it is well known that freeness can be
determined from combinatorics and the splitting type. More precisely,

\begin{prop}\label{thm:TeraoAsSplittingType} Let \(\cl{A}\) and \(\cl{B}\) be hyperplane
  arrangements in \(\bb{P}^n\), and suppose \(\cl{A}\) and \(\cl{B}\)
  have isomorphic intersection  lattices. Suppose that \(D_0(\cl{A})\) is free, then \(D_0(\cl{B})\)
  is free if and only if it has the same splitting type as \(\cl{A}\).
\end{prop}

\begin{proof} By a theorem of Terao \(c_2(D_0(\cl{A}))\) is determined
  solely by \(L_{\cl{A}}\). The result is now a consequence of the criterion that
  \(D_0(\cl{A})\) is free if and only if \(c_2(D_0(\cl{A})) = a_1
  a_2\)  where \((a_1,a_2)\) is the splitting type, see for instance \cite{BR}. 
\end{proof}

This characterization allows us to generalize a theorem of
\cite{schenck} 
which was stated only for balanced free arrangements in
\(\bb{P}^2_{\bb{C}}\). Here balanced means free
arrangements with splitting type \((a,a)\) or \((a,a+1)\).

\begin{theorem}  For a finitely generated graded module
  \(M\), let \(\alpha(M)\) denote the initial degree of \(M\), that is
\[\alpha(M) := \inf\{d \in \bb{Z} \mid [M]_d \neq 0 \}\]

  Let \(\cl{A}\) and \(\cl{B}\) be combinatorially equivalent line arrangements
  \(\bb{P}^2_{\bb{K}}\). If
  \(D_0(\cl{A})\) is free, then
  \[\alpha( D_0(\cl{B})) \leq \alpha( D_0(\cl{A}))\]
  with equality if and only if \(\cl{B}\) is free.

  In particular, if \(\cl{A}\) is free with exponents \((1,a,b)\) and
  \(\cl{B}\) is not free, then \(D_0(\cl{B})\) has a generator in
  degree \(<a\) and all other minimal generators are in degree \( > b\).
\end{theorem}

\begin{remark} Note if \(\bb{K} \subseteq \bb{C}\), the following argument can
  be slightly simplified by applying \cref{thm:UnexpectednessIsGlobal}.
\end{remark}

\begin{proof}
  The reverse direction is immediate as in that case \(D_0(\cl{A})\)
  and \(D_0(\cl{B})\) are isomorphic.

  To prove the forward implication, we apply the characterization given in
  \ref{thm:TeraoAsSplittingType}. Hence, assume that \(D_0(\cl{A})\) has splitting
  type \((a_1,a_2)\), and \(D_0(\cl{B})\) has splitting type
  \((b_0,b_1)\), with \((b_0,b_1) \neq (a_0,a_1)\). It suffices to
  show that \(\alpha(D_0(\cl{B})) < \alpha(D_0(\cl{A}))\).

  By \cite{Yuz} freeness is an open
  property. Hence, if \(\cl{B}\) is not free it lies on closed subvariety
  of \(V_{L(\cl{A})}\) the variety parameterizing arrangements with
  intersection lattice isomorphic to \(L(\cl{A})\). We can view \(D(\cl{B})\)
  as the kernel of the linear map \(\Der(S) \to \prod_{H
    \in \cl{B}} S/(I(H))\) which maps \(\theta \mapsto
  (\theta(\ell_H))_{\{H \in \cl{B}\}}\), so by lower semicontinuity of rank
  we may conclude that \(\dim [D_0(\cl{B})]_d \geq \dim  [D_0(\cl{A})]_d\) for all \(d\). Applying the same argument 
  to the restriction of \(D_0(\cl{B})_d\) to a general line, we see that \(b_0 < a_0
  \leq a_1 < b_1\). As \(\dim [D_0(\cl{B})]_{a_0} \geq \dim
  [D_0(\cl{A})]_{a_0} > 0 \), we can apply
  \cref{prop:DivisibilityProp} to get
  
  \[\alpha(D_0(\cl{B})) = b_0 < a_0 =
    \alpha(D_0(\cl{A})).\]

  The final sentence follows from this and \cref{prop:2ndGeneratorCriterion}.
\end{proof}

One corollary of the above theorem is an extension of a theorem of
\cite{FV} over \(\bb{C}\), to positive characteristic.

\begin{cor}If \(\cl{A} \subseteq \bb{P}^2_{\bb{K}}\) is a free
  arrangement with splitting type \((a_1,a_2)\) and some point \(P \in
  \bb{P}^2_{\bb{K}}\) is incident to at least \(a_1\) lines of
  \(\cl{A}\). Then any arrangement over \(\bb{P}^2_{\bb{K}}\)
  combinatorially equivalent to \(\cl{A}\) is also free.
\end{cor}

\begin{proof} Let \(\cl{B}\) be an arrangement that is combinatorially
  equivalent to \(\cl{A}\), and let \((b_1,b_2)\) denote it's
  splitting type. Dualizing the problem statement, we see that there
  is a line \(H\) containing at least \(a_1\) points of
  \(\cl{B}^{\perp}\) so by \cref{thm:UnexpectedCurveChar} and the
  prior theorem we have
  \(a_1-1 \leq b_1 \leq a_1\).

  Let \(h\) denote the linear form defining this line.
  Furthermore for a general \(Q\), let \(g\) denote the product of
  linear forms through \(Q\) and each point not on \(H\), then \(h g
  \in [I^{\gg}(\cl{B}^{\perp})]_{a_2+2}\). If \(b_1 = a_1-1\), then
   \(hg\) would correspond by \cref{thm:GlobalIso}
  to a minimal generator of \(D_0(\cl{B})\) in degree
  \(b_2=a_2+1\). The prior theorem together with
  \cref{thm:TeraoAsSplittingType} now gives a contradiction.
\end{proof}

We now close by discussing connections between Terao's conjecture and
a conjecture due to Dirac. It was conjectured in 
\cite{Dirac} that for every finite set of points non collinear points
\(Z \subseteq \bb{P}^2_{\bb{R}}\), that there is always some \(Q \in Z\) so
that
\[|L_Q(Z)| = \left|\left\{\Spn(Q,P) \mid P \in Z \setminus
        Q\right\}\right| \geq \frac{|Z|}{2}.\]

  However, some counterexamples have been found to the original
  formulation (see \cite{Grun}). This has lead to two reformulations of the original
  conjecture which we reprint below.

\begin{conj}[Weak Dirac Problem] Determine the smallest constant \(C\), so that for every
  finite set of noncollinear points \(Z \subseteq \bb{P}^2_{\bb{R}}\),
  there exists some \(Q \in Z\) where
  \[|L_Q(Z)|  \geq  \frac{|Z|}{C} \]
\end{conj}

\begin{conj}[Strong Dirac Conjecture] There exists some constant
  \(c_0>0\) so that for every set of finite noncollinear points \(Z
  \subseteq \bb{R}^2\), there exists some \(Q \in Z\) so that 
  \[|L_Q(Z)| \geq \frac{Z}{2} -c_0\]
\end{conj}

Counterexamples have been found to Dirac's Original Conjecture for
every odd \(n = |Z|\) with the exception of those \(n\) of the form
\(n = 12k + 11\) with \(k \geq 4\) (see \cite{AIKN}). Despite that the known
counterexamples only barely break the original conjecture bound. Most
satisfy the Strong Dirac Conjecture with \(c_0 = 1/2\) and all but
finitely many satisfy the conjecture with \(c_0 = 3/2\).

We now show that any minimal counterexample to Terao's Conjecture
for real line arrangements must itself be a counterexample to the
original Conjecture of G.Dirac, and must be extremal in the regards to
the other two cases.

\begin{theorem}\label{thm:TeraoToDirac} Let \(\cl{A}\) and \(\cl{B} \subseteq \bb{P}^2\) be
  real (or complex) line arrangements, which form a counter example to
  Terao's conjecture. Meaning \(L_{\cl{A}} \iso L_{\cl{B}}\), but
  \(D_0(\cl{A})\) is free with splitting type \((a_1,a_2)\) where as
  \(D_0(\cl{B})\) is not free. Furthermore, suppose
  there is no pair of lines \((L,L') \in \cl{A} \times \cl{B}\) we can remove to get
  subarrangements \(\cl{A'} =\cl{A} \setminus \{L\}\) and \(\cl{B'} =
  \cl{B} \setminus \{L'\}\) forming a smaller
  counterexample.

  Then letting \(\cl{A}^{\perp}\) be the set of points dual to
  \(\cl{A}\) we have
  \[|L_P(Z)| \leq a_1 \leq \left\lfloor\frac{ |Z| - 1}{2}
    \right\rfloor.\]
\end{theorem}

Our proof of the above theorem relies on the following proposition 
which seems useful in it's own right.
It is related to Terao's well known Addition-Deletion Theorem

\begin{prop} Let \(\cl{A}_z \subseteq \bb{P}^2_{\bb{C}}\) be a free line
  arrangement with splitting type \((a_1,a_2)\)  and \(Z\) the dual
  set of points. If there is some \(P \in Z\)
  with \(|L_P(Z)| > a_1 + 1\), then \(|L_P(Z)| = a_2+1\) and
  \(\cl{A}_{W}\) is free where \(W = Z \setminus \{P\}\).
\end{prop}

\begin{proof} By theorem B of \cite{BR}, \(c_2(D_0(\cl{A}_Z)) \geq a_1 a_2\), and
  \(\cl{A}\) is free if and only if equality holds. Furthermore, if
  \(L_P(Z) > a_1 + 1\), then letting \(F \in [I^{\gg}(Z)]\) it follows
  by \cref{prop:LpBound} that \(\varepsilon_P(F) = 0\). Then
  by \cref{lemma:VanishingEquivalence} the
  linear form, \(\ell_p\), defining the line dual to \(P\) must divide
  \(F\). However, then we necessarily have \(F/\ell_P \in
  [I^{\gg}(W)]_{a_1} \iso [D_0(\cl{A})]_{a_1-1}\) so \(\cl{A}_W\) must have
  splitting type \((a_1-1,a_2)\).

  We note that it suffices to show that \(D_0(\cl{A}_W)\) is free,
  since Terao's Famous Addition-Deletion Formula then ensures that
  \(L_P(W) = a_2+1\). Yet this follows
  since if \(F\) and \(G\) freely generate \(D_0(\cl{A}_Z)\), then \(F
  / \ell_P\) and \(G\) must generate \(D_0(\cl{A}_W)\).
\end{proof}

\begin{proof}[Proof of \cref{thm:TeraoToDirac}] By the preceding
  proposition there exists no \(P \in Z\) with \(|L_P(Z)| > a_1+1\).
  Furthermore, by Terao's Addition-Deletion formula there is no \(P
  \in L_P(Z)\) with \(|L_p(Z)| = a_1+1\), since then letting \(\cl{A'}
  = \cl{A} \setminus \{\ell_0 = 0\}\) we would get a smaller
  counterexample. Hence, for all \(P \in Z\) we have
  \[|L_p(Z)| \leq a_1 \leq \frac{|Z|-1}{2}\]
  establishing the result.
\end{proof}

\PFEIndent

\bibliographystyle{plain}
\bibliography{Final_Cut}

\begin{thebibliography}{10}

\bibitem{AIKN}
Jin Akiyama, Hiro Ito, Midori Kobayashi, and Gisaku Nakamura.
\newblock Arrangements of {$n$} points whose incident-line-numbers are at most
  {$n/2$}.
\newblock {\em Graphs Combin.}, 27(3):321--326, 2011.

\bibitem{AH}
J.~Alexander and A.~Hirschowitz.
\newblock Polynomial interpolation in several variables.
\newblock {\em J. Algebraic Geom.}, 4(2):201--222, 1995.

\bibitem{STK}
The Stacks~Project Authors.
\newblock {\em The Stacks Project}, 2020.

\bibitem{BMSS}
Thomas Bauer, Grzegorz Malara, Tomasz Szemberg, and Justyna Szpond.
\newblock Quartic unexpected curves and surfaces.
\newblock {\em manuscripta mathematica}, 11 2018.

\bibitem{BR}
Cristina Bertone and Margherita Roggero.
\newblock Splitting type, global sections and {C}hern classes for torsion free
  sheaves on {${\bf P}^n$}.
\newblock {\em J. Korean Math. Soc.}, 47(6):1147--1165, 2010.

\bibitem{Boj.}
R.~Bojanowski.
\newblock Zastosowania uog{\'o}lnionej nier{\'o}wno{\'s}ci
  bogomolova-miyaoka-yau. master thesis (in polish).
\newblock 2003.

\bibitem{CHMN}
D.~Cook, II, B.~Harbourne, J.~Migliore, and U.~Nagel.
\newblock Line arrangements and configurations of points with an unexpected
  geometric property.
\newblock {\em Compos. Math.}, 154(10):2150--2194, 2018.

\bibitem{DIV}
Roberta Di~Gennaro, Giovanna Ilardi, and Jean Vall\`es.
\newblock Singular hypersurfaces characterizing the {L}efschetz properties.
\newblock {\em J. Lond. Math. Soc. (2)}, 89(1):194--212, 2014.

\bibitem{Dirac}
G.~A. Dirac.
\newblock Collinearity properties of sets of points.
\newblock {\em Quart. J. Math. Oxford Ser. (2)}, 2:221--227, 1951.

\bibitem{DFHMST}
Marcin Dumnicki, Lucja Farnik, Brian Harbourne, Grzegorz Malara, Justyna
  Szpond, and Halszka Tutaj-Gasinska.
\newblock A matrixwise approach to unexpected hypersurfaces.
\newblock 2019.

\bibitem{Eureka}
J.~Edmonds.
\newblock {\em Combinatorial Optimization -- Eureka, You Shrink}, volume 2570
  of {\em Lecture Notes in Computer Science}, chapter Submodular Functions and
  Certain Polyhedra, pages 11--26.
\newblock Springer-Verlag Berlin Heidelberg.

\bibitem{GSyz}
D.~Eisenbud.
\newblock {\em The Geometry of Syzygies: a Second Course in Algebraic Geometry
  and Commutative Algebra}, volume 229 of {\em Graduate Texts in Mathematics}.
\newblock Springer - New York, 2005.

\bibitem{FV}
Daniele Faenzi and Jean Vall\`es.
\newblock Logarithmic bundles and line arrangements, an approach via the
  standard construction.
\newblock {\em J. Lond. Math. Soc. (2)}, 90(3):675--694, 2014.

\bibitem{Grun}
B~{Gr\"{u}nbaum}.
\newblock {\em Arrangments and Spreads}, volume~10 of {\em CBMS Regional
  Conference Series in Mathematics}.
\newblock American Mathematical Society, 1972.

\bibitem{Han}
Zeye Han.
\newblock A note on the weak {D}irac conjecture.
\newblock {\em Electron. J. Combin.}, 24(1):Paper 1.63, 5, 2017.

\bibitem{HMUZ}
B.~Harbourne, J.~Migliore, U.~Nagel, and Z.~Teitler.
\newblock Unexpected hypersurfaces and where to find them.
\newblock {\em Michigan Math. J.}, page (to appear), 2018.

\bibitem{HMT}
B.~Harbourne, J.~Migliore, and H.~Tutaj-Gasińska.
\newblock New constructions of unexpected hypersurfaces in \(\mathbb{P}^n\).
\newblock 2019.

\bibitem{Hart}
Robin Hartshorne.
\newblock Stable reflexive sheaves.
\newblock {\em Math. Ann.}, 254(2):121--176, 1980.

\bibitem{Hoch}
Melvin Hochster.
\newblock Criteria for equality of ordinary and symbolic powers of primes.
\newblock {\em Math. Z.}, 133:53--65, 1973.

\bibitem{Lang.}
Adrian Langer.
\newblock Logarithmic orbifold {E}uler numbers of surfaces with applications.
\newblock {\em Proc. London Math. Soc. (3)}, 86(2):358--396, 2003.

\bibitem{Migl08}
Juan~C. Migliore.
\newblock The geometry of the weak {L}efschetz property and level sets of
  points.
\newblock {\em Canad. J. Math.}, 60(2):391--411, 2008.

\bibitem{NT}
U.~{Nagel} and B.~{Trok}.
\newblock {Segre's Regularity Bound for Fat Point Schemes}.
\newblock {\em Ann. Sc. Norm. Super. Pisa Cl. Sci.}, page (to appear), March
  2020.

\bibitem{Nara}
H.~Narayanan.
\newblock The principal lattice of partitions of a submodular function.
\newblock {\em Linear Algebra Appl.}, 144:179--216, 1991.

\bibitem{ot1992}
Peter Orlik and Hiroaki Terao.
\newblock {\em Arrangements of hyperplanes}.
\newblock Springer-Verlag, 1992.

\bibitem{Oxley}
J.~Oxley.
\newblock {\em Matroid Theory}, volume~21 of {\em Oxford Graduate Texts in
  Mathematics}.
\newblock Oxford University Press, 1992.

\bibitem{Pok}
P.~{Pokora}.
\newblock Hirzebruch-type inequalities viewed as tools in combinatorics.
\newblock {\em ArXiv e-prints}, 2018.

\bibitem{schenck}
Henry~K. Schenck.
\newblock Elementary modifications and line configurations in {$\mathbb{P}^2$}.
\newblock {\em Comment. Math. Helv.}, 78(3):447--462, 2003.

\bibitem{Yuz}
Sergey Yuzvinsky.
\newblock Free and locally free arrangements with a given intersection lattice.
\newblock {\em Proc. Amer. Math. Soc.}, 118(3):745--752, 1993.

\end{thebibliography}

\end{document}